\documentclass[a4paper,11pt,leqno,notitlepage,twoside]{amsart}
  
\usepackage[T1]{fontenc}
\usepackage[latin1]{inputenc}
\usepackage{amsmath}
\usepackage{amssymb}
\usepackage{latexsym} 
\usepackage{epsf}
\usepackage{amsthm}
\usepackage{amsfonts}
\usepackage{graphicx}
\usepackage{url}
\usepackage[right]{eurosym}
\usepackage{ifpdf}
\usepackage{array}
\usepackage{marvosym}
\usepackage{wasysym}

\numberwithin{equation}{section}
  
\theoremstyle{plain} 
\newtheorem{theorem}{Theorem}[section] 
 
\newtheorem{proposition}[theorem]{Proposition}
\newtheorem{lemma}[theorem]{Lemma}

\newtheorem*{theo}{Theorem}
\newtheorem*{theoA}{Theorem A}
\newtheorem*{theoB}{Theorem B}

\newtheorem*{propoC}{Proposition C}
\newtheorem*{propoD}{Proposition D}

\theoremstyle{definition} 
\newtheorem{remark}{Remark}[section]
 
\newtheorem*{hypothesis}{\textbf{Hypothesis} $\mathbf{H_2(\theta)}$}

\title[Real zeros and size of Rankin-Selberg $L$-functions in the level aspect]{Real zeros and size of Rankin-Selberg $L$-functions in the level aspect} 
\author[G. Ricotta]{G. Ricotta}
\thanks{ricotta@math.univ-montp2.fr \newline
\vrule width 12pt height 0pt depth 0pt 2000 \textit{Mathematics Subject Classification}. Primary 11M41}
\date{Version of \today}

\begin{document}

\begin{abstract} 
In this paper, some asymptotic formulas are proved for the harmonic mollified second moment of a family of Rankin-Selberg $L$-functions. One of the main new input is a substantial improvement of the admissible length of the mollifier which is done by solving a shifted convolution problem by a spectral method on average. A first consequence is a new subconvexity bound for Rankin-Selberg $L$-functions in the level aspect. Moreover, infinitely many Rankin-Selberg $L$-functions having at most eight non-trivial real zeros are produced and some new non-trivial estimates for the analytic rank of the family studied are obtained.
\end{abstract} 
 
\maketitle

\tableofcontents 

\section{Introduction and statement of the results}

This paper is motivated by the striking result of J.B. Conrey and K. Soundararajan proven in \cite{CoSo}:
\begin{theo}[J.B. Conrey-K. Soundararajan (2002)]
\label{zeros}
There exists infinitely many (at least $20\%$ in a suitable sense) primitive quadratic Dirichlet characters $\chi$ whose Dirichlet L-function $L(\chi,s):=\sum_{n}\chi(n)n^{-s}$ does not vanish on the critical segment $\left[0,1\right]$.
\end{theo}
The family of L-functions considered in \textbf{\cite{CoSo}} is $\mathcal{G}:=\cup_{X\in\left\{2^m, m\in\mathbb{N}\right\}}\mathcal{G}(X)$ with
\begin{equation*}
\mathcal{G}(X):=\left\{L\left(\chi_{-8d},.\right), 2\nmid d, \mu^2(d)=1, X\leq d\leq 2X\right\}
\end{equation*}
where $\chi_{-8d}(n):=\left(\frac{-8d}{n}\right)$ is the Kronecker symbol. The proof, which is based on the mollification method, exploits the following properties of the family $\mathcal{G}$:
\begin{itemize}
\item
the functional equation of each $L$-function of this family has the same sign;
\item
this sign equals $+1$ and consequently the order of vanishing at the critical point $\frac{1}{2}$ of each $L$-function is an even integer;
\item
the symmetry type of this family is symplectic - this entails that the first zero is repelled from the real axis and justifies the method used by the authors.
\end{itemize}
K. Soundararajan announced at the Journées Arithmétiques $2003$ in Graz a similar result for the families $\mathcal{H}_\pm:=\cup_{K\in\left\{2^m, m\in\mathbb{N}^*\right\}}\mathcal{H}_\pm(K)$ with
\begin{eqnarray*}
\mathcal{H}_+(K) & := & \left\{L\left(f,.\right), f\in S_k^p(1), K\leq k\leq 2K, k\equiv 0\mod 4\right\}, \\
\mathcal{H}_-(K) & := & \left\{L\left(f,.\right), f\in S_k^p(1), K\leq k\leq 2K, k\equiv 2\mod 4\right\}
\end{eqnarray*}
where $S_k^p(1)$ denotes the set of primitive cusp forms of level $1$, weight $k$ and trivial nebentypus. It is then natural to try to generalize these results to other families of L-functions. Throughout this article, $g$ will be a \textbf{fixed} primitive (arithmetically normalized namely with first Fourier coefficient equal to one) cusp form of square-free level $D$, weight $k_{g}$ and trivial nebentypus $\varepsilon_{D}$, and $f$ will be a \textbf{varying} primitive cusp form of level $q$, weight $k$ and trivial nebentypus $\varepsilon_{q}$ denoted by $f\in S_k^p(q)$. We prove a result cognate to that of \textbf{\cite{CoSo}} for the family of Rankin-Selberg $L$-functions $\mathcal{F}:=\cup_{\substack{q\in\mathcal{P} \\ q\nmid D}}\mathcal{F}(q)$ where:
\begin{equation*}
\forall q\in\mathcal{P}, \quad\mathcal{F}(q):=\left\{L(f\times g,.),f\in S_{k}^{p}(q)\right\}.
\end{equation*}
From now on, $L(f\times g,.)$ is the Rankin-Selberg L-function described in section 4 of \textbf{\cite{KoMiVa}} associated to the pair $(f,g)$ and $\mathcal{P}$ denotes the set of prime numbers. The family $\mathcal{F}$ has the same properties (at least conjecturally) as the family $\mathcal{G}$. The challenge lies in the fact that the analytic conductor $Q(f\times g)$ say of any $L(f\times g,.)$ in $\mathcal{F}(q)$ is large by comparison with the size of $\left\vert\mathcal{F}(q)\right\vert$; one has
\begin{equation*}
\frac{\log{Q(f\times g)}}{\log{\left\vert\mathcal{F}(q)\right\vert}}\rightarrow 2 \quad \text{as} \quad q\rightarrow +\infty
\end{equation*}
while for the families $\mathcal{G}$ and $\mathcal{H}_\pm$, one has
\begin{equation*}
\frac{\log{Q(\chi_{-8d})}}{\log{\left\vert\mathcal{G}(X)\right\vert}}\rightarrow 1, \frac{\log{Q(f)}}{\log{\left\vert\mathcal{H}_\pm(K)\right\vert}}\rightarrow 1\, \text{as respectively} \, X, K\rightarrow +\infty.
\end{equation*}
In particular, the second moment in our case (whose evaluation is necessary to apply the \textbf{mollification method}) is already critical (in the sense of \cite{Mi2}); this is not the case of the families $\mathcal{G}$ and $\mathcal{H}_\pm$, for which the fourth moment is critical. Moreover, the L-functions of the family $\mathcal{F}$ are Euler products of degree four (rather than one or two) which significantly increases the combinatorial analysis. \newline\newline
For $q$ in $\mathcal{P}$ and $k\geq 2$ an even integer, we define the following harmonic averaging operator
\begin{equation*}
\forall q\in\mathcal{P}, \quad A_q^{h}[\alpha]:=\sum_{f\in S_{k}^{p}(q)}^{\qquad h}\alpha_{f}:=\sum_{f\in S_{k}^{p}(q)}\omega_{q}(f)\alpha_{f}
\end{equation*}
for sequences of complex numbers indexed by $S_k^p(q)$ and with the harmonic weight $\omega_{q}(f):=\frac{\varGamma(k-1)}{(4\pi)^{k-1}\langle f,f\rangle_q}$ ($\langle .,.\rangle_q$ is the Petersson scalar product on the space of cusp forms of level $q$, weight $k$ and trivial nebentypus). We also define the harmonic probability measure on $S_{k}^{p}(q)$ by
\begin{equation*}
\mu_q^{h}(E):=\frac{1}{A_q^{h}[1]}\sum_{f\in E}\omega_{q}(f)
\end{equation*}
for any subset $E$ of $S_{k}^{p}(q)$. With these notations, our analogue of the theorem of J.B. Conrey and K. Soundararajan is:
\begin{theoA}
\label{theoA}
Let $g$ be a primitive cusp form of square-free level $D$, weight $k_{g}\geq 22$ and trivial nebentypus. As $q\rightarrow +\infty$ among primes and $f$ ranges over the set of primitive cusp forms of level $q$, weight $k\geq k_g+6$ and trivial nebentypus, there are infinitely many (at least $1.8\%$ in a suitable sense) $f$ in $S_k^p(q)$ such that $L(f\times g,.)$ has at most eight non-trivial real zeros. More precisely, for $q$ a prime coprime with $D$, and $k\geq k_g+6$, we have:
\begin{equation*}
\mu_q^{h}\left(\left\{f\in S_k^p(q), L(f\times g,.) \, \text{has at most $8$ zeros in } \left[0,1\right]\right\}\right)\geq 0.018+o_g(1).
\end{equation*}
\end{theoA}
\begin{remark}
Under the Ramanujan-Petersson-Selberg conjecture (confer $H_2(0)$ next page), we would obtain $4\%$ of $L(f\times g,.)$ having at most $6$ non-trivial real zeros. However, even this strong and deep hypothesis does not seem to give the existence of infinitely many Rankin-Selberg $L$-functions having no zeros in $\left[0,1\right]$ by the present method.
\end{remark}
\begin{remark}
In the course of the proof of theorem A, we also prove that the analytic rank of the family $\mathcal{F}$ is bounded on average. More precisely, set
\begin{equation}
\label{ar}
r(f\times g):=\text{ord}_{s=\frac{1}{2}}L(f\times g,s),
\end{equation}
one has
\begin{equation*}
\label{harmonicrank}
\frac{1}{A_q^h[1]}A_q^h[r(.\times g)]\leq 9.82+o_g(1)
\end{equation*}
and we can replace the constant $9.82$ by $7.66$ under Ramanujan-Petersson-Selberg conjecture. Moreover, following the method of \cite{H-BMi}, one can even show the exponential decay of the analytic rank of the family $\mathcal{F}$ namely there exists some absolute constants $B, C>0$ such that:
\begin{equation*}
\frac{1}{A_q^{h}[1]}A_q^{h}\left[\exp{\left(B r(.\times g)\right)}\right]\leq C.
\end{equation*}
\end{remark}
The proof of theorem A relies on some asymptotic formulas for the harmonic mollified second moment of the family $\mathcal{F}$, which is defined by
\begin{equation}
\mathcal{W}^{h}(g;\mu):=A_q^{h}\left[\left\vert\mathcal{L}\left(.\times g,\frac{1}{2}+\mu\right)\right\vert^{2}\right]
\end{equation}
where for $f\in S_k^p(q)$ and $s\in\mathbb{C}$ we have set
\begin{equation*}
\mathcal{L}(f\times g,s):=L(f\times g,s)M(f\times g,s);
\end{equation*}
here, $M(f\times g,.)$ is some Dirichlet polynomial (the so-called mollifier) of the following shape
\begin{equation*}
M(f\times g,s):=\sum_{1\leq\ell\leq L}\frac{x_\ell(g,s)}{\ell^s}\lambda_f(\ell)
\end{equation*}
where the length $L\geq 1$ has to be as large as possible. Here, the $\left(\lambda_f(\ell)\right)_{1\leq\ell\leq L}$ are Hecke eigenvalues of $f$ and the $\left(x_\ell(g,s)\right)_{1\leq\ell\leq L}$ are well chosen mollifying coefficients depending on $s$, $g$ on some parameter $0<\Upsilon<1$ and on some polynomial $P$ satisfying $P(0)=P^\prime(0)=P^\prime(\Upsilon)=0$ and $P(\Upsilon)=1$ (see section \ref{proof}). Our key technical result is an asymptotic formula for $\mathcal{W}^{h}(g;\mu)$ when
\begin{equation*}
\frac{\varepsilon_0}{\log{q}}\leq\vert\mu\vert\ll\frac{1}{\log{q}}
\end{equation*}
for some small absolute constant $\varepsilon_0>0$. Given $u$ and $v$ two real numbers and $\Delta>0$, we define:
\begin{multline*}
\mathcal{V}(u,v):=1+\frac{\exp{(-u)}}{\Delta}\left(\frac{\sinh{u}}{u}-\frac{\sin{v}}{v}\right) \\
\times\int_{0}^{\Upsilon}\exp{\left(-2u\Delta(1-x)\right)}\left\vert P^{\prime}(x)+\frac{P^{\prime\prime}(x)}{2(u+iv)\Delta}\right\vert^{2}\mathrm{d}x.
\end{multline*}
Our \textbf{main first result} is an asymptotic formula for $\mathcal{W}^{h}(g;\mu)$ in terms of $\mathcal{V}(u,v)$; namely for
\begin{equation*}
\Delta:=\frac{\log{L}}{\log{\left(q^2\right)}}
\end{equation*}
which we call the relative (logarithmic) length of the mollifier, one has
\begin{multline}
\label{afsecond}
\mathcal{W}^{h}(g;\mu)=\mathcal{V}(\log{\left(q^2\right)}\Re{(\mu)},\log{\left(q^2\right)}\Im{(\mu)})+\text{Errsec}(q,L;\mu) \\
+O_{k,g}\left(\frac{1}{q^\delta}+\frac{1}{\log{q}}\right.\begin{cases}
L^{-2\Re{(\mu)}(1-\Upsilon)} & \text{ if } \Re{(\mu)}\geq 0, \\
q^{-2\Re{(\mu)}}L^{-4\Re{(\mu)}} & \text{ otherwise}
\end{cases}
\Bigg)
\end{multline}
for $\delta>0$ an absolute constant and $\text{Errsec}(q,L;\mu)$ some error term:
\begin{equation}
\label{error}
\text{Errsec}(q,L;\mu)=O_{k,g}\left(\frac{1}{q^\alpha}\right)
\end{equation}
for some $\alpha>0$ as soon as $\Delta$ is small enough in which case $\Delta$ is said to be \textbf{effective}.
\begin{remark}
The asymptotic for the harmonic mollified second moment of this family is the same as the asymptotic for the mollified second moment of the family of Dirichlet L-functions considered by J.B. Conrey and K. Soundararajan. This is consistent with the Random Matrix Model, as these two families are expected to have the same symmetry type.
\end{remark}
\begin{remark}
\label{extension}
In fact, we also prove that \eqref{afsecond} holds with some weaker assumptions on $\mu$; namely when $\mu$ satisfies  $\frac{\varepsilon_0}{\log{q}}\leq\vert\mu\vert$, $-\frac{c}{\log{q}}\leq\Re{(\mu)}\leq\frac{f_1(q)}{\log{q}}$ and $\left\vert\Im{(\mu)}\right\vert\leq\frac{f_2(q)}{\log{q}}$ for some $\varepsilon_0>0$, $c>0$ and some non-negative functions $f_1$, $f_2$ with the following properties:
\begin{equation*}
\lim_{q\rightarrow +\infty}f_1(q)=+\infty,\,\,\,f_1(q)=o(\log{q}), \,\,\,f_2(q)=O(\log{q}).
\end{equation*}
In this case, \eqref{afsecond} becomes:
\begin{multline}
\label{afsecondbis}
\mathcal{W}^{h}(g;\mu)=\mathcal{V}(\log{\left(q^2\right)}\Re{(\mu)},\log{\left(q^2\right)}\Im{(\mu)})+\text{Errsec}(q,L;\mu) \\
+O_{k,g}\left(\frac{1}{q^\delta}+\frac{f_1(q)+f_2(q)}{\log{q}}\right.\begin{cases}
L^{-2\Re{(\mu)}(1-\Upsilon)} & \text{ if } \Re{(\mu)}\leq 0, \\
q^{-2\Re{(\mu)}}L^{-4\Re{(\mu)}} & \text{ otherwise}.
\end{cases}
\Bigg)
\end{multline}
\end{remark}
Our task now is to produce effective positive $\Delta$. The existence of such $\Delta$ is a consequence of the work of E. Kowalski, P. Michel and J. Vanderkam (\cite{KoMiVa}) and their result leads to:
\begin{propoC}
Let $g$  be a primitive cusp form of square-free level $D$ and trivial nebentypus. Assume that $q$ is prime, coprime with $D$. If $\vert\mu\vert\ll\frac{1}{\log{q}}$ then for any natural integer $L\geq 1$,
\begin{equation}
\label{Errsec1}
\text{Errsec}(q,L;\mu)=O_{\varepsilon,k,g}\left((qL)^\varepsilon\left(L^{\frac{5}{2}}q^{-\frac{1}{12}}+L^{\frac{21}{4}}q^{-\frac{1}{4}}\right)\right)
\end{equation}
for any $\varepsilon>0$.
In particular, every $\Delta<\frac{1}{60}=0.01666...$ is effective.
\end{propoC}
This is a consequence of an asymptotic formula for the harmonic twisted second moment of this family given by
\begin{equation}
\mathcal{M}_{g}^{h}(\mu;\ell):=A_q^{h}\left[L\left(.\times g,\frac{1}{2}+\mu\right)L\left(.\times g,\frac{1}{2}+\overline{\mu}\right)\lambda_{.}(\ell)\right]
\end{equation}
where $\mu\in\mathbb{C}$, $q\in\mathcal{P}$, $\ell\geq 1$ and $\lambda_{.}(\ell)$ is a Hecke eigenvalue. It is shown in \cite{KoMiVa} that (confer Theorem \ref{komiva} in this paper):
\begin{theo}[E. Kowalski-P. Michel-J. Vanderkam (2002)]
Let $g$  be a primitive cusp form of square-free level $D$ and trivial nebentypus and $\mu$ be a complex number. Assume that $q$ is prime, coprime with $D$. If $\vert\Re{(\mu)}\vert\ll\frac{1}{\log{q}}$ then for any natural integer $1\leq \ell<q$,
\begin{equation}
\label{twist}
(qD)^{2\Re{(\mu)}}\mathcal{M}_{g}^{h}(\mu;\ell)=\text{MT}(\mu)+\text{Errtwist}(q,\ell;\mu)
\end{equation}
where $\text{MT}(\mu)$ stands for the main term and is described in section \ref{around} and a bound for the error term is given by
\begin{equation}
\label{errtwist}
\text{Errtwist}(q,\ell;\mu)=O_{\varepsilon,k,g}\left((q\ell)^{\varepsilon}(1+\vert\Im{(\mu)}\vert)^{B}\left(\ell^{\frac{3}{4}}q^{-\frac{1}{12}}+\ell^{\frac{17}{8}}q^{-\frac{1}{4}}\right)\right)
\end{equation}
for some absolute constant $B>0$ and for any $\varepsilon>0$.
\end{theo}
Nevertheless, this is not sufficient to obtain Theorem A\footnote{with $\Delta<\frac{1}{60}$ we would obtain a positive proportion of $L(f\times g,.)$ having at most $22$ zeros on $[0,1]$.}. Our \textbf{second main input} is a large improvement of the effective value of $\Delta$ by the introduction of the spectral theory of automorphic forms. To state our result, we introduce the following hypothesis which measures the approximation towards the Ramanujan-Petersson-Selberg conjecture.
\begin{hypothesis}
For any cuspidal automorphic form $\pi$ on $GL_{2}(\mathbb{Q})\backslash GL_{2}(\mathbb{A}_{\mathbb{Q}})$ with local Hecke parameters $\alpha_{\pi}^{(1)}(p)$, $\alpha_{\pi}^{(2)}(p)$ for $p<\infty$ and $\mu_{\pi}^{(1)}(\infty)$, $\mu_{\pi}^{(2)}(\infty)$ at infinity, the following bounds are available:
\begin{eqnarray*}
\vert\alpha_{\pi}^{(j)}(p)\vert & \leq & p^{\theta}, \, j=1,2, \\
\left\vert\Re{\left(\mu_{\pi}^{(j)}(\infty)\right)}\right\vert & \leq & \theta, \, j=1,2,
\end{eqnarray*}
provided $\pi_{p}$, $\pi_{\infty}$ are unramified, respectively.
\end{hypothesis}
We say that $\theta$ is \textbf{admissible} if $H_2(\theta)$ is satisfied. At the moment, the smallest admissible value of $\theta$ is $\theta_{0}=\frac{7}{64}$ thanks to the works of H. Kim, F. Shahidi and P. Sarnak (confer \cite{KiSh} and \cite{KiSa}).\newpage
\begin{propoD}
\label{propoD}
Let $\alpha$ be in $\left]0,1\right[$. Let $g$  be a primitive cusp form of square-free level $D$, weight $k_{g}>1+\frac{5}{2(1-\alpha)}$ and trivial nebentypus and $\mu$ be a complex number. Assume that $q$ is prime, coprime with $D$ and that $k\geq k_g+6$. If $\,\theta$ is admissible and $\vert\Re{(\mu)}\vert\ll\frac{1}{\log{q}}$ then for any natural integer $\ell\geq 1$,
\begin{equation}
\label{twistbis}
\text{Errtwist}(q,\ell;\mu)=O_{\varepsilon,k,g}\left((q\ell)^{\varepsilon}(1+\vert\Im{(\mu)}\vert)^{B}\left(\ell^{2+\theta}q^{-\left(\frac{1}{2}-\theta\right)}+\ell^{\frac{9}{4}+\frac{\theta}{2}-\alpha}q^{-\left(\alpha-\frac{1}{2}-\theta\right)}\right)\right)
\end{equation}
and for any natural integer $L\geq 1$,
\begin{equation}
\label{Errsec2}
\text{Errsec}(q,L;\mu)=O_{\varepsilon,k,g}\left((qL)^{\varepsilon}(1+\vert\Im{(\mu)}\vert)^{B}\left(L^{2+2\theta}q^{-\left(\frac{1}{2}-\theta\right)}+L^{\frac{11}{2}+\theta-2\alpha}q^{-\left(\alpha-\frac{1}{2}-\theta\right)}\right)\right)
\end{equation}
for some absolute constant $B>0$ and for any $\varepsilon>0$. Consequently, under $H_2(\theta)$, every $\Delta<\Delta_{\text{max}}(\theta):=\frac{1-2\theta}{4(5+2\theta)}$ is effective granted that $k$ and $k_g$ are large enough.
\end{propoD}
\begin{remark}
We note that:
\begin{eqnarray*}
\Delta_{max}(\theta_{0}) & = &\frac{25}{668}=0.03742... \\
\Delta_{max}(0) & = &\frac{1}{20}=0.05. 
\end{eqnarray*}
\end{remark}
The error term in \eqref{twist} comes from the resolution of a shifted convolution problem by the authors, which builds on the $\delta$-symbol method of W. Duke, J.B. Friedlander and H. Iwaniec (\cite{DuFrIw}). This error term is improved using a technique of P. Sarnak (confer \cite{Sa}) which makes systematic use of spectral theory of automorphic forms (see section \ref{scp}). However, this method alone would only enable us to take $\Delta<\frac{1-2\theta}{8(4+\theta)}$ and we have to supplement it by additional refinements (in particular by considering the shifted convolution problem on average and detecting cancellations throughout large sieve inequalities) which lead to an effective length of $\frac{1-2\theta}{4(7+2\theta)}$. Finally, Proposition D is obtained thanks to an estimate of triple products on average over the spectrum of B. Krötz and R.J. Stanton (\textbf{\cite{KrSt}} and see also \textbf{\cite{Ko2}}). \newline\newline
Another consequence of our refinements is an improvement over the previously known subconvexity bounds for Rankin-Selberg $L$-functions in the level aspect obtained by the \textbf{amplification method}:
\begin{theoB}
\label{theoB}
Let $g$ be a primitive cusp form of square-free level $D$, weight $k_{g}\geq 20$ and trivial nebentypus. Let us assume that $q$ is a prime large enough and that $k\geq k_g+6$. If $\theta$ is admissible then for any natural integer $j$ and any $f$ in $S_{k}^{p}(q)$, we have
\begin{equation}
\left\vert L^{(j)}\left(f\times g,\frac{1}{2}+it\right)\right\vert\ll_{\varepsilon,k,j,g}\left(1+\vert t\vert\right)^{B}q^{\frac{1}{2}-\omega(\theta)+\varepsilon},
\end{equation}
for any $\varepsilon>0$ where $t$ is real, the exponent $B$ is absolute and $\omega(\theta):=\frac{1-2\theta}{4(9+4\theta)}$.
\end{theoB} 
\begin{remark}
In \cite{KoMiVa}, a subconvex bound is obtained but with $\omega(\theta)$ replaced by $\frac{1}{80}=0.0125$. Note that $\omega(\theta_0)=\frac{25}{1208}=0.020695...$ and that $\omega(0)=\frac{1}{36}=0.027777...$
\end{remark}

One may wonder what happens when one tries to remove the harmonic weights in Theorem A. In \cite{KoMi}, E. Kowalski and P. Michel provided a general technique to deduce asymptotic formulas for the natural average $A_q[\alpha]:=\sum_{f\in S_{k}^{p}(q)}\alpha_{f}$ from asymptotic formulas for the harmonic average as long as the coefficients $\alpha_{f}$ do not increase or oscillate too much as $q$ goes to infinity. In our case, one can deduced the same asymptotic formula as in \eqref{afsecond} for
\begin{equation*}
\mathcal{W}(g;\mu):=\frac{1}{A_q[1]}A_q\left[\left\vert\mathcal{L}\left(.\times g,\frac{1}{2}+\mu\right)\right\vert^2\right]
\end{equation*}
but with the length of the mollifier strictly smaller than $\omega(\theta)$. In other words, it seems that getting rid of the harmonic weights has a cost in this situation.  \newline

\noindent{\textbf{Notations.}}
From now on, $\mu$ will denote a complex number and $\tau=\Re{(\mu)}$, $t=\Im{(\mu)}$, $\delta=i\Im{(\mu)}$. We also set $\mu_1:=\mu$ and $\mu_2:=\overline{\mu}$. In several places, given an Euler product $L(s)=\prod_{p\in\mathcal{P}}L_{p}(s)$ , we write $L_{(N)}(s):=\prod_{p\mid N}L_{p}(s)$ and $L^{(N)}(s):=\prod_{p\nmid N}L_{p}(s)$ for any natural integer $N$. We set: $\log_2(x):=\log{(\log{x})}$. $\tau(n)$ equals the number of divisors of $n$ and $\mu(n)$ is the Möbius function at $n$. We will denote by $\varepsilon$ and $B>0$ some absolute positive constants whose definition may vary from line to line. The notations $f(q)\ll_{A}g(q)$ or $f(q)=O_A(q)$ mean that $\vert f(q)\vert$ is smaller than a constant which only depends on $A$ times $g(q)$ at least for $q$ large enough. Similarly, $f(q)=o(1)$ means that $\lim_{q\rightarrow +\infty}f(q)=0$. Finally, if $\mathcal{E}$ is a property, the Kronecker symbol $\delta_{\mathcal{E}}$ equals $1$ if $\mathcal{E}$ is satisfied and $0$ else. \newline\newline
For all background and notations about classical modular forms and Rankin-Selberg $L$-functions, we refer the reader to sections 3 and 4 of \textbf{\cite{KoMiVa}} and to Appendix C.\newline\newline
\noindent{\textbf{Acknowledgments.}}
I sincerely thank my advisor, Professor Philippe Michel, for all his comments and remarks which got the better of my doubts. I also think of Professors Etienne Fouvry and Emmanuel Kowalski for their advices and encouragements. I wish to thank the Fields Institute of Toronto, where part of this work was done, for the excellent working conditions. I also acknowledge the referee for a careful reading of the manuscript.

\section{A review of classical modular forms}

\label{modforms}
In this section, we recall general facts about modular forms. The main reference is \cite{Iw2}. For $N\geq 1$, we consider $\varGamma_{0}(N)$ the congruence subgroup of level $N$ and $\varepsilon_{N}$ the trivial Dirichlet character of modulus $N$. All elements of $GL_{2}^{+}(\mathbb{R})$ act on the upper-half plane $\mathbb{H}$ by linear-fractional transformations and this defines an action of the group $SL_{2}(\mathbb{R})$ on it. For $\gamma=\begin{pmatrix} a & b \\ c & d \end{pmatrix}$ in $GL_{2}^{+}(\mathbb{R})$ and $z$ in $\mathbb{H}$ we set $j(\gamma,z):=cz+d$. Let $m$ be an even natural integer. For $\gamma$ in $GL_{2}^{+}(\mathbb{R})$ and $h:\mathbb{H}\rightarrow \mathbb{C}$, we define:
\begin{equation*}
\forall z\in\mathbb{H}, \quad h\arrowvert_{\gamma}^{m}(z):=\frac{(\text{det}(\gamma))^{\frac{m}{2}}}{j(\gamma,z)^{m}}h(\gamma.z).
\end{equation*}
This formula clearly defines an action of $SL_{2}(\mathbb{R})$ on the space of complex valued functions on $\mathbb{H}$, which is said to be of weight $m$. \newline

\subsection{Cusp forms}

A holomorphic function $h:\mathbb{H}\mapsto \mathbb{C}$ which satisfies:
\begin{equation*}
\forall\gamma\in\varGamma_{0}(N),\quad h\arrowvert_{\gamma}^{m}=h
\end{equation*}
and is holomorphic at the cusps of $\varGamma_{0}(N)$ is a \emph{modular form} of level $N$, weight $m$ and trivial nebentypus $\varepsilon_{N}$. Such a modular form is a \emph{cusp form} if $y^{\frac{m}{2}}h(z)$ is bounded on the upper-half plane. We denote by $S_{m}(N)$ this set of cusp forms which is equipped with the Petersson inner product:
\begin{equation*}
\left\langle h_{1},h_{2}\right\rangle_N=\int_{\varGamma_0(N)\setminus\mathbb{H}}y^{m}h_{1}(z)\overline{h_{2}(z)}\frac{\mathrm{d}x\mathrm{d}y}{y^{2}}.
\end{equation*}
One can obtain the Fourier expansion at infinity of each such cusp form $h$:
\begin{equation*}
\forall z\in\mathbb{H},\quad h(z)=\sum_{n\geq1}\psi_{h}(n)n^{\frac{m-1}{2}}e(nz)
\end{equation*}
where $e(z):=\exp{(2i\pi z)}$.

\subsection{Hecke operators}

For every natural integer $\ell\geq 1$, the \emph{Hecke operator} of weight $m$, nebentypus $\varepsilon_{N}$ and rank $\ell$ on $S_{m}(N)$ is defined by:
\begin{equation*}
\forall z\in\mathbb{H}, \quad\left(T_{\ell}(h)\right)(z):=\frac{1}{\sqrt{\ell}}\sum_{ad=\ell}\varepsilon_{N}(a)\sum_{0\leq b<d}h\left(\frac{az+b}{d}\right).
\end{equation*}
Thus, we remark that $T_{\ell}$ is independent of $m$ and we can prove that it is hermitian if $\text{gcd}(\ell,N)=1$. Moreover, we can show that the algebra spanned by the Hecke operators is a commutative one. More precisely, we have the following composition property:
\begin{equation}
\label{compo}
\forall(\ell_1,\ell_2)\in (\mathbb{N}^*)^2, \quad T_{\ell_{1}}\circ T_{\ell_{2}}=\sum_{d\mid\left(\ell_{1},\ell_{2}\right)}\varepsilon_{N}(d)T_{\frac{\ell_{1}\ell_{2}}{d^{2}}}.
\end{equation}
A cusp form which is also an eigenfunction of the $T_{\ell}$ for $\text{gcd}(\ell,N)=1$ is called a \emph{Hecke} cusp form and an orthonormal basis of $S_{m}(N)$ made of Hecke cusp forms is called a \emph{Hecke eigenbasis}.

\noindent{\textbf{Atkin-Lehner theory.}}
The main reference of this part is \cite{AtLe}. Briefly speaking, we obtain, with the previous notations, a splitting of $S_{m}(N)$ in $S_{m}^{o}(N)\oplus^{\perp_{\langle .,\rangle_N}}S_{m}^{n}(N)$ where:
\begin{eqnarray}
\label{oldnew}
S_{m}^{o}(N)\! & = & \! \text{Vect}_{\mathbb{C}}\left\{g(dz), N^{\prime}\mid N, d\mid \frac{N}{N^\prime}, d\neq 1, g\in S_{m}(N^{\prime})\right\}, \\
S_{m}^{n}(N)\! & = & \! \left(S_{m}^{o}(N)\right)^{\perp_{\langle .,.\rangle_N}}
\end{eqnarray}
where "$o$" stands for old and "$n$" for new. These two spaces are invariant under the action of the Hecke operators $T_{l}$ for $\text{gcd}(l,N)=1$. A \emph{primitive} cusp form $h$ is a Hecke cusp form which is new and satisfies:
\begin{equation*}
\psi_{h}(1)=1.
\end{equation*}
Such an element $h$ is automatically an eigenfunction of the other Hecke operators and also of the Atkin-Lehner operators which will be defined later and satisfies $\psi_{h}(\ell)=\lambda_{h}(\ell)$ for all integer $\ell$ where $T_\ell(h)=\lambda_h(\ell)h$ ($\lambda_{h}(\ell)$ is the \emph{Hecke eigenvalue} of rank $\ell$). The set of primitive cusp forms will be denoted by $S_{m}^{p}(N)$. Let $h$ be a cusp form with Hecke eigenvalues $\left(\lambda_{h}(\ell)\right)_{\left(\ell,N\right)=1}$. The composition property \eqref{compo} of the Hecke operators entails that for all $\ell_{1}$ and for $\text{gcd}\left(\ell_{2},N\right)=1$:
\begin{eqnarray}
\label{compoeigen}
\psi_{h}(\ell_{1})\lambda_{h}(\ell_{2}) & = & \sum_{d\mid\left(\ell_{1},\ell_{2}\right)}\varepsilon_{N}(d)\psi_{h}\left(\frac{\ell_{1}\ell_{2}}{d^{2}}\right), \\
\psi_{h}(\ell_{1}\ell_{2}) & = & \sum_{d\mid\left(\ell_{1},\ell_{2}\right)}\mu(d)\varepsilon_{N}(d)\psi_{h}\left(\frac{\ell_{1}}{d}\right)\lambda_{h}\left(\frac{\ell_{2}}{d}\right)
\end{eqnarray}
and this relation holds for all $\ell_{1}, \ell_{2}$ if $h$ is primitive. The adjointness relation is:
\begin{equation}
\label{adjointness}
\forall\text{gcd}(\ell,N)=1,\quad\lambda_{h}(\ell)=\overline{\lambda_{h}(\ell)}, \quad \psi_{h}(\ell)=\overline{\psi_{h}(\ell)}
\end{equation}
and this remains true for all $\ell$ if $h$ is a primitive cusp form.

\subsection{Bounds for Hecke eigenvalues of cusp forms}

Let $h$ be a primitive cusp form of level $N$, weight $m$ and trivial nebentypus $\varepsilon_{N}$. Remember that:
\begin{equation*}
\forall\ell\in\mathbb{N}^{*},\quad T_{\ell}h=\lambda_{h}(\ell)h.
\end{equation*}
For a prime $p$, let $\alpha_{h,1}(p)$ and $\alpha_{h,2}(p)$ be the complex roots of the following quadratic equation:
\begin{equation*}
X^{2}-\lambda_{h}(p)X+\varepsilon_{N}(p)=0.
\end{equation*}
It follows from the work of Eichler-Shimura-Igusa and Deligne that the Ramanujan-Petersson bound holds true:
\begin{equation}
\label{individualh}
\vert \alpha_{h,1}(p)\vert, \vert \alpha_{h,2}(p)\vert\leq 1 \,\,\,\textrm{and so}\,\,\,\forall\ell\geq 1, \vert\lambda_{h}(\ell)\vert\leq \tau(\ell).
\end{equation}
Setting $\sigma_h(n):=\sum_{d\mid n}\vert\lambda_h(d)\vert$, it entails that:
\begin{equation}
\forall X>0,\quad\sum_{n\leq X}\sigma_h(n)^2\ll_{\varepsilon, h} X^{1+\varepsilon}.
\end{equation}
for all $\varepsilon>0$.

\subsection{Atkin-Lehner operators}

The results of this part were established by A. Atkin and J. Lehner. We assume that $N=N_{1}N_{2}$ with $\text{gcd}\left(N_{1},N_{2}\right)=1$. Let $x, y, z, w$ four integers satisfying:
\begin{eqnarray*}
y & \equiv & 1\mod\left(N_{1}\right), \\
x & \equiv & 1\mod\left(N_{2}\right), \\
N_{1}^{2}xw-Nyz & = & N_{1}.
\end{eqnarray*}
If $\omega_{\vrule width 0pt height 7pt depth 10pt N_{ 1}}=\begin{pmatrix}
xN_{1} & y \\ 
zN & wN_{1}
\end{pmatrix}$ then $W_{\vrule width 0pt height 7pt depth 10pt N_{1}}=\left\arrowvert_{\vrule width 0pt height 5pt depth 0pt\omega_{\vrule width 0pt height 0pt depth 0pt N_{1}}}^{\vrule width 0pt height 5pt depth 0pt m}\right.$ is a linear endomorphism of $S_m(N)$ independent of $x$, $y$, $z$ and $w$. If $N_{1}=N$ then $W_{\vrule width 0pt height 7pt depth 10pt N_{1}}$ is the classical Fricke involution given by $\omega_{\vrule width 0pt height 7pt depth 10pt N}=\begin{pmatrix}
0 & 1 \\
-N & 0
\end{pmatrix}$. The following proposition holds:
\begin{proposition}[A. Atkin-W. Li (1970)]
\label{refAtLe}
If $N_{1}\mid N$ and $\text{gcd}\left(N_{1},\frac{N}{N_{1}}\right)=1$ then:
\begin{equation*}
\forall h\in S_{m}^n(N),\quad W_{N_{1}}h=\eta_{h}(N_{1})h
\end{equation*}
where $\eta_{h}(N_{1})=\pm 1$. 
\end{proposition}

\section{A review of Rankin-Selberg L-functions}

\label{L-func}
Throughout this section, $g_{1}$ belongs to $S_{k_{1}}^{p}(D_{1})$ and $g_{2}$ belongs to $S_{k_{2}}^{p}(D_{2})$.

\subsection{About Rankin-Selberg L-functions}
The \emph{Rankin-Selberg $L$-function} of $g_{1}$ and $g_{2}$ is the following $L$-function defined a priori for $\Re{(s)}>1$ by:
\begin{equation*}
L(g_{1}\times g_{2},s):=\zeta^{(D_{1}D_{2})}(2s)\sum_{l\geq 1}\frac{\lambda_{g_{1}}(l)\lambda_{g_{2}}(l) }{l^{s}}.
\end{equation*}
It admits an Eulerian product $L(g_{1}\times g_{2},.):=\prod_{p\in\mathcal{P}}L_{p}(g_{1}\times g_{2},.)$ where:
\begin{equation*}
\forall p\in\mathcal{P}, \forall s\in\mathbb{C},\quad L_{p}(g_{1}\times g_{2},s)=\prod_{1\leq i,j\leq 2}\left(1-\alpha_{g_{1},i}(p)\alpha_{g_{2},j}(p)p^{-s}\right)^{-1}.
\end{equation*}
By Rankin-Selberg theory, $L(g_{1}\times g_{2},.)$ admits a meromorphic continuation to the complex plane with at most simple poles at $s=0,1$ which occur only if $g_{1}=g_{2}$. This L-function satisfies a functional equation. When $\text{gcd}\left(D_{1},D_{2}\right)=1$, it takes the following form. We set:
\begin{equation*}
\forall s\in\mathbb{C}, \quad\Lambda\left(g_{1}\times g_{2},s\right):=\left(\frac{D_{1}D_{2}}{4\pi^{2}}\right)^{s}\varGamma\left(s+\frac{\left\vert k_{1}-k_{2}\right\vert}{2}\right)\varGamma\left(s+\frac{k_{1}+k_{2}}{2}-1\right)L\left(g_{1}\times g_{2},s\right).
\end{equation*}
The functional equation is then:
\begin{equation*}
\forall s\in\mathbb{C}, \quad \Lambda\left(g_{1}\times g_{2},s\right)=\varepsilon\left(g_{1}\times g_{2}\right)\Lambda\left(g_{1}\times g_{2},1-s\right)
\end{equation*}
where the sign of the functional equation in our case is $\varepsilon(g_{1}\times g_{2})=1$.

\subsection{About symmetric square L-functions}
Closely related to $L(g_{1}\times g_{1},.)$ is the following Dirichlet series defined for $\Re{(s)}>1$:
\begin{equation*}
L(\text{Sym}^{2}(g_{1}),s):=\zeta^{(D_{1})}(2s)\sum_{l\geq 1}\frac{\lambda_{g_{1}}(l^{2})}{l^{s}}.
\end{equation*}
The Eulerian product of $L(\text{Sym}^{2}(g_{1}),.)$ is given by $\prod_{p\in\mathcal{P}}L_{p}(\text{Sym}^{2}(g_{1}),.)$ with:
\begin{equation*}
\forall p\in\mathcal{P}, \forall s\in\mathbb{C}, \quad L_{p}(\text{Sym}^{2}(g_{1}),s)=\prod_{1\leq i\leq j\leq 2}\left(1-\alpha_{g_{1},i}(p)\alpha_{g_{1},j}(p)p^{-s}\right)^{-1}.
\end{equation*}
Hence, we get $L(g_{1}\times g_{1},s)=\zeta^{(D_{1})}(s)L(\text{Sym}^{2}g_{1},s)$ for all complex number $s$.

\section{Proof of Theorem A and estimates for the analytic rank}

\label{proof}

\subsection{Principle of the proof}

Let $b>0$, $c>0$ and $\sigma_0>1$ some real numbers. $\mathcal{B}(q)$ will denote the rectangular box with vertices $\left(\frac{1}{2}-\frac{c}{\log{q}},\pm\frac{b}{\log{q}}\right)$ and $\left(\sigma_0,\pm\frac{b}{\log{q}}\right)$. Let $N$ be a natural integer and $S_{k}^{p}(q,N)$ be the set of primitive cusp forms $f$ in $S_{k}^{p}(q)$ whose Rankin-Selberg L-function $L(f\times g,.)$ admits
\begin{itemize}
\item
a zero of order $2n_1$ at $\frac{1}{2}$,
\item
and  $n_2$ zeros (counted with multiplicities) in $\left]\frac{1}{2},1\right]$
\end{itemize}
such that $2n_1+2n_2\geq 2(N+1)$. Let us remark that $S_k^p(q)\backslash S_k^p(q,N)$ is precisely the set of modular forms $f$ in $S_k^p(q)$ whose Rankin-Selberg $L$-function $L(f\times g,.)$ has at most $2N$ zeros in $\left[0,1\right]$. We are producing some $N$ such that (as $q$ tends to infinity among the primes)
\begin{equation*}
\mu_q^{h}\left(S_{k}^{p}(q,N)\right)\leq s_0(N)+o_g(1)
\end{equation*}
with $s_0(N)<1$ a constant which depends only on $N$ and conclude that for at least $100(1-s_{0}(N))$ percent of primitive cusp forms of weight $k$ and trivial nebentypus, $L(f\times g,.)$ has at most $2N$ non-trivial zeros on the real axis (and in fact in a small box $\mathcal{B}(q)$).

\subsection{Selberg's lemma}
\label{sectionselberg}
\begin{lemma}
\label{selberg}
Let $\psi$ be a holomorphic function which does not vanish on a half plane $\Re{(z)}\geq W$. Let $\mathcal{B}$ be the rectangular box of vertices $W_{0}\pm iH$, $W_{1}\pm iH$ where $H>0$ and $W_{0}<W<W_{1}$. We have:
\begin{multline*}
4H\sum_{\substack{\beta+i\gamma\in\mathcal{B} \\
\psi(\beta+i\gamma)=0}}
\cos{\left(\frac{\pi\gamma}{2H}\right)}\sinh{\left(\frac{\pi(\beta-W_{0})}{2H}\right)}=\int_{-H}^{H}\cos{\left(\frac{\pi t}{2H}\right)}\log{\vert \psi(W_{0}+it)\vert}\mathrm{d}t \\
+\int_{W_{0}}^{W_{1}}\sinh{\left(\frac{\pi(\alpha-W_{0}}{2H}\right)}\log{\vert \psi(\alpha+iH)\psi(\alpha-iH)\vert}\mathrm{d}\alpha \\
-\Re{\left(\int_{-H}^{H}\cos{\left(\pi\frac{W_{1}-W_{0}+it}{2iH}\right)}(\log{\psi})(W_{1}+it)\right)}\mathrm{d}t.
\end{multline*}
\end{lemma}
A proof of this lemma is given in \textbf{\cite{CoSo}} and relies on the fact that
\begin{equation*}
\int_{\partial\mathcal{B}}k(s)(\log{f})(s)\mathrm{d}s=0
\end{equation*}
with $k(s):=\cos{\left(\pi\frac{s-W_0}{2iH}\right)}$. Let us mention the properties which will be useful to us:
\begin{itemize}
\item
$k$ is purely imaginary on $\Im{(s)}=H$ and satisfies over there $k(s)=-k\left(\overline{s}\right)$,
\item
$\Re{(k)}\geq 0$ in $\mathcal{B}$.
\end{itemize}

\subsection{The successive steps}
\label{step}
We follow the method of J.B. Conrey and K. Soundararajan (\textbf{\cite{CoSo}}). Lemma \ref{selberg} applied to the box $\mathcal{B}(q)$ and the function $L(f\times g,.)$ entails that
\begin{equation}
\label{selbergappli}
4b\sinh{\left(\frac{\pi c}{2b}\right)}r(f\times g)+4b\sum_{\substack{\beta\geq\frac{1}{2}-\frac{c}{\log{q}} \\
\beta\neq\frac{1}{2} \\
L(f\times g,\beta)=0}}
\sinh{\left(\frac{\pi\left(c+\log{q}\left(\beta-\frac{1}{2}\right)\right)}{2b}\right)}\leq\sum_{j=1}^{3}I_{f}^{q}(j)
\end{equation}
where
\begin{equation*}
I_{f}^{q}(1)=\int_{-b}^{b}\cos{\left(\frac{\pi t}{2b}\right)}\log{\left\vert L\left(f\times g,\frac{1}{2}-\frac{c}{\log{q}}+i\frac{t}{\log{q}}\right)\right\vert}\mathrm{d}t,
\end{equation*}
\begin{equation*}
I_{f}^{q}(2)=\int_{-c}^{\left(\sigma_{0}-\frac{1}{2}\right)\log{q}}\sinh{\left(\frac{\pi(x+c)}{2b}\right)}
\log{\left\vert L\left(f\times g,\frac{1}{2}+\frac{x}{\log{q}}+i\frac{b}{\log{q}}\right)\right\vert^{2}}\mathrm{d}x,
\end{equation*}
\begin{equation*}
I_{f}^{q}(3)=-\Re\left(\int_{-b}^{b}\cos{\left(\pi\frac{\left(\sigma_{0}-\frac{1}{2}\right)\log{q}+c+it}{2ib}\right)}\left(\log{L(f\times g,.)}\right)\left(\sigma_{0}+i\frac{t}{\log{q}}\right)\mathrm{d}t\right).
\end{equation*}
One can show (confer \cite{Ri}) that if $f$ belongs to $S_k^p(q,N)$ then the left-hand side of \eqref{selbergappli} is larger than $(N+1)\times 8b\sinh{\left(\frac{\pi c}{2b}\right)}$. Thus,
\begin{equation*}
\mu_q^{h}\left(S_{k}^{p}(q,N)\right)\leq\frac{1}{N+1}\frac{1}{8b\sinh{\left(\frac{\pi c}{2b}\right)}}\frac{1}{A_q^{h}[1]}A_q^{h}\left[\sum_{j=1}^{3}I_{.}^{q}(j)\right].
\end{equation*}
The concavity of the $\log$ function leads to
\begin{equation}
\label{firstzero}
\mu_q^{h}\left(S_{k}^{p}(q,N)\right)\leq\frac{1}{N+1}\frac{1}{8b\sinh{\left(\frac{\pi c}{2b}\right)}}\left(J_{1}^{q,h}+J_{2}^{q,h}+\frac{1}{A_q^{h}[1]}A_q^{h}[I_{.}^{q}(3)]\right)
\end{equation}
where
\begin{eqnarray*}
J_{1}^{q,h} & := & \int_{0}^{b}\cos{\left(\frac{\pi t}{2b}\right)}\log{\left(\frac{1}{A_q^{h}[1]}A_q^{h}\left[\left\vert L\left(.\times g,\frac{1}{2}+\frac{-c+it}{\log{q}}\right)\right\vert^{2}\right]\right)}\mathrm{d}t, \\
J_{2}^{q,h} & := & \int_{-c}^{\left(\sigma_{0}-\frac{1}{2}\right)\log{q}}\sinh{\left(\frac{\pi(x+c)}{2b}\right)}\log{\left(\frac{1}{A_q^{h}[1]}A_q^{h}\left[\left\vert L\left(.\times g,\frac{1}{2}+\frac{x+ib}{\log{q}}\right)\right\vert^{2}\right]\right)}\mathrm{d}x.
\end{eqnarray*}
Similarly, we have from \eqref{selbergappli}:
\begin{equation}
\label{rank0}
\frac{1}{A_q^{h}[1]}A_q^{h}\left[r(.\times g)\right]\leq 2\frac{1}{8b\sinh{\left(\frac{\pi c}{2b}\right)}}\left(J_{1}^{q,h}+J_{2}^{q,h}+\frac{1}{A_q^{h}[1]}A_q^{h}[I_{.}^{q}(3)]\right).
\end{equation}
We need the right-hand side of \eqref{firstzero} and \eqref{rank0} to be small. Unfortunately, the weight function $\sinh{\left(\frac{\pi(x+c)}{2b}\right)}$ which appears in $J_{2}^{q,h}$ grows exponentially on the horizontal sides of the box. This problem is solved by mollifying: one replaces $L(f\times g,.)$ by $\mathcal{L}(f\times g,.)$ such that the exponential growth of $\sinh{\left(\frac{\pi(x+c)}{2b}\right)}$ is balanced by the exponential decay of $\log{\left(\frac{1}{A_q^{h}[1]}\mathcal{W}^{h}\left(g;\frac{x+ib}{\log{q}}\right)\right)}$.
\begin{remark}
Naively, one would like to be able to choose a kernel $k$ having the properties listed above in section \ref{sectionselberg} and such that the corresponding weight function (in $J_{2}^{q,h}$) does not grow exponentially on the horizontal sides of $\mathcal{B}(q)$. Unfortunately, as K. Soundararajan remarked at the Journées Arithmétiques 2003 in Graz, such kernel does not exist. So, it appears that the mollification step is a necessity.
\end{remark}

\subsection{Choosing the mollifier}

Note that on the half plane $\Re{(s)}>1$, $L(f\times g,s)=\sum_{n\geq 1}a_{f\times g}(n)n^{-s}$ where
\begin{equation*}
a_{f\times g}(n)=\sum_{n=n_1^{2}n_2}\varepsilon_{q}(n_1)\varepsilon_{D}(n_1)\lambda_{f}(n_2)\lambda_{g}(n_2)
\end{equation*}
satisfies $\vert a_{f\times g}(n)\vert\ll_\varepsilon n^\varepsilon$ for any $\varepsilon>0$. We need the Dirichlet coefficients of the inverse of $L(f\times g,.)$:
\begin{lemma}
\label{inverse}
For $\,\Re{(s)}>1$ one has
\begin{equation*}
\frac{1}{L(f\times g,s)}\!\!\!=\!\!\!K(g,2s)\!\!\!\sum_{\ell=\ell_{1}\ell_{2}^{2}\ell_{3}^{3}}\!\!\!\frac{\mu^{2}(\ell_{1}\ell_{2}\ell_{3})\mu(\ell_{1}\ell_{3})
\varepsilon_{q}(\ell_{3})\varepsilon_{D}(\ell_{2}\ell_{3})\lambda_{g}(\ell_{1}\ell_{3})}{K_{(\ell)}(g,s)\ell^{s}}\!\!\!\lambda_{f}(\ell_{1}\ell_{2}^{2}\ell_{3})
\end{equation*}
where $K(g,s):=\prod_{p\in\mathcal{P}}K_p(g,s)$ is an absolutely convergent Euler product on $\Re{(s)}>2$ given by
\begin{equation*}
\forall p\in\mathcal{P}, \quad K_{p}(g,s):=1+\varepsilon_{q}(p)\lambda_{g}(p^{2})p^{-s}+\varepsilon_{qD}(p)p^{-2s}.
\end{equation*}
\end{lemma}
\noindent{\textbf{Proof of lemma \ref{inverse}.}} 
We give no details. Setting $L(f\times g,s)^{-1}:=\sum_{\ell\geq 1}u_\ell\ell^{-s}$, one shows that $u_\ell=0$ except if $\ell=\ell_1\ell_2^2\ell_3^3\ell_4^4$ with $\ell_1$,  $\ell_2$, $\ell_3$, $\ell_4 $ square-free numbers pairwise coprime: 
\begin{equation*}
u_\ell=\mu(\ell_1\ell_3)\varepsilon_{qD}(\ell_3\ell_4)\lambda_{f}(\ell_1\ell_3)\lambda_{g}(\ell_1\ell_3)\sum_{\ell_2=\ell_2^\prime \ell_2^{\prime\prime}}\varepsilon_{q}(\ell_2^{\prime\prime})\varepsilon_{D}(\ell_2^\prime)\lambda_{f}(\ell_{2}^{\prime 2})\lambda_{g}(\ell_{2}^{\prime\prime 2}).
\end{equation*}
Note that $K(g,s)$ is an absolutely convergent Euler on $\Re{(s)}>2$ as:
\begin{equation*}
\forall p\in\mathcal{P}, \quad K_{p}(g,s):=1+O\left(\frac{1}{p^{\Re{(s)}}}\right).
\end{equation*}
\begin{flushright}
$\blacksquare$
\end{flushright}
Let  $0<\Upsilon<1$ be a real number and $P$ be a polynomial satisfying $P(\Upsilon)=1$, $P(0)=P^{\prime}(0)=P^{\prime}(\Upsilon)=0$. Let $L\geq 1$ be a natural integer. We set:
\begin{equation*}
F_{L}^{\Upsilon}(\ell)=\left\{\begin{array}{cl}
1 & \mbox{if } 1\leq \ell\leq L^{1-\Upsilon} \\
P\left(\frac{\log{\left(\frac{L}{\ell}\right)}}{\log{L}}\right) & \mbox{if } L^{1-\Upsilon}\leq \ell \leq L \\
0 & \mbox{else.}
\end{array}
\right.
\end{equation*}
The mollifier we choose is
\begin{multline*}
M(f\times g,s)=K(g,2s)\sum_{\ell=\ell_{1}\ell_{2}^{2}\ell_{3}^{3}}\frac{\mu^{2}(\ell_{1}\ell_{2}\ell_{3})\mu(\ell_{1}\ell_{3})\varepsilon_{D}(\ell_{3})\varepsilon_{q}(\ell_{2}\ell_{3})\lambda_{g}(\ell_{1}\ell_{3})}{K_{(\ell)}(g,2s)^{-1}\ell^s} \\
\lambda_{f}(\ell_1\ell_2^2\ell_3)F_{L}^{\Upsilon}(\ell_{1}\ell_{2}^{2}\ell_{3}) \\
=\sum_{\ell\geq 1}\frac{x_{\ell}(g,s)}{\ell^{s}}\lambda_{f}(\ell)\vrule width 7.2cm height 0pt depth 0pt
\end{multline*}
where
\begin{equation*}
x_{\ell}(g,s)=K(g,2s)\sum_{\ell=\ell_{1}\ell_{2}^{2}\ell_{3}}
\frac{\mu^{2}(\ell_{1}\ell_{2}\ell_{3})\mu(\ell_{1}\ell_{3})\varepsilon_{q}(\ell_{3})\varepsilon_{D}(\ell_{2}\ell_{3})\lambda_{g}(\ell_{1}\ell_{3})F_{L}^{\Upsilon}(\ell)}{K_{(\ell)}(g,2s)\ell_{3}^{2s}}
\end{equation*}
so that $M(f\times g,.)$ is a Dirichlet polynomial of length at most $L$ which approximates $L(f\times g,.)^{-1}$. From the shape of the mollifier, we immediately deduce
\begin{equation}
\label{loin}
\mathcal{L}(f\times g,s)=1+O_{\varepsilon}\left(L^{(1-\Upsilon)(1+\varepsilon-\Re{(s)})}\right)
\end{equation}
on $\Re{(s)}>1+\varepsilon$ for any $\varepsilon>0$. As a consequence, $\mathcal{L}(f\times g,.)$ has no zeros to the right of $1+\frac{\log_{2}{q}}{\log{q}}$ (at least for $q$ large enough). Moreover, for $q$ large enough, $\vert\mathcal{L}(f\times g,s)-1\vert<1$ on $\Re{(s)}>1+\varepsilon$ and we choose the branch of the logarithm given by
\begin{equation*}
\left(\log\mathcal{L}(f\times g,.)\right)(s):=\sum_{n\geq 1}\frac{(-1)^{n+1}}{n}(\mathcal{L}(f\times g,s)-1)^n
\end{equation*}
on $\Re{(s)}>1+\varepsilon$. We are going to give a useful integral expression of the coefficients $x_\ell(g,s)$ of the mollifier following a technique introduced in \textbf{\cite{KoMiVa2}}. To each polynomial $A(X)=\sum_{k\geq 0}a_{k}X^{k}$ and to each real number $M$, we associate the following transform:
\begin{equation*}
\forall s\in\mathbb{C}, \quad\widehat{A_{M}}(s)=\sum_{k\geq 0}a_{k}\frac{k!}{(s\log{M})^{k}}.
\end{equation*}
We have the following result:
\begin{lemma}
\label{lemmainte}
Let $m\geq 1$ be a natural integer.
\begin{equation*}
\frac{1}{2i\pi\log{M}}\int_{(3)}\left(\frac{M}{m}\right)^{s}\widehat{A_{M}}(s)\frac{\mathrm{d}s}{s^{2}}=\delta_{m<M}\left(\int^{(1)}A\right)\left(\frac{\log{\left(\frac{M}{m}\right)}}{\log{M}}\right)
\end{equation*}
where $\int^{(1)}A$ is the first antiderivative of $A$ without constant of integration.
\end{lemma}
\noindent{\textbf{Proof of lemma \ref{lemmainte}.}}
By linearity, it is enough to prove this lemma for $A(X)=X^k$ with $k\in\mathbb{N}^*$. Setting $y=\frac{M}{m}$, it consists in proving that
\begin{equation*}
\frac{1}{2i\pi}\int_{(3)}y^{s}\frac{\mathrm{d}s}{s^{k+2}}=\delta_{y>1}\frac{\log^{k+1}{(y)}}{(k+1)!}
\end{equation*} 
which is standard using suitable contour shifts (confer \textbf{\cite{KoMiVa2}}).
\begin{flushright} $\blacksquare$ \end{flushright}

To the polynomial $P$, we associate $R(X):=P((1-\Upsilon)X+\Upsilon)-1$ and we have the integral expression of the coefficients of the mollifier:
\begin{proposition}
\label{intemolli}
Let $\ell\geq 1$ be a natural integer and $s$ be a complex number.
\begin{multline*}
x_\ell(g,s)=\frac{1}{2i\pi\log{L}}\int_{(3)}L^{\left(1-\Upsilon\right)s}\left(\widehat{P^{\prime}_{L}}(s)L^{\Upsilon s}-\frac{1}{1-\Upsilon}\widehat{R^{\prime}_{L^{1-\Upsilon}}}(s)\right)K(g,2s) \\
\times\sum_{\ell=\ell_{1}\ell_{2}^{2}\ell_{3}}
\frac{\mu^{2}(\ell_{1}\ell_{2}\ell_{3})\mu(\ell_{1}\ell_{3})\varepsilon_{q}(\ell_{3})\varepsilon_{D}(\ell_{2}\ell_{3})\lambda_{g}(\ell_{1}\ell_{3})}{K_{(\ell)}(g,2s)\ell_{3}^{2s}\ell^s}\frac{\mathrm{d}s}{s^2}.
\end{multline*}
\end{proposition}
\noindent{\textbf{Proof of proposition \ref{intemolli}.}}
The main point is that we have:
\begin{equation*}
F_{L}^{\Upsilon}(\ell)=\frac{1}{2i\pi\log{L}}\int_{(3)}\left(\frac{L^{1-\Upsilon}}{\ell}\right)^s\left(\widehat{P^{\prime}_{L}}(s)L^{\Upsilon s}-\frac{1}{1-\Upsilon}\widehat{R^{\prime}_{L^{1-\Upsilon}}}(s)\right)\frac{\mathrm{d}s}{s^2}.
\end{equation*}
The previous integral expression of $F_L^\Upsilon(\ell)$ is a direct consequence of lemma \ref{lemmainte}.
\begin{flushright}
$\blacksquare$
\end{flushright}

\subsection{End of the proof}
We repeat the same procedure as in section \ref{step} but with the mollified Rankin-Selberg $L$-function instead of the Rankin-Selberg $L$-function itself. Then, \eqref{firstzero} and \eqref{rank0} become
\begin{equation}
\label{firstzerobis}
\mu_q^{h}\left(S_{k}^{p}(q,N)\right)\leq\frac{1}{N+1}\frac{1}{8b\sinh{\left(\frac{\pi c}{2b}\right)}}\left(J_{1}^{q,h}+J_{2}^{q,h}+\frac{1}{A_q^{h}[1]}A_q^{h}[I_{.}^{q}(3)]\right)
\end{equation}
and
\begin{equation}
\label{rank1}
\frac{1}{A_q^{h}[1]}A_q^{h}\left[r(.\times g)\right]\leq 2\frac{1}{8b\sinh{\left(\frac{\pi c}{2b}\right)}}\left(J_{1}^{q,h}+J_{2}^{q,h}+\frac{1}{A_q^{h}[1]}A_q^{h}[I_{.}^{q}(3)]\right).
\end{equation}
where
\begin{eqnarray*}
J_{1}^{q,h} & := & \int_{0}^{b}\cos{\left(\frac{\pi t}{2b}\right)}\log{\left(\frac{1}{A_q^{h}[1]}\mathcal{W}^{h}\left(g;\frac{-c+it}{\log{q}}\right)\right)}\mathrm{d}t, \\
J_{2}^{q,h} & := & \int_{-c}^{\left(\sigma_{0}-\frac{1}{2}\right)\log{q}}\sinh{\left(\frac{\pi(x+c)}{2b}\right)}\log{\left(\frac{1}{A_q^{h}[1]}\mathcal{W}^{h}\left(g;\frac{x+ib}{\log{q}}\right)\right)}\mathrm{d}x.
\end{eqnarray*}
We set $\tilde{b}:=2b$, $\tilde{c}:=2c$ and we choose $\sigma_{0}:=1+\frac{\log_{2}{q}}{\log{q}}$ and we assume that $\Delta$ is effective. Theorem \ref{propoE} leads to:
\begin{equation}
\label{right}
\frac{1}{A_q^{h}[1]}A_q^{h}[I_{.}^{q}(3)]\ll\left(\log{q}\right)^{\left(\frac{\pi}{\tilde{b}}-4\Delta(1-\Upsilon)\right)}q^{\frac{\pi}{2\tilde{b}}-2\Delta(1-\Upsilon)}.
\end{equation}
This is an error term under the following assumption on the height of the box:
\begin{equation*}
\tilde{b}>\frac{\pi}{4\Delta(1-\Upsilon)}.
\end{equation*}
Proposition D leads to:
\begin{equation*}
J_{1}^{q,h}=\frac{1}{2}\int_{0}^{\tilde{b}}\cos{\left(\frac{\pi t}{2\tilde{b}}\right)}\log{\left(\mathcal{V}\left(-\tilde{c},t\right)\right)}\mathrm{d}t+O_g\left(\frac{1}{q^\delta}+\frac{1}{\log{q}}\right).
\end{equation*}
Let $0<\beta<1$ be some real number. We set:
\begin{equation*}
J_{2}^{q,h}=\int_{-c}^{\log^\beta{(q)}}\cdots+\int_{\log^\beta{(q)}}^{\left(\sigma_0-\frac{1}{2}\right)\log{q}}\cdots:=J_{2,1}^{q,h}+J_{2,2}^{q,h}.
\end{equation*}
Our choice of the height of the box (so that all integrals converge), Proposition D , Remark \ref{extension} and \eqref{afsecond} entail that
\begin{eqnarray*}
J_{2,1}^{q,h} & \leq & \frac{1}{2}\int_{0}^{+\infty}\sinh{\left(\frac{\pi x}{2\tilde{b}}\right)}\log{\mathcal{V}\left(x-\tilde{c},\tilde{b}\right)}\mathrm{d}x+\frac{\exp{\left(\frac{\pi\log^{\beta}{(q)} }{\tilde{b}}\right)}}{q^\delta}+\frac{1}{\log^{1-\beta}{(q)}}, \\
J_{2,2}^{q,h} & \ll & \exp{\left(-\left(4\Delta(1-\Upsilon)-\frac{\pi}{\tilde{b}}\right)\right)\log^{\beta}(q)}
\end{eqnarray*}
and we conclude that
\begin{multline*}
\mu_q^{h}\left(S_{k}^{p}(q,N)\right)\leq\frac{1}{N+1}\frac{1}{8\tilde{b}\sinh{\left(\frac{\pi\tilde{c}}{2\tilde{b}}\right)}}\Bigg(\int_{0}^{\tilde{b}}\cos{\left(\frac{\pi t}{2\tilde{b}}\right)}\log{\left(\mathcal{V}\left(-\tilde{c},t\right)\right)}\mathrm{d}t \\
+\int_{0}^{\infty}\sinh{\left(\frac{\pi x}{2\tilde{b}}\right)}\log{\left(\mathcal{V}\left(x-\tilde{c},\tilde{b}\right)\right)}\mathrm{d}x\Bigg)+o_g(1)
\end{multline*}
and that
\begin{multline*}
\frac{1}{A_q^{h}[1]}A_q^{h}\left[r(.\times g)\right]\leq 2\frac{1}{8\tilde{b}\sinh{\left(\frac{\pi\tilde{c}}{2\tilde{b}}\right)}}\Bigg(\int_{0}^{\tilde{b}}\cos{\left(\frac{\pi t}{2\tilde{b}}\right)}\log{\left(\mathcal{V}\left(-\tilde{c},t\right)\right)}\mathrm{d}t \\
+\int_{0}^{\infty}\sinh{\left(\frac{\pi x}{2\tilde{b}}\right)}\log{\left(\mathcal{V}\left(x-\tilde{c},\tilde{b}\right)\right)}\mathrm{d}x\Bigg)+o_g(1).
\end{multline*}
We choose under $H_2(\theta)$, $P(x)=3\left(\frac{x}{\Upsilon}\right)^{2}-2\left(\frac{x}{\Upsilon}\right)^{3}$, $\tilde{b}=\frac{\pi}{4\Delta(1-\Upsilon)-10^{-10}}$, $\Delta=\Delta_{max}(\theta)-10^{-10}$ and we minimize \footnote{The program (inte.mws) is available at http://www.dms.umontreal.ca/$\sim$ricotta.} the right-hand side by a numerical choice of the remaining parameters. Under $H_2(\theta_0)$\footnote{Under $H_2(0)$, we get $\frac{3.83}{N+1}$ and $7.66$ if we choose $\Upsilon=0.45$ and $\tilde{c}=23.7$.}, the choice $\Upsilon=0.44$, $\tilde{c}=23$ gives
\begin{equation*}
\mu_q^{h}\left(S_{k}^{p}(q,N)\right)\leq\frac{4.91}{N+1}+o_g(1)
\end{equation*}
and
\begin{equation*}
\frac{1}{A_q^{h}[1]}A_q^{h}\left[r(.\times g)\right]\leq 9.82+o_g(1).
\end{equation*}
Finally, $N$ must be $4$.
\begin{flushright}
$\blacksquare$
\end{flushright}

\section{The harmonic mollified second moment near the critical point}

\label{around}

\subsection{The second harmonic twisted moment}

E. Kowalski, P. Michel and J. Vanderkam computed this moment under some sensible conditions on $D$ and $k, q$.\newline\newline
We recall here some notations used in \textbf{\cite{KoMiVa}}. For $z, s$ some complex numbers, we set
\begin{equation*}
G_{g,z}(s):=\frac{\left(4\pi^{2}\right)^{z}}{\varGamma\left(\frac{1}{2}+z+\frac{\vert k-k_g\vert}{2}\right)\varGamma\left(\frac{1}{2}+z+\frac{k+k_g}{2}-1\right)}\left(\frac{\xi\left(\frac{1}{2}+s-z\right)}{\xi\left(\frac{1}{2}\right)}\right)^{5}\frac{P_{g}(s)}{P_{g}(z)}
\end{equation*}
where $\xi(s)=s(1-s)\pi^{-\frac{s}{2}}\varGamma\left(\frac{s}{2}\right)\zeta(s)$ and $P_{g}(s)$ is an even polynomial whose coefficients are real and depend only on $k$ and $k_{g}$ chosen such that the function  $P_{g}(s)\varGamma\left(\frac{1}{2}+s+\frac{\vert k-k_g\vert}{2}\right)\varGamma\left(\frac{1}{2}+s+\frac{k+k_g}{2}-1\right)$ is analytic on $\Re{(s)}>-A$ where $A>\frac{1}{2}$. We observe that
\begin{equation}
\label{transfo}
\forall (z,s)\in\mathbb{C}^2,\quad G_{g,z}(-s)=\varepsilon_z(f\times g)G_{g,-z}(s)
\end{equation}
with
\begin{equation*}
\varepsilon_z(f\times g):=\frac{\left(4\pi^{2}\right)^{z}\varGamma\left(\frac{1}{2}-z+\frac{\vert k-k_g\vert}{2}\right)\varGamma\left(\frac{1}{2}-z+\frac{k+k_g}{2}-1\right)}{\left(4\pi^{2}\right)^{-z}\varGamma\left(\frac{1}{2}+z+\frac{\vert k-k_g\vert}{2}\right)\varGamma\left(\frac{1}{2}+z+\frac{k+k_g}{2}-1\right)}
\end{equation*}
Then, we set:
\begin{equation*}
H_{g,z}(s):=\left(4\pi^{2}\right)^{-s}\varGamma\left(\frac{1}{2}+s+\frac{\vert k-k_g\vert}{2}\right)\varGamma\left(\frac{1}{2}+s+\frac{k+k_g}{2}-1\right)G_{g,z}(s).
\end{equation*}
So, it follows that $H_{g,z}(z)=1$.\newline\newline
We need to introduce some extra notations. For any $(\alpha,\beta)=(\pm\mu,\pm\overline{\mu})$ and any natural integer $\ell\geq 1$, we set
\begin{equation*}
\mathsf{M}_{g}((\alpha,\beta);\ell):=\frac{\varphi(q)}{q\sqrt{\ell}}\, \text{res}_{u=\alpha}\,\frac{1}{2i\pi}\int_{(3)}J_{g}(u,v;(\alpha,\beta);\ell)(qD)^{u+v}\frac{\mathrm{d}v}{(u-\alpha)(v-\beta)}
\end{equation*}
with
\begin{multline*}
J_{g}(u,v;(\alpha,\beta);\ell):=H_{g,\alpha}(u)\zeta^{(D)}(1+2u)H_{g,\beta}(v)\zeta^{(qD)}(1+2v) \\
\times\nu_{g}(\ell;u,v)\frac{L(g\times g,1+u+v)}{\zeta^{D}(2(1+u+v))}
\end{multline*}
where
\begin{equation}
\label{nnnnug}
\nu_{g}(\ell;u,v):=\sum_{\delta\varepsilon=\ell}\frac{1}{\delta^{u}\varepsilon^{v}}\prod_{p\mid\delta\varepsilon}\left(\sum_{k\geq 0}\frac{\lambda_{g}(p^{k+v_{p}(\delta)})\lambda_{g}(p^{k+v_{p}(\varepsilon)})}{p^{k(1+u+v)}}\right)\left(\sum_{k\geq 0}\frac{\lambda_{g}(p^{k})\lambda_{g}(p^{k})}{p^{k(1+u+v)}}\right)^{-1}. 
\end{equation}
Finally, we set:
\begin{eqnarray*}
\varepsilon_{f\times g}(\mu,\overline{\mu}) & = & 1, \\
\varepsilon_{f\times g}(-\mu,\overline{\mu}) & = & \varepsilon_{\mu}(f\times g), \\
\varepsilon_{f\times g}(\mu,-\overline{\mu}) & = &  \varepsilon_{\overline{\mu}}(f\times g), \\
\varepsilon_{f\times g}(-\mu,-\overline{\mu}) & = & \varepsilon_{\mu}(f\times g)\varepsilon_{\overline{\mu}}(f\times g).
\end{eqnarray*}
One has:
\begin{theorem}[E. Kowalski-P. Michel-J. Vanderkam (2002)]
\label{komiva}
Let $g$  be a primitive cusp form of square-free level $D$ and trivial nebentypus. Assume that $q$ is prime, coprime with $D$. If $\vert\tau\vert\ll\frac{1}{\log{q}}$ then for any natural integer $1\leq\ell<q$,
\begin{equation*}
\left(qD\right)^{2\tau}\mathcal{M}_{g}^{h}(\mu;\ell)=\sum_{(\alpha,\beta)=(\pm\mu,\pm\overline{\mu})}\varepsilon_{f\times g}(\alpha,\beta)\mathsf{M}_{g}((\alpha,\beta);\ell)+\text{Errtwist}(q,\ell;\mu)
\end{equation*}
where
\begin{equation*}
\text{Errtwist}(q,\ell;\mu)=O_{\varepsilon,k,g}\left((q\ell)^{\varepsilon}\left(1+\vert t\vert\right)^{B}\left(\ell^{a_{1}}q^{-b_{1}}+\ell^{a_{2}}q^{-b_{2}}\right)\right)
\end{equation*}
for any $\varepsilon>0$ with $a_{1}=\frac{3}{4}, b_{1}=\frac{1}{12}$ and $a_{2}=\frac{17}{8}, b_{2}=\frac{1}{4}$.
\end{theorem}
\begin{remark}
Actually, theorem \ref{komiva} was only proved for $k<12$ so that $S_k(q)$ has no old forms. We explain in appendix \ref{kgrd} how to remove this condition using a technique of H. Iwaniec, W. Luo and P. Sarnak (\textbf{\cite{IwLuSa}}).
\end{remark}
Coming back to the notation of the introduction, we have:
\begin{equation*}
MT(\mu):=\sum_{(\alpha,\beta)=(\pm\mu,\pm\overline{\mu})}\varepsilon_{f\times g}(\alpha,\beta)\mathsf{M}_{g}((\alpha,\beta);\ell).
\end{equation*}

\subsection{The harmonic mollified second moment near $\frac{1}{2}$}

By opening the square and using multiplicative properties of Hecke eigenvalues, one gets
\begin{multline*}
\mathcal{W}^{h}(g;\mu)=(qD)^{-2\tau}\sum_{\ell_{1},\ell_{2}\geq 1}\frac{1}{\ell_{1}^{\frac{1}{2}+\mu_1}\ell_{2}^{\frac{1}{2}+\mu_2}}\sum_{d\geq 1}\frac{\varepsilon_{q}(d)}{d^{1+2\tau}}x_{d\ell_{1}}\left(g,\frac{1}{2}+\mu_1\right)x_{d\ell_{2}}\left(g,\frac{1}{2}+\mu_2\right) \\
\times(qD)^{2\tau}\mathcal{M}_{g}^{h}(\mu;\ell_{1}\ell_{2})
\end{multline*}
where one has set $\mu_1:=\mu$ and $\mu_2:=\overline{\mu}$. Our next step is to evaluate $\mathcal{W}^{h}(g;\mu)$ for $\mu$ within a distance $O\left(\frac{1}{\log{q}}\right)$ of the origin (Proposition C). We set for  $(\alpha,\beta)=(\pm\mu,\pm\overline{\mu})$:
\begin{multline}
\label{Wg}
\mathsf{W}_{g}(\alpha,\beta):=(qD)^{-2\tau}\sum_{\ell_{1},\ell_{2}\geq 1}\frac{1}{\ell_{1}^{\frac{1}{2}+\mu_{1}}\ell_{2}^{\frac{1}{2}+\mu_{2}}}\sum_{d\geq 1}\frac{\varepsilon_{q}(d)}{d^{1+\mu_{1}+\mu_{2}}} \\
\times x_{d\ell_{1}}\left(g,\frac{1}{2}+\mu_{1}\right)x_{d\ell_{2}}\left(g,\frac{1}{2}+\mu_{2}\right)\mathsf{M}_{g}((\alpha,\beta);\ell_{1}\ell_{2}).
\end{multline}
Theorem \ref{komiva} leads to:
\begin{proposition}
\label{firststep}
Let $g$ be a primitive cusp form of square-free level $D$ and trivial nebentypus. Assume that $q$ is prime, coprime with $D$. If $\vert\tau\vert\ll\frac{1}{\log{q}}$ then for any natural integer $1\leq L<\sqrt{q}$,
\begin{equation}
\mathcal{W}^{h}(g;\mu)=\sum_{(\alpha,\beta)=(\pm\mu,\pm\overline{\mu})}\varepsilon_{f\times g}(\alpha,\beta)\mathsf{W}_{g}(\alpha,\beta)+\text{Errsec}(q,L;\mu)
\end{equation}
where
\begin{equation*}
\text{\text{Errsec}}(q,L;\mu):=\!\!(qD)^{-2\tau}\!\!\!\sum_{\substack{1\leq \ell_1, \ell_2, d, \\
d\ell_1\leq L, \\
d\ell_2\leq L}}\!\!\!\!\frac{x_{d\ell_1}\left(g,\frac{1}{2}+\mu_1\right)x_{d\ell_2}\left(g,\frac{1}{2}+\mu_2\right)}{\ell_1^{\frac{1}{2}+\mu_1}\ell_2^{\frac{1}{2}+\mu_2}d^{1+\mu_1+\mu_2}}\text{Errtwist}(q,\ell_1\ell_2;\mu)
\end{equation*}
satisfies
\begin{equation}
\text{Errsec}(q,L;\mu)=O_{\epsilon,k,g}\left((qL)^{\varepsilon}(1+\vert t\vert)^{B}\left(L^{2\left(a_{1}+\frac{1}{2}\right)}q^{-b_{1}}+L^{2\left(a_{2}+\frac{1}{2}\right)}q^{-b_{2}}\right)\right)
\end{equation}
for any $\varepsilon>0$. As a consequence, $\Delta<\inf{\left(\frac{b_{1}}{4\left(a_{1}+\frac{1}{2}\right)},\frac{b_{2}}{4\left(a_{2}+\frac{1}{2}\right)}\right)}=\frac{1}{60}$ is effective.
\end{proposition}
\noindent{\textbf{Proof of proposition \ref{firststep}.}}
We only have to check the order of magnitude of the error term. We get:
\begin{multline*}
\vert \text{Errsec}(q,L;\mu)\vert\leq\sum_{d\geq 1}\frac{1}{d^{1+2\tau}}\sum_{\ell_{1},\ell_{2}\geq 1}\frac{\left\vert x_{d\ell_{1}}\left(g,\frac{1}{2}+\mu_{1}\right)\right\vert}{\ell_{1}^{\frac{1}{2}+\tau}}\frac{\left\vert x_{d\ell_{2}}\left(g,\frac{1}{2}+\mu_{2}\right)\right\vert}{\ell_{2}^{\frac{1}{2}+\tau}} \\
\times\left((\ell_1\ell_2)^{a_{1}}q^{-b_{1}}+(\ell_1\ell_2)^{a_{2}}q^{-b_{2}}\right).
\end{multline*}
As $\left\vert x_{d\ell_{i}}\left(g,\frac{1}{2}+\mu_{i}\right)\right\vert\ll_{\varepsilon}L^{\varepsilon}\sum_{d\ell_{i}=m_{1}m_{2}^{2}m_{3}}1$, $(d\ell_{i})^{-\tau}=O(1)$ and $a_{i}-\frac{1}{2}\geq 0$, we remark that
\begin{equation*}
\vert \text{Errsec}(q,L;\mu)\vert\ll_{\varepsilon}(qL)^{\varepsilon}\sum_{d\geq 1}\frac{1}{d}\sum_{1\leq \ell_{1},\ell_{2}\leq\frac{L}{d}}\left(L^{2(a_{1}-\frac{1}{2})}q^{-b_{1}}+L^{2(a_{2}-\frac{1}{2})}q^{-b_{2}}\right)
\end{equation*}
which leads to the result. \begin{flushright} $\blacksquare$ \end{flushright}
We study now the main term of the second harmonic mollified moment.
\begin{proposition}
\label{secondstep}
Let $g$ be a primitive cusp form of square-free level $D$ and trivial nebentypus. If $\vert\tau\vert\ll\frac{1}{\log{q}}$ then there exists $\delta>0$ such that
\begin{equation*}
\sum_{(\alpha,\beta)=(\pm\mu,\pm\overline{\mu})}\varepsilon_{f\times g}(\alpha,\beta)\mathsf{W}_{g}(\alpha,\beta)=\sum_{(\alpha,\beta)=(\pm\mu,\pm\overline{\mu})}\Psi(\alpha,\beta)\mathcal{V}_{(\alpha,\beta)}(\mu)+O_{\varepsilon,g}(q^{-\delta})
\end{equation*}
where for any $(\alpha,\beta)=(\pm\mu,\pm\overline{\mu})$
\begin{equation*}
\mathcal{V}_{(\alpha,\beta)}(\mu):=\sum_{\ell\geq 1}\nu_{g}(\ell;\alpha,\beta)\sum_{\ell_{1}\ell_{2}=\ell}\sum_{d\geq 1}\frac{1}{d^{1+\mu_{1}+\mu_{2}}}\frac{x_{d\ell_{1}}\left(g,\frac{1}{2}+\mu_{1}\right)}{\ell_{1}^{1+\mu_{1}}}\frac{x_{d\ell_{2}}\left(g,\frac{1}{2}+\mu_{2}\right)}{\ell_{2}^{1+\mu_{2}}}
\end{equation*}
and
\begin{multline*}
\Psi(\alpha,\beta):=\frac{\varphi(q)}{q}(qD)^{-2\tau+\alpha+\beta}\varepsilon_{f\times g}(\alpha,\beta)L(g\times g,1+\alpha+\beta) \\
\times\frac{\zeta^{(D)}(1+2\alpha)\zeta^{(qD)}(1+2\beta)}{\zeta^{(D)}(2(1+\alpha+\beta))}.
\end{multline*}
\end{proposition}
\noindent{\textbf{Proof of proposition \ref{secondstep}.}}
According to \eqref{Wg} and the integral expression of the coefficients of the mollifier (confer proposition \ref{intemolli}), one gets for $(\alpha,\beta)=(\pm\mu,\pm\overline{\mu})$:
\begin{multline*}
\mathsf{W}_{g}(\alpha,\beta)=\frac{\phi(q)}{q\left(log{L}\right)^{2}}\text{res}_{u=\alpha}\frac{1}{(2i\pi)^{3}}\int_{(3)}\int_{(3)}\int_{(3)} \\
\times m_g(u,v,s_1,s_2)\frac{\mathrm{d}s_{1}}{s_{1}^{2}}\frac{\mathrm{d}s_{2}}{s_{2}^{2}}\frac{\mathrm{d}v}{(u-\alpha)(v-\beta)}
\end{multline*}
where we have set:
\begin{multline*}
m_g(u,v,s_1,s_2):=(qD)^{u+v}H_{g,\alpha}(u)H_{g,\beta}(v)h_{g}(u,v,s_1,s_2) \\
\times L^{(1-\Upsilon)s_{1}}\left(\widehat{P^{\prime}_{L}}(s_{1})L^{\Upsilon s_{1}}-\frac{1}{1-\Upsilon}\widehat{R^{\prime}_{L^{1-\Upsilon}}}(s_{1})\right)L^{(1-\Upsilon)s_{2}}\left(\widehat{P^{\prime}_{L}}(s_{2})L^{\Upsilon s_{2}}-\frac{1}{1-\Upsilon}\widehat{R^{\prime}_{L^{1-\Upsilon}}}(s_{2})\right) \\
\times\frac{\zeta^{(D)}(1+2u)\zeta^{(qD)}(1+2v)}{\zeta^{(D)}(1+s_1+2\mu_1)\zeta^{(D)}(1+s_2+2\mu_2)\zeta^{(D)}(2(1+u+v))} \\
\times\frac{L(g\times g,1+u+v)L(g\times g,1+s_1+s_2+\mu_1+\mu_2)}{L(g\times g,1+u+s_2+\mu_2)L(g\times g,1+v+s_2+\mu_2)} \\
\times\frac{L(g\times g,1+s_1+2\mu_1)L(g\times g,1+s_2+2\mu_2)}{L(g\times g,1+u+s_1+\mu_1)L(g\times g,1+v+s_1+\mu_1)}.
\end{multline*}
Here, $h_{g}$ satisfies $h_g(u,v,s_1,s_2)=h_g(v,u,s_1,s_2)$ and defines an holomorphic function given by an absolutely convergent Euler product if $u$, $v$, $s_1$ and $s_2$ all have real part greater than some small negative real number ($-10^{-6}$ say). Thus, the pole at $u=\alpha$ is simple and so:
\begin{equation*}
\mathsf{W}_{g}(\alpha,\beta)=\frac{\phi(q)}{q\left(log{L}\right)^{2}}\frac{1}{(2i\pi)^{3}}\int_{(3)}\int_{(3)}m_g(\alpha,v,s_1,s_2)\frac{\mathrm{d}s_{1}}{s_{1}^{2}}\frac{\mathrm{d}s_{2}}{s_{2}^{2}}\frac{\mathrm{d}v}{v-\beta}.
\end{equation*}
As a function of $v$, the integrand has three simple poles at $v=\beta$, $-\alpha$ and $0$. We shift the $v$-contour to $\left(-\frac{1}{2}+\varepsilon\right)$ hitting these three poles and we remark that the remaining integral is bounded by $q^{-\delta}$ for some $\delta>0$. Thus, at an admissible cost, we have
\begin{equation*}
\sum_{(\alpha,\beta)=(\pm\mu,\pm\overline{\mu})}\!\!\!\!\!\varepsilon_{f\times g}(\alpha,\beta)\mathsf{W}_{g}(\alpha,\beta)\!\!=\!\!\sum_{(\alpha,\beta)=(\pm\mu,\pm\overline{\mu})}\!\!\!\!\!\!\varepsilon_{f\times g}(\alpha,\beta)\!\!\left(r_1(\alpha,\beta)+r_2(\alpha,\beta)+r_3(\alpha,\beta)\right)
\end{equation*}
where $r_1(\alpha,\beta)$ (respectively $r_2(\alpha,\beta)$, $r_3(\alpha,\beta)$) is the contribution of the residue at $v=\beta$ (respectively $-\alpha$, $0$) which comes from $\mathsf{W}_g(\alpha,\beta)$. We remark that
\begin{eqnarray*}
\varepsilon_{f\times g}(\alpha,\beta)r_2(\alpha,\beta) & = & -\varepsilon_{f\times g}(-\alpha,-\beta)r_2(-\alpha,-\beta), \\
\varepsilon_{f\times g}(\alpha,\beta)r_3(\alpha,\beta) & = & -\varepsilon_{f\times g}(\alpha,-\beta)r_3(\alpha,-\beta)
\end{eqnarray*}
according to \eqref{transfo}. Summing up, we get at an admissible cost
\begin{equation*}
\varepsilon\left(g;\mu\right)\mathsf{W}(g;\mu)=\sum_{(\alpha,\beta)=(\pm\mu,\pm\overline{\mu})}\varepsilon_{f\times g}(\alpha,\beta)r_1(\alpha,\beta)
\end{equation*}
which is exactly the main term in proposition \ref{secondstep}. \begin{flushright} $\blacksquare$ \end{flushright}

We set for any integers $m, n\geq 1$ and for any $(\alpha,\beta)=(\pm\mu,\pm\overline{\mu})$:
\begin{eqnarray*}
V_g(m,n;\alpha,\beta) & := & \prod_{\substack{p\in\mathcal{P} \\ p\mid\mid m \\ p\mid\mid n}}\frac{\nu_{g}(p^3;\alpha,\beta)}{\nu_{g}(p;\alpha,\beta)\nu_{g}(p^2;\alpha,\beta)}, \\
W_g(m,n;\alpha,\beta) & := & \prod_{\substack{p\in\mathcal{P} \\ p\mid\mid m \\ p\mid\mid n}}\frac{\nu_{g}(p^2;\alpha,\beta)}{\nu_{g}(p;\alpha,\beta)^2}.
\end{eqnarray*}
\begin{lemma}
\label{v}
Let $\mu\in\mathbb{C}$ and $(\alpha,\beta)=(\pm\mu,\pm\overline{\mu})$. We have
\begin{equation*}
\mathcal{V}_{(\alpha,\beta)}(\mu)=\sum_{w\geq 1}\frac{1}{w^{1+\mu_{1}+\mu_{2}}}\sum_{uv\mid w}\tau_{(\alpha,\beta)}(u,v)S_{u,v,w}(\alpha,\beta;\mu_{1})S_{u,v,w}(\alpha,\beta;\mu_{2})
\end{equation*}
where for $z\in\left\{\mu_1,\mu_2\right\}$
\begin{eqnarray*}
\tau_{(\alpha,\beta)}(u,v) & = & \frac{\mu(u)\nu_{g}(u^{2};\alpha,\beta)\nu_{g}(v^{2};\alpha,\beta)V_g(u,v;\alpha,\beta)^2}{uv}, \\
S_{u,v,w}(\alpha,\beta;z) & = & \sum_{\ell\geq 1}\frac{\nu_{g}(\ell;\alpha,\beta)V_g(\ell,v;\alpha,\beta)W_g(\ell,u;\alpha,\beta)}{\ell^{1+z}}x_{w\ell}\left(g,\frac{1}{2}+z\right).
\end{eqnarray*}
\end{lemma}
\noindent{\textbf{Proof of lemma \ref{v}.}}
One gets setting $\ell_{1}=ka$ and $\ell_{2}=kb$ with $a\wedge b=1$:
\begin{equation*}
\mathcal{V}_{(\alpha,\beta)}(\mu)=\sum_{k\geq 1}\sum_{a\wedge b=1}\frac{\nu_{g}(k^{2}ab;\alpha,\beta)}{k^{1+\mu_{1}+\mu_{2}}a^{1+\mu_{1}}b^{1+\mu_{2}}}\sum_{d\geq 1}\frac{x_{dka}\left(g,\frac{1}{2}+\mu_{1}\right)x_{dkb}\left(g,\frac{1}{2}+\mu_{2}\right)}{d^{1+\mu_{1}+\mu_{2}}}.
\end{equation*}
As $\nu_{g}$ is a multiplicative function, $ka$ and $kb$ are cube-free integers and $a\wedge b=1$, we have:
\begin{equation*}
\nu_{g}(k^{2}ab;\alpha,\beta)=\nu_{g}(k^{2};\alpha,\beta)\nu_{g}(a;\alpha,\beta)\nu_{g}(b;\alpha,\beta)V_g(a,k;\alpha,\beta)V_g(b,k;\alpha,\beta).
\end{equation*}
Hence,
\begin{multline*}
\mathcal{V}_{(\alpha,\beta)}(\mu)=\sum_{k,c,d\geq 1}\frac{\nu_{g}(k^{2};\alpha,\beta)\mu(c)}{(kcd)^{1+\mu_{1}+\mu_{2}}kc} \\
\times\sum_{a\geq 1}\frac{\nu_{g}(ac;\alpha,\beta)V_g(ac,k;\alpha,\beta)}{a^{1+\mu_{1}}}x_{dka}\left(g,\frac{1}{2}+\mu_{1}\right) \\
\times\sum_{b\geq 1}\frac{\nu_{g}(bc;\alpha,\beta)V_g(bc,k;\alpha,\beta)}{b^{1+\mu_{2}}}x_{dkb}\left(g,\frac{1}{2}+\mu_{2}\right).
\end{multline*}
Once again, we get
\begin{equation*}
\nu_{g}(ac;\alpha,\beta)=\nu_{g}(a;\alpha,\beta)\nu_{g}(c;\alpha,\beta)W_g(a,c;\alpha,\beta)
\end{equation*}
which leads to the right expression stated in the lemma. \begin{flushright} $\blacksquare$ \end{flushright}
We set for any $(\alpha,\beta)=(\pm\mu,\pm\overline{\mu})$
\begin{eqnarray*}
\mathcal{V}_{(\alpha,\beta)}^{\leq}(\mu) & := & \sum_{1\leq w\leq L^{1-\Upsilon}}\frac{1}{w^{1+\mu_{1}+\mu_{2}}}\sum_{uv\mid w}\tau_{(\alpha,\beta)}(u,v)S_{u,v,w}(\alpha,\beta;\mu_{1})S_{u,v,w}(\alpha,\beta;\mu_{2}), \\ \mathcal{V}_{(\alpha,\beta)}^{>}(\mu) & := & \sum_{L^{1-\Upsilon}<w\leq L}\frac{1}{w^{1+\mu_{1}+\mu_{2}}}\sum_{uv\mid w}\tau_{(\alpha,\beta)}(u,v)S_{u,v,w}(\alpha,\beta;\mu_{1})S_{u,v,w}(\alpha,\beta;\mu_{2})
\end{eqnarray*}
and we refer to these by the summation of respectively the short range and long range terms. 

\subsubsection{Contribution of the short range terms. \newline}

\noindent{\textbf{Treatment of $\mathrm{S_{u,v,w}(\alpha,\beta;z)}$ when $\mathrm{1\leq w\leq L^{1-\mathnormal{\Upsilon}}}$.}}
We set for any complex number $z$ and any natural integer $\ell\geq 1$
\begin{equation*}
\phi_{z}(\ell):=\sum_{\ell=\ell_{1}\ell_{2}^{2}\ell_{3}}\frac{\mu^{2}(\ell_{1}\ell_{2}\ell_{3})\mu(\ell_{1}\ell_{3})\varepsilon_{D}(\ell_{2}\ell_{3})}{\ell_{3}^{1+2z}}\lambda_{g}(\ell_{1}\ell_{3})
\end{equation*}
so that
\begin{equation*}
x_{\ell}\left(g,\frac{1}{2}+z\right)=K(g,1+2z)\sum_{\ell\geq 1}\phi_{z}(\ell)F_{L}^{\Upsilon}(\ell)K_{(\ell)}(g,1+2z)^{-1}.
\end{equation*}
We also set for any integers $u, v, w\geq 1$ with $uv\mid w$, any real $y>0$, any complex number $s$ and any polynomial $R$:
\begin{multline}
\label{Tw}
T_{u,v,w}(s;\alpha,\beta,z)=K(g,1+2z)\sum_{\ell\geq 1}\frac{\nu_{g}(\ell;\alpha,\beta)V_g(\ell,v;\alpha,\beta)W_g(\ell,u;\alpha,\beta)}{\ell^{1+s+z}} \\
\times\sum_{w\ell=\ell_{1}\ell_{2}^{2}\ell_{3}}\frac{\mu^{2}(\ell_{1}\ell_{2}\ell_{3})\mu(\ell_{1}\ell_{3})\varepsilon_{D}(\ell_{2}\ell_{3})}{\ell_{3}^{1+2z}}\lambda_{g}(\ell_{1}\ell_{3})K_{(w\ell)}(g,1+2z)^{-1},
\end{multline}
\begin{multline}
\label{TwR}
T_{u,v,w,y,R}(\alpha,\beta,z)= K(g,1+2z)\sum_{1\leq \ell\leq\frac{y}{w}}\frac{\nu_{g}(\ell;\alpha,\beta)V_g(\ell,v;\alpha,\beta)W_g(\ell,u;\alpha,\beta)}{\ell^{1+z}} \\
\times\sum_{w\ell=\ell_{1}\ell_{2}^{2}\ell_{3}}\frac{\mu^{2}(\ell_{1}\ell_{2}\ell_{3})\mu(\ell_{1}\ell_{3})\varepsilon_{D}(\ell_{2}\ell_{3})}{\ell_{3}^{1+2z}}\lambda_{g}(\ell_{1}\ell_{3})K_{(w\ell)}(g,1+2z)^{-1}R\left(\frac{\log{\left(\frac{y}{w\ell}\right)}}{\log{y}}\right).
\end{multline}
Finally, we define for any prime number $p$ and any complex number $s$
\begin{equation*}
\label{L0}
L_{p}^0(s,\alpha,\beta,z)=K_{p}(g,1+2z)+\frac{\nu_{g}(p;\alpha,\beta)\phi_z(p)}{p^{1+s+z}}+\frac{\nu_{g}(p^2;\alpha,\beta)\phi_z(p^2)}{p^{2(1+s+z)}}
\end{equation*}
and also
\begin{multline*}
L_{p}^1(u,v,w;s,\alpha,\beta,z)=K_{p}(g,1+2z)\Bigg(1+ \\
\times\frac{\nu_g(p,u;\alpha,\beta)V_g(p,v;\alpha,\beta)W_g(p,u;\alpha,\beta)\phi_z(wp)}{p^{1+s+z}\phi_z(w)}\Bigg).
\end{multline*}
We will need a nice zero-free region for $L(g\times g,1+.)$ (confer \textbf{\cite{KoMiVa}}):
\begin{lemma}
\label{zero-free}
Given $g$ as above, there exists $c_{g}>0$ depending only on $g$ such that the function $L(g\times g,1+.)$ has no zeros in the domain
\begin{equation*}
\left\{s\in\mathbb{C},\Re{(s)}\geq\frac{-c_{g}}{\log{\left(2+\vert\Im{(s)}\vert\right)}}\right\}.
\end{equation*}
Moreover, this function, its inverse and its derivatives up to any order $\alpha$ are bounded in modulus in this domain by $C_{g,\alpha,\delta}\left(1+\Im{(s)}\right)^{\delta}$ for any $\delta>0$.
\end{lemma}
This will be useful in the following lemma:
\begin{lemma}
\label{lemmaI}
Let $z\in\left\{\mu_{1},\mu_{2}\right\}$, $(\alpha,\beta)=(\pm\mu,\pm\overline{\mu})$, $y>w$ and $uv\mid w$. We have:
\begin{multline}
\label{Twbis}
T_{u,v,w}(s;\alpha,\beta,z)=\phi_{z}(w)K_{(w)}(g,1+2z)^{-1}h_{1}(u,v,w;s,\alpha,\beta,z) \\
\times\frac{L^{(q)}(Sym^{2}(g),1+2z)}{L(g\times g,1+s+z+\alpha)L(g\times g,1+s+z+\beta)}
\end{multline}
where $h_{1}$ is a holomorphic function when all the complex variables have real part greater than some small negative real number given by an absolutely convergent Euler product $h_{1}(u,v,w;s,\alpha,\beta,z):=\prod_{p\in\mathcal{P}}h_{1,p}(u,v,w;s,\alpha,\beta,z)$ with
\begin{multline}
\label{h1}
\forall p\in\mathcal{P},\quad h_{1,p}(u,v,w;s,\alpha,\beta,z)=\frac{L_{p}(g\times g,1+s+z+\alpha)L_{p}(g\times g,1+s+z+\beta)}{L_{p}^{(q)}(Sym^{2}(g),1+2z)} \\
\times\begin{cases}
L_{p}^0(s,\alpha,\beta,z) & \text{if} \quad p\nmid w, \\
L_{p}^1(u,v,w;s,\alpha,\beta,z) & \text{if} \quad p\mid\mid w, \\
K_p(g,1+2z) & \text{if} \quad p^2\mid\mid w.
\end{cases}
\end{multline}
As a consequence, if $\mu$ is a bounded complex number satisfying $\vert\tau\vert\ll\frac{1}{\log{q}}$ and $R$ a polynomial satisfying $R(0)=R^{\prime}(0)=0$ then
\begin{multline}
\label{TwRbis}
T_{u,v,w,y,R}(\alpha,\beta,z)=\Bigg\{\text{res}_{s=0}\,\phi_{z}(w)K_{(w)}(g,1+2z)^{-1}h_{1}(u,v,w;s,\alpha,\beta,z) \\
\times\frac{L^{(q)}(Sym^{2}(g),1+2z)}{sL(g\times g,1+s+z+\alpha)L(g\times g,1+s+z+\beta)}\sum_{\ell\geq 0}\frac{1}{(s\log{y})^{\ell}}R^{(\ell)}\left(\frac{\log{\left(\frac{y}{w}\right)}}{\log{y}}\right)\Bigg\} \\
+O_{g}\left(\frac{\vert\phi_z(w)\vert}{\log^2{y}}\left(\frac{y}{w}\right)^{-\left(\tau+\inf{\left(\Re(\alpha),\Re(\beta)\right)}\right)}\exp{\left(-A_{0}\sqrt{\log{\left(\frac{y}{w}\right)}}\right)}\right)
\end{multline}
for some $A_{0}>0$.
\end{lemma}
\noindent{\textbf{Proof of lemma \ref{lemmaI}.}}
The equation \eqref{Twbis} follows by comparing two Euler products. According to lemma \ref{vug} and its definition, the function $h_1$ is given by an Euler product of the following shape (everything was made for and the Ramanujan-Petersson bound for Hecke eigenvalues of $g$ are available)
\begin{equation}
\label{shapeh1}
\forall p \in\mathcal{P},\quad h_{1,p}(u,v,w;s,\alpha,\beta,z)=1+O_{u,v,w}\left(\sum_{i\in I}\frac{1}{p^{2+\Re{(a_is+b_i\alpha+c_i\beta+d_iz)}}}\right)
\end{equation}
for some finite index set $I$ and some integers $a_i$, $b_i$, $c_i$ and $d_i$. Thus, if all the complex variables have some slightly negative real parts such that
\begin{equation*}
\forall i\in I,\quad\Re{(a_is+b_i\alpha+c_i\beta+d_iz)}\geq -1+\delta
\end{equation*}
for some fixed $\delta>0$ then this Euler product absolutely converges and defines a holomorphic function. To get \eqref{TwRbis}, we use the Taylor expansion of $R$:
\begin{multline*}
T_{u,v,w,y,R}(\alpha,\beta,z)=\sum_{j\geq 2}\frac{R^{(j)}(0)}{(\log{y})^{j}}\frac{1}{2i\pi}\int_{(3)}K_{(w)}(g,1+2z)^{-1}h_{1}(u,v,w;s,\alpha,\beta,z) \\
\times\phi_{z}(w)\frac{L^{(q)}(Sym^{2}(g),1+2z)}{L(g\times g,1+s+z+\alpha)L(g\times g,1+s+z+\beta)}\left(\frac{y}{w}\right)^{s}\frac{\mathrm{d}s}{s^{j+1}}.
\end{multline*}
According to the assumptions on $\mu$, we can find $F_1>0$ such that $\Re{(z+\alpha)}$ and $\Re{(z+\beta)}\geq\frac{-F_1}{\log{y}}$. We move the integral to the line $\Re{(s)}=\frac{F_1+1}{\log{\left(\frac{y}{w}\right)}}$ without crossing any pole and then we cut the integral at the segment $\left[\frac{F_1+1}{\log{\left(\frac{y}{w}\right)}}-iT,\frac{F_1+1}{\log{\left(\frac{y}{w}\right)}}+iT\right]$ at an admissible cost of $O\left(\vert\phi_z(w)\vert\log{\left(\frac{y}{w}\right)}^2T^{-2}\right)$ where $T:=\exp{\left(\sqrt{\log{\left(\frac{y}{w}\right)}}\right)}$. We move the previous line segment to
\begin{equation*}
\left[-\inf{(\Re{(z+\alpha)},\Re{(z+\beta)})}-\frac{F_2}{\log{T}}-iT,-\inf{(\Re{(z+\alpha)},\Re{(z+\beta)})}-\frac{F_2}{\log{T}}+iT\right]
\end{equation*}
where $F_2>0$ is chosen such that this line segment is included in a free-zero area for $L(g\times g,1+.+z+\alpha)L(g\times g,1+.+z+\beta)$ given by lemma \ref{zero-free}. We cross a multiple pole at $s=0$ whose residue is precisely the main term in \eqref{TwRbis}. The remaining integrals contribute as
\begin{equation*}
O_{g}\left(\vert\phi_z(w)\vert\left(\log{\left(\frac{Ty}{w}\right)}\right)^2\left(T^{-2}+\left(\frac{y}{w}\right)^{-\inf{(\Re{(z+\alpha)},\Re{(z+\beta)})}-\frac{F_2}{\log{T}}}\right)\right).
\end{equation*}
\begin{flushright}
$\blacksquare$
\end{flushright}
This leads directly to:
\begin{proposition}
\label{Swsmall}
Let $z\in\left\{\mu_{1},\mu_{2}\right\}$, $(\alpha,\beta)=(\pm\mu,\pm\overline{\mu})$ and $u,v\geq 1$. If $\mu$ is a bounded complex number satisfying $\vert\tau\vert\ll\frac{1}{\log{q}}$ and $1\leq w\leq L^{1-\Upsilon}$ then
\begin{multline}
\label{Swsmallbis}
S_{u,v,w}(\alpha,\beta,z)=\delta_{(z+\alpha)(z+\beta)\neq 0}\phi_{z}(w)K_{(w)}(g,1+2z)^{-1} \\
\times h_{1}(u,v,w;0,\alpha,\beta,z)\frac{L^{(q)}(Sym^{2}(g),1+2z)}{L(g\times g,1+z+\alpha)L(g\times g,1+z+\beta)} \\
+O_{g}\left(\frac{\vert\phi_z(w)\vert}{\log^2{q}}\left(\frac{L^{1-\Upsilon}}{w}\right)^{-\left(\tau+\inf{\left(\Re(\alpha),\Re(\beta)\right)}\right)}\exp{\left(-A_{0}\sqrt{\log{\left(\frac{L^{1-\Upsilon}}{w}\right)}}\right)}\right).
\end{multline}
\end{proposition}
\noindent{\textbf{Proof of proposition \ref{Swsmall}.}}
Let $Q(X):=1-P(\Upsilon+(1-\Upsilon)X)$. We remark that:
\begin{equation*}
S_{u,v,w}(\alpha,\beta,z)=T_{u,v,w,L,P}(\alpha,\beta,z)+T_{u,v,w,L^{1-\Upsilon},Q}(\alpha,\beta,z)
\end{equation*}
When applying lemma \ref{lemmaI} twice, the reader may remark that the only contribution comes from the values of $P$ and $Q$ and that the other main terms coming from the values of the derivatives of $P$ and $Q$ cancel each other; this concludes the proof. \begin{flushright} $\blacksquare$ \end{flushright}
\noindent{\textbf{Treatment of $\mathrm{\mathcal{V}_{(\alpha,\beta)}^{\leq}(\mu)}$.}}
We set:
\begin{multline}
\label{tildeL}
L_{p}(s,\alpha,\beta)=\left(\prod_{z\in\left\{\mu_1,\mu_2\right\}}L_p^0(0,\alpha,\beta,z)\right) \\
+p^{-(1+s)}\Bigg\{\left(\prod_{z\in\left\{\mu_1,\mu_2\right\}}\frac{\phi_z(p)L_p^1(1,1,p;0,\alpha,\beta,z)}{K_p(g,1+2z)}\right) \\
+\nu_g(p^2;\alpha,\beta)p^{-1}\left(\prod_{z\in\left\{\mu_1,\mu_2\right\}}\frac{\phi_z(p)L_p^1(1,p,p;0,\alpha,\beta,z)}{K_p(g,1+2z)}\right) \\
-\nu_g(p;\alpha,\beta)^2p^{-1}\left(\prod_{z\in\left\{\mu_1,\mu_2\right\}}\frac{\phi_z(p)L_p^1(p,1,p;0,\alpha,\beta,z)}{K_p(g,1+2z)}\right)\Bigg\} \\
+\left(\prod_{z\in\left\{\mu_1,\mu_2\right\}}\phi_z(p^2)\right)p^{-2(1+s)}\Bigg\{1+\nu_g(p^2;\alpha,\beta)p^{-1} \\
-\nu_g(p;\alpha,\beta)^2p^{-1}-\nu_g(p;\alpha,\beta)^2\nu_g(p^2;\alpha,\beta)V_g(p,p;\alpha,\beta)^2p^{-2} \\
+\nu_g(p^4;\alpha,\beta)p^{-2}\Bigg\}.
\end{multline}
\begin{lemma}
\label{lemmaII}
Let $\mu\in\mathbb{C}$, $(\alpha,\beta)=(\pm\mu,\pm\overline{\mu})$ and $z\in\left\{\mu_1,\mu_2\right\}$.
\begin{multline}
\label{lemmaII1}
\sum_{w\geq 1}\frac{1}{w^{1+s}}\sum_{uv\mid w}\tau_{(\alpha,\beta)}(u,v) \\
\times\prod_{z\in\left\{\mu_{1},\mu_{2}\right\}}\left(\frac{\phi_{z}(w)}{K_{(w)}(g,1+2z)}h_{1}(u,v,w;0,\alpha,\beta,z)\right)= \\
L(g\times g,1+s)h_{2}(s,\alpha,\beta)
\end{multline}
where $h_{2}$ is a holomorphic function when all the complex variables have some real part greater than some small negative real number given by an absolutely convergent Euler product $h_{2}(s,\alpha,\beta):=\prod_{p\in\mathcal{P}}h_{2,p}(s,\alpha,\beta)$ with
\begin{multline}
\label{h2bis}
\forall p\in\mathcal{P},\quad h_{2,p}(s,\alpha,\beta)=\prod_{z\in\left\{\mu_1,\mu_2\right\}}\left(\frac{L_p(g\times g,1+z+\alpha)L_p(g\times g,1+z+\beta)}{L_p^{(q)}(\text{Sym}^2(g),1+2z)}\right) \\
\times L_p(g\times g,1+s){L}_{p}(s,\alpha,\beta).
\end{multline}
As a consequence, if $\vert\tau\vert\ll\frac{1}{\log{q}}$ then
\begin{multline}
\label{lemmaII2}
\sum_{1\leq w\leq x}\frac{1}{w^{1+\mu_1+\mu_2}}\sum_{uv\mid w}\tau_{(\alpha,\beta)}(u,v) \\
\times\prod_{z\in\left\{\mu_{1},\mu_{2}\right\}}\left(\frac{\phi_{z}(w)}{K_{(w)}(g,1+2z)}h_{1}(u,v,w;0,\alpha,\beta,z)\right)= \\
L(g\times g,1+\mu_{1}+\mu_{2})h_{2}(\mu_{1}+\mu_{2},\alpha,\beta)\left(1-x^{-2\tau}\right)+O_{g}\left(x^{-2\tau}\right).
\end{multline}
\end{lemma}
\noindent{\textbf{Proof of lemma \ref{lemmaII}.}}
The first part \eqref{lemmaII1} comes from a computation of Euler products. Once again, the shape of the Euler product which defines $h_2$ is
\begin{equation}
\label{shapeh2}
\forall p \in\mathcal{P},\quad h_{2,p}(u,v,w;s,\alpha,\beta,z)=1+O_{u,v,w}\left(\sum_{i\in I}\frac{1}{p^{2+\Re{(a_is+b_i\alpha+c_i\beta+d_iz)}}}\right)
\end{equation}
for some finite index set $I$ and some integers $a_i$, $b_i$, $c_i$ and $d_i$. Thus, for exactly the same reasons as in the proof of lemma \ref{lemmaI} (confer \eqref{shapeh1}), this Euler product is absolutely convergent when all the complex variables have some slightly negative real parts. To get \eqref{lemmaII2}, according to explicit Perron's formula, our sum equals
\begin{equation*}
\frac{1}{2i\pi}\int_{A-iT}^{A+iT}L(g\times g,1+s+2\tau)h_{2}(s+2\tau,\alpha,\beta)x^{s}\frac{\mathrm{d}s}{s}+O\left(\frac{x^{A-2\tau}}{T}\right)
\end{equation*}
where $A>-2\tau$ and $T>0$ will be chosen later. We shift the contour to $\Re{(s)}=-A$ hitting some poles at $s=0$ and $s=-2\tau$. The remaining integrals contribute as $O_{g}\left(\frac{x^{A-2\tau}}{T}+x^{-A-2\tau}T\right)$. We choose $T=x^{A}$ in order to justify the error term in \eqref{lemmaII2}. The residues of the crossed poles are
\begin{multline*}
L(g\times g,1+2\tau)h_{2}(2\tau,\alpha,\beta)-\frac{\text{res}_{s=1}L(g\times g,s)}{2\tau}h_{2}(0,\alpha,\beta)x^{-2\tau}= \\
L(g\times g,1+2\tau)h_{2}(2\tau,\alpha,\beta)\left(1-x^{-2\tau}\right)+O_{g}\left(x^{-2\tau}\right).
\end{multline*}
\begin{flushright}
$\blacksquare$
\end{flushright}
In the following proposition, we estimate $\mathcal{V}_{(\alpha,\beta)}^{\leq}(\mu)$.
\begin{proposition}
\label{Vsmall}
Let $\mu\in\mathbb{C}$ and $(\alpha,\beta)=(\pm\mu,\pm\overline{\mu})$. If $\vert\mu\vert\ll\frac{1}{\log{q}}$ then
\begin{multline*}
\mathcal{V}_{(\alpha,\beta)}^{\leq}(\mu)=\delta_{(\alpha,\beta)=(\mu,\overline{\mu})}\,\,h_{2}(\mu_{1}+\mu_{2},\mu,\overline{\mu})L(g\times g,1+2\tau)\left(1-L^{-2\tau(1-\Upsilon)}\right) \\
\times\left(\prod_{z\in\left\{\mu_{1},\mu_{2}\right\}}\frac{L^{(q)}(\text{Sym}^{2}(g),1+2z)}{L(g\times g,1+\mu+z)L(g\times g,1+\overline{\mu}+z)}\right) \\
+O_{g}\left(\frac{1}{\log^4{(q)}}\left(L^{-(1-\Upsilon)\left(\tau+\inf{\left(\Re{(\alpha)},\Re{(\beta)}\right)}\right)}+L^{-2(1-\Upsilon)\left(\tau+\inf{\left(\Re{(\alpha)},\Re{(\beta)}\right)}\right)}\right)\right) \\
+\delta_{(\alpha,\beta)=(\mu,\overline{\mu})}O_{g}\left(\frac{1}{\log^4{(q)}}L^{-2\tau(1-\Upsilon)}\right).
\end{multline*}
\end{proposition}
\noindent{\textbf{Proof of proposition \ref{Vsmall}.}}
Since $(z+\alpha)(z+\beta)\neq 0$:
\begin{equation*}
\left(L(g\times g,1+z+\alpha)L(g\times g,1+z+\beta)\right)^{-1}\ll_g\left(\log{q}\right)^{-2},
\end{equation*}
the proposition follows from proposition \ref{Swsmall} and lemma \ref{v}.
\begin{flushright}
$\blacksquare$
\end{flushright}
\noindent{\textbf{Treatment of the short range terms.}}
We sum up the estimate of the short range terms in the following theorem:
\begin{theorem}
\label{Wsmall}
Let $\mu\in \mathbb{C}$. If  $\frac{\varepsilon_0}{\log{q}}\leq\vert\mu\vert\ll\frac{1}{\log{q}}$ for some absolute constant $\varepsilon_0>0$ then
\begin{multline}
\label{Wsmallbis}
\sum_{(\alpha,\beta)=(\pm\mu,\pm\overline{\mu})}\Psi(\alpha,\beta)\mathcal{V}_{(\alpha,\beta)}^{\leq}(\mu)=1-L^{-2\tau(1-\Upsilon)} \\
+O_{g}\left(\frac{1}{q^\delta}+\frac{1}{\log{q}}\right.\begin{cases}
L^{-2\tau(1-\Upsilon)} & \text{ if } \tau:=\Re{(\mu)}\geq 0, \\
q^{-2\tau}L^{-4\tau(1-\Upsilon)} & \text{ otherwise}
\end{cases}
\Bigg)
\end{multline}
for some $\delta>0$.
\end{theorem}
\noindent{\textbf{Proof of theorem \ref{Wsmall}. }}As $\Psi(\alpha,\beta)\ll\log^3{(q)} \, \,q^{-2\tau+\alpha+\beta}$ (the worst case being $(\alpha,\beta)=(\mu,\overline{\mu})$), proposition \ref{Vsmall} implies that:
\begin{multline*}
\sum_{(\alpha,\beta)=(\pm\mu,\pm\overline{\mu})}\Psi(\alpha,\beta)\mathcal{V}_{(\alpha,\beta)}^{\leq}(\mu)=\Psi(\mu,\overline{\mu})h_{2}(\mu_{1}+\mu_{2},\mu,\overline{\mu})L(g\times g,1+\mu_{1}+\mu_{2}) \\
\times\left(1-L^{-2\tau(1-\Upsilon)}\right)\left(\prod_{z\in\left\{\mu_{1},\mu_{2}\right\}}\frac{L^{(q)}(Sym^{2}(g),1+2z)}{L(g\times g,1+\mu+z)L(g\times g,1+\overline{\mu}+z)}\right) \\
+O_{g}\left(\frac{1}{\log{q}}\right.\begin{cases}
L^{-2\tau(1-\Upsilon)} & \text{ if } \tau\geq 0, \\
q^{-2\tau}L^{-4\tau(1-\Upsilon)} & \text{ otherwise.}
\end{cases}
\Bigg)
\end{multline*}
The main term of the previous equality equals:
\begin{equation*}
\frac{\varphi(q)}{q}\frac{\zeta_q(1+2\overline{\mu})}{L_q(\text{Sym}^2g,1+2\mu)L_q(\text{Sym}^2g,1+2\overline{\mu})}\frac{h_{2}(2\tau,\mu,\overline{\mu})}{\zeta^{(D)}(2(1+2\tau))}\left(1-L^{-2\tau(1-\Upsilon)}\right).
\end{equation*}
According to proposition \ref{arithconst}, we have
\begin{equation}
\label{calcul1}
\frac{h_{2}(2\tau,\mu,\overline{\mu})}{\zeta^{(D)}(2(1+2\tau))}=1+O_g\left(\frac{1}{q^\delta}\right)
\end{equation}
for some $\delta>0$.
\begin{flushright}
$\blacksquare$
\end{flushright}
\begin{remark}
Equation \eqref{calcul1} is the result of tedious but elementary computations which are carried out in Appendix D. One may find it surprising that this apparently rather complicated Euler product turns out to a have very simple expression (which in fact is crucial for the method to work). This however is a consequence of our choice of mollifier. It is very plausible that a more conceptual explanation of this phenomenon can be gotten from the random matrix model for the family $\mathcal{F}$ of Rankin-Selberg $L$-functions and the vertical Sato-Tate laws satisfied by the Hecke eigenvalues of modular forms. For this, we refer to the recent work of J.B. Conrey, D. Farmer, J. Keating, M. Rubinstein and N. Snaith (\textbf{\cite{CoFaKeRuSn}}) who formulate very precise conjectures for the moments of central value for many families of $L$-functions (although not for our peculiar family, which certainly can be investigated along the same lines) and the talk of C. Hughes at the Newton Institute on amplified and mollified moments of families of $L$-functions (\textbf{\cite{Hu}}).
\end{remark}


\subsubsection{Contribution of the long range terms. \newline}
\noindent{\textbf{Treatment of $\mathrm{\mathcal{V}_{(\alpha,\beta)}^{>}(\mu)}$.}}
Arguing along the same lines, we obtain
\begin{proposition}
\label{Vlarge}
Let $\mu\in\mathbb{C}$ and $(\alpha,\beta)=(\pm\mu,\pm\overline{\mu})$. If $\vert\mu\vert\ll\frac{1}{\log{q}}$ then:
\begin{multline}
\label{Vlargebis}
\mathcal{V}_{(\alpha,\beta)}^{>}(\mu)=h_{2}(\mu_{1}+\mu_{2},\alpha,\beta)(\mu_1+\mu_2)L(g\times g,1+\mu_1+\mu_2) \\
\times\left(\text{res}_{s=1}L(g\times g,s)\right)^{-4}\left(\prod_{z\in\left\{\mu_{1},\mu_{2}\right\}}L^{(q)}(Sym^{2}(g),1+2z)\right) I_{\alpha,\beta}(L,\Upsilon,P;\mu) \\
+O_{g}\left(\frac{1}{\log^4{(q)}}\left(L^{-2\left(\tau+\inf{\left(\Re{(\alpha)},\Re{(\beta)}\right)}\right)}+L^{-\left( \tau+\inf {\left(\Re{(\alpha)},\Re{(\beta)} \right)} \right)}+L^{-2\tau}+L^{-2\tau(1-\Upsilon)}\right)\right)
\end{multline}
with
\begin{multline}
\label{inte}
I_{\alpha,\beta}(L,\Upsilon,P;\mu):=\log{L}\int_{0}^{\Upsilon}L^{-2\tau(1-x)} \\
\times\left(\prod_{z\in\left\{\mu_{1},\mu_{2}\right\}}\left((z+\alpha)(z+\beta)P(x)+\frac{(2z+\alpha+\beta)}{\log{L}}P^{\prime}(x)+\frac{1}{\log^2{(L)}}P^{\prime\prime}(x)\right)\right)\mathrm{d}x.
\end{multline}
\end{proposition}
\noindent{\textbf{Treatment of the long range terms.}}
Firstly, we compute an expression for the previous integrals $I_{\alpha,\beta}(L,\Upsilon,P;\mu)$ which are obtained by some integration by parts knowing that $P(0)=P^{\prime}(0)=P^{\prime}(\Upsilon)=0$ and $P(\Upsilon)=1$. The results are given in the following lemma: 
\begin{lemma}
\label{lemmaIIII}
Let $(\alpha,\beta)=(\pm\mu,\pm\overline{\mu})$.
\begin{multline*}
I_{\alpha,\beta}(L,\Upsilon,P;\mu)=\frac{4\mu\overline{\mu}}{\log{L}}\int_{0}^{\Upsilon}L^{-2\tau(1-x)}\left\vert P^{\prime}(x)+\frac{P^{\prime\prime}(x)}{2\mu\log{L}}\right\vert^{2}\mathrm{d}x \\
+\delta_{(\alpha,\beta)=(\mu,\overline{\mu})}8\mu\overline{\mu}\tau L^{-2\tau(1-\Upsilon)}.
\end{multline*}
\end{lemma}
We sum up the contribution of the long range terms in the following theorem:
\begin{theorem}
\label{Wlarge}
Let $\mu\in\mathbb{C}$. If $\frac{\varepsilon_0}{\log{q}}\leq\vert\mu\vert\ll\frac{1}{\log{q}}$ for some $\varepsilon_0>0$ then
\begin{multline}
\label{Wlargebis}
\mathcal{W}^{h}_{>}(\mu)=\left(\frac{q^{2\tau}-q^{-2\tau}}{2\tau\log{L}}-\frac{q^{2\delta}-q^{-2\delta}}{2\delta\log{L}}\right)\int_{0}^{\Upsilon}L^{-2\tau(1-x)}\left\vert P^{\prime}(x)+\frac{P^{\prime\prime}(x)}{2\mu\log{L}}\right\vert^{2}\mathrm{d}x \\
+L^{-2\tau(1-\Upsilon)}+O_{g}\left(\frac{1}{q^\delta}+\frac{1}{\log{q}}\right.\begin{cases}
L^{-2\tau(1-\Upsilon)} & \text{ if } \tau:=\Re{(\mu)}\geq 0, \\
q^{-4\tau}L^{-4\tau} & \text{ otherwise}
\end{cases}
\Bigg)
\end{multline}
for some $\delta>0$.
\end{theorem}
\noindent{\textbf{Proof of theorem \ref{Wlarge}. }}As $\Psi(\alpha,\beta)\ll\log^3{(q)} \, \,q^{-2\tau+\alpha+\beta}$, proposition \ref{Vlarge} implies that
\begin{multline*}
\sum_{(\alpha,\beta)=(\pm\mu,\pm\overline{\mu})}\Psi(\alpha,\beta)\mathcal{V}_{(\alpha,\beta)}^{>}(\mu)=\sum_{(\alpha,\beta)=(\pm\mu,\pm\overline{\mu})}\widetilde{\Psi}(\alpha,\beta)q^{-2\tau+\alpha+\beta}I_{\alpha,\beta}(L,\Upsilon,P;\mu) \\
+O_{g}\left(\frac{1}{\log{q}}\right.\begin{cases}
L^{-2\tau(1-\Upsilon)} & \text{ if } \tau\geq 0, \\
q^{-4\tau}L^{-4\tau} & \text{ otherwise}
\end{cases}
\Bigg)
\end{multline*}
with
\begin{multline*}
\widetilde{\Psi}(\alpha,\beta):=\left(\text{res}_{s=1}L(g\times g,s)\right)^{-4}q^{2\tau-(\alpha+\beta)}f(\alpha,\beta)2\tau L(g\times g,1+2\tau)h_2(2\tau,\alpha,\beta) \\
\left(\prod_{z\in\left\{\mu_1,\mu_2\right\}}L^{(q)}(\text{Sym}^2(g),1+2z)\right).
\end{multline*}
Moreover, $\widetilde{\Psi}(\alpha,\beta)=\frac{1}{4\alpha\beta(\alpha+\beta)}\frac{h_2(0,0,0)}{\zeta^{(D)}(2)}+O_g\left(\log^3{(q)}\right)$ and $I_{\alpha,\beta}(L,\Upsilon,P;\mu)\ll\frac{L^{-2\tau(1-\Upsilon)}}{\log^4{(q)}}$. According to proposition \ref{arithconst}, we have
\begin{equation}
\label{calcul2}
\frac{h_{2}(0,0,0)}{\zeta^{(D)}(2)}=1+O_g\left(\frac{1}{q^\delta}\right)
\end{equation}
for some $\delta>0$. Thus, the contribution of the long range terms is
\begin{equation*}
\sum_{(\alpha,\beta)=(\pm\mu,\pm\overline{\mu})}\frac{q^{-2\tau+\alpha+\beta}}{4\alpha\beta(\alpha+\beta)}I_{\alpha,\beta}(L,\Upsilon,P;\mu)
\end{equation*}
which is exactly the main term of \eqref{Wlargebis} according to lemma \ref{lemmaIIII}.
\begin{flushright}
$\blacksquare$
\end{flushright}


\section{Averaged shifted convolution problems}

\label{scp}
This section is the central part of this paper. We give here a way of solving shifted convolution problems on average.

\subsection{Introduction and first result}

Let $\Psi:\mathbb{R}_+^*\times\mathbb{R}_+^*\rightarrow\mathbb{R}$ be a smooth compactly supported function
\begin{equation}
\label{Hsupport}
\text{Supp}(\Psi)\subset\left[Z,2Z\right]\times\left[\frac{Y}{2},2Y\right]
\end{equation}
for some real numbers $Z, Y>0$ satisfying
\begin{equation}
\label{Hcontrole}
\exists P>0, \forall (\alpha,\beta)\in\mathbb{N}^2,\forall (z,y)\in\left(\mathbb{R}_+^*\right)^2, \quad z^\alpha y^\beta\frac{\partial^{\alpha+\beta}\Psi}{\partial z^{\alpha}\partial y^{\beta}}(z,y)\ll_{\alpha, \beta}P^{\alpha+\beta}.
\end{equation}
Let $a_1, a_2\geq 1$ with $a_1a_2<q$ be some natural integers. One considers the shifted convolution problem
\begin{equation*}
\forall h\in\mathbb{Z}^*, \quad S_h(\Psi,g;a_1,a_2):=\sum_{a_1m-a_2n=h}\lambda_g(m)\lambda_g(n)\Psi(a_1m,a_2n)
\end{equation*}
and the shifted convolution problem on average
\begin{equation*}
\Sigma_r(\Psi,g;a_1,a_2):=\sum_{\substack{h\in\mathbb{Z}^* \\
h\equiv 0\mod r}}S_h(\Psi,g;a_1,a_2)
\end{equation*}
for any natural integer $r\geq 1$. Note that the $h$-sum is of length $\sup{\left(Z,Y\right)}$. Solving the shifted convolution problem (respectively the shifted convolution problem on average) consists in finding a non-trivial bound for $S_h(\Psi,g;a_1,a_2)$ (respectively $\Sigma_r(\Psi,g;a_1,a_2)$). The $\delta$-symbol method of W. Duke, J. Friedlander and H. Iwaniec (confer \textbf{\cite{DuFrIw}} and \textbf{\cite{KoMiVa}}) leads to:
\begin{theorem}
\label{delta}
Let $h\in\mathbb{Z}^*$ and $r\in\mathbb{N}^*$. If $a_1\wedge a_2=1$ and $\Psi$ satisfies \eqref{Hsupport} and \eqref{Hcontrole} then
\begin{equation*}
S_h(\Psi,g;a_1,a_2)\ll_{\varepsilon ,g}P^{\frac{5}{4}}(Z+Y)^{\frac{1}{4}}(YZ)^{\frac{1}{4}+\varepsilon}
\end{equation*}
for any $\varepsilon>0$. Thus,
\begin{equation*}
\Sigma_r(\Psi,g;a_1,a_2)\ll_{\varepsilon ,g}P^{\frac{5}{4}}(Z+Y)^{\frac{1}{4}}(YZ)^{\frac{1}{4}+\varepsilon}\frac{\sup{(Z,Y)}}{r}
\end{equation*}
for any $\varepsilon>0$.
\end{theorem}

\subsection{The spectral method on average}

For some background and notations about Maass forms we refer to appendix \ref{Mas}. All is based on the analytic properties of the following Dirichlet series (confer \cite{Sa} and \cite{Mi1})
\begin{equation*}
D_{h}(g,a_1,a_2;s):=\sum_{a_1m-a_2n=h}\lambda_{g}(m)\lambda_{g}(n)\left(\frac{\sqrt{a_1a_2mn}}{a_1m+a_2n}\right)^{k_{g}-1}(a_1m+a_2n)^{-s}
\end{equation*}
which is linked to our problem by Mellin's inversion formula
\begin{equation}
\label{link}
S_{h}(g,a_1,a_2)=\frac{1}{2i\pi}\int_{(2)}D_{h}(g,a_1,a_2;s)\widehat{\Psi}(h,s)\,\,\mathrm{d}s
\end{equation}
where
\begin{equation}
\label{Fhat}
\widehat{\Psi}(h,s)=\int_{\sup{\left(\vert h\vert,h+Na_2\right)}}^{h+4Na_2}\Psi\left(\frac{u+h}{2},\frac{u-h}{2}\right)\left(4+\frac{2h}{u-h}-\frac{2h}{u+h}\right)^{\frac{k_{g}-1}{2}}u^{s}\frac{\mathrm{d}u}{u}.
\end{equation}
Note that $\widehat{\Psi}(h,s)=0$ if $\vert h\vert\gg\sup{(Z,Y)}$. The spectral method consists in getting a non-trivial individual estimate of $D_{h}(g,a_1,a_2;s)$ whereas the spectral method on average takes care of the extra average over $h$.
\begin{lemma}
\label{FFhat}
If $\,\Psi$ satisfies \eqref{Hsupport} and \eqref{Hcontrole} then:
\begin{equation*}
\widehat{\Psi}(h,s)\ll_{\eta}\left(\frac{\sup{(Z,Y)}}{\inf{(Z,Y)}}\right)^{\frac{k_{g}-1}{2}+\eta -1}\sup{(Z,Y)}^{\Re{(s)}}\frac{P^\eta}{\vert s\vert^{\eta}}
\end{equation*}
for any natural integer $\eta$.
\end{lemma}
\noindent{\textbf{Proof of lemma \ref{FFhat}.}}
According to the support properties of $\Psi$ and by $\eta$ integration by parts, we have
\begin{equation*}
\widehat{\Psi}(h,s)=\int_{\pm h+O(\inf{(Z,Y)})}^{\pm h+O(\inf{(Z,Y)})}\frac{u^{s+\eta-1}}{s(s+1)\cdots (s+\eta-1)}\gamma^{(\eta)}(u)\mathrm{d}u
\end{equation*}
with $\gamma(u)=\Psi\left(\frac{u+h}{2},\frac{u-h}{2}\right)\left(4+\frac{2h}{u-h}-\frac{2h}{u+h}\right)^{\frac{k_{g}-1}{2}}$. One shows with \eqref{Hcontrole} that
\begin{equation*}
\gamma^{(\eta)}(u)\ll_{\eta}\left(\frac{\sup{(Z,Y)}}{\inf{(Z,Y)}}\right)^{\frac{k_{g}-1}{2}}\frac{1}{\inf{(Z,Y)}^{\eta}}P^\eta
\end{equation*}
which is enough for the proof.
\begin{flushright}
$\blacksquare$
\end{flushright}
One defines the following Maass cusp forms of level $Da_1a_2$, weight $0$ and trivial nebentypus
\begin{equation*}
V(z):=(a_1y)^{\frac{k_g}{2}}g(a_1z)(a_2y)^{\frac{k_g}{2}}\overline{g(a_2z)}
\end{equation*}
and
\begin{equation*}
U_h(z,s):=\sum_{\gamma\, \in \, \left(\varGamma_0(Da_1a_2)\right)_\infty\backslash\varGamma_0(Da_1a_2)}\left(\Im{(\gamma.z)}\right)^s e(-h\Re{(\gamma.z)}).
\end{equation*}
A straightforward computation gives:
\begin{equation*}
D_{h}(g,a_1,a_2;s)=\frac{(2\pi)^{s+k_g-1}}{\Gamma(s+k_g-1)\sqrt{a_1a_2}}(U_h(.,s),\overline{V}).
\end{equation*}
Let $\beta:=\left(u_j\right)_{j\geq 1}$ a Hecke eigenbasis of $\mathcal{C}_0(Da_1a_2)$ satisfying $\left(\Delta_0+\lambda_j\right)\,\, u_j=0$ with $\lambda_j:=\frac{1}{4}+r_j^2$ and made of eigenfunctions of the reflexion operator namely: $\forall n\in\mathbb{Z}^*, \rho_j(-n)=\varepsilon_j\rho_j(n)$ for some $\varepsilon_j\in\left\{\pm 1\right\}$. Parseval's equality leads to:
\begin{multline}
\label{parseval}
D_{h}(g,a_1,a_2;s)=\frac{(2\pi)^{s+k_g-1}}{\Gamma(s+k_g-1)\sqrt{a_1a_2}} \\
\Bigg\{\sum_{j\geq 1}\frac{\sqrt{\vert h\vert}\,\overline{\rho_j(-h)}}{2\pi^{s-\frac{1}{2}}\vert h\vert^{s-\frac{1}{2}}}\Gamma\left(\frac{s-\frac{1}{2}+ir_j}{2}\right)\Gamma\left(\frac{s-\frac{1}{2}-ir_j}{2}\right)\left( u_j,\overline{V}\right) \\
+\frac{1}{4\pi}\sum_{\kappa\in \text{Cusp}\left(\varGamma_0(Da_1a_2)\right)}\int_{-\infty}^{\infty}\frac{\sqrt{\vert h\vert}\,\overline{\rho_{\kappa}(-h,t)}}{2\pi^{s-\frac{1}{2}}\vert h\vert^{s-\frac{1}{2}}}\Gamma\left(\frac{s-\frac{1}{2}+it}{2}\right)\Gamma\left(\frac{s-\frac{1}{2}-it}{2}\right) \\
\left( E_\kappa\left(.,\frac{1}{2}+it\right),\overline{V}\right)\mathrm{d}t\Bigg\}.
\end{multline}
The hypothesis $H_2(\theta)$ described in the introduction is a very natural one as it allows us to control the size of the discrete part of the right-hand side of \eqref{parseval}. In fact, H. Iwaniec, W. Luo and P. Sarnak (\textbf{\cite{IwLuSa}}) showed that if $\theta$ is admissible then it is possible to choose $\beta$ with:
\begin{equation}
\label{Fouriercontrol}
\forall h\in\mathbb{Z}^*, \quad \rho_j(h)\ll_{\varepsilon}\frac{(\vert h\vert a_1a_2(1+\vert r_j\vert))^{\varepsilon}}{\sqrt{a_1a_2}}\cosh{\left(\frac{\pi r_j}{2}\right)}\vert h\vert^{\theta-\frac{1}{2}}
\end{equation}
for any $j\geq 1$ and for any $\varepsilon>0$\footnote{P. Michel provided a useful averaged version over the spectrum of this upper-bound in \textbf{\cite{Mi1}}.}. P. Sarnak (\textbf{\cite{Sa}}) proved the following highly non-trivial individual estimate for the triple products $\left( u_j,\overline{V}\right)$:
\begin{equation}
\label{scalarproduct1}
\forall j\geq 1, \quad \left( u_j,\overline{V}\right)\ll_{g}\sqrt{a_1a_2}\left(1+\vert r_j\vert\right)^{k_g+1}\exp{\left(\frac{-\pi\vert r_j\vert}{2}\right)}
\end{equation}
The crucial fact is that the exponential growth in $j$ of $\rho_j(h)$ is balanced by the exponential decay in $j$ of $\left( u_j,\overline{V}\right)$. Using this, P. Sarnak proved that $D_{h}(g,a_1,a_2;s)$ admits an holomorphic continuation to $\Re{(s)}>\frac{1}{2}+\theta+\varepsilon$ under $H_2(\theta)$ for any $\varepsilon>0$ (\textbf{\cite{Sa}}). The continuous analogue being true, one obtains thanks to Weyl's law for the spectrum and an estimate for the number of cusps of the congruence subgroup $\varGamma_0(Da_1a_2)$ the following estimate for the triple products on average over the spectrum
\begin{multline}
\label{averagesarnak}
\sum_{\vert r_j\vert\leq R}\left\vert\left(u_j,\overline{V}\right)\right\vert^2\exp{\left(\pi\vert r_j\vert\right)}+\frac{1}{4\pi}\sum_{\kappa\in\text{Cusp}(\varGamma_0(Da_1a_2))}\\
\times\int_{-R}^{R}\left\vert\left(E_\kappa\left(.,\frac{1}{2}+it\right),\overline{V}\right)\right\vert^2\exp{\left(\pi\vert t\vert\right)}\mathrm{d}t\ll_{g, \varepsilon}(a_1a_2R)^\varepsilon(a_1a_2)^2R^{2k_g+\mathbf{x}}
\end{multline}
with $\mathbf{x=4}$. In fact, B. Krötz and R.J. Stanton (\textbf{\cite{KrSt}}) obtained the same estimate but with $\mathbf{x=0}$. Note that the optimality of this last estimate with respect to the parameter $R$ was already proved by A. Good (\textbf{\cite{Go}}). Moreover E. Kowalski (\textbf{\cite{Ko2}}) computed the dependency in the level of $g$. We can now state:
\begin{theorem}
\label{spectralaverageKrSt}
Let $r=q^\alpha\tilde{r}\in\mathbb{N}^*$ with $\alpha\in\mathbb{N}$ and $\tilde{r}\wedge q=1$. If $ \theta$ is admissible and $\Psi$ satisfies \eqref{Hsupport} and \eqref{Hcontrole} then
\begin{multline*}
\Sigma_r(\Psi,g;a_1,a_2)\ll_{\varepsilon ,g}q^\varepsilon\left(\frac{\sup{(Z,Y)}}{\inf{(Z,Y)}}\right)^{\frac{k_{g}-1}{2}+1+\varepsilon}\frac{(a_1a_2)^{\frac{1}{2}}}{q^{\frac{\alpha}{2}}\tilde{r}^{\frac{1}{2}+\theta}}P^{2+\varepsilon}\sup{(Z,Y)}^{1+\theta+\varepsilon} \\
\times\left(\sup{\left(1,\frac{\sup{(Z,Y)}^{\frac{1}{2}+\varepsilon}}{q^{\frac{\alpha}{2}+\varepsilon}\sqrt{a_1a_2}}\right)}+\frac{\delta_{q\mid r}}{q^{\frac{1}{2}+\theta(\alpha+1)}}\sup{\left(1,\frac{\sup{(Z,Y)}^{\frac{1}{2}+\varepsilon}}{q^{\frac{\alpha+1}{2}+\varepsilon}\sqrt{a_1a_2}}\right)}\right)
\end{multline*}
for any $\varepsilon>0$.
\end{theorem}
\begin{remark}
The shifted convolution problem is said to be balanced when $Y$ and $Z$ are of the same size and unbalanced else. In the balanced case, theorem \ref{spectralaverageKrSt} is better than theorem \ref{delta} whereas it is not the case in the unbalanced case. At least two reasons for that:
\begin{itemize}
\item
in theorem \ref{delta}, $Y$ and $Z$ are almost symmetric parameters,
\item
in theorem \ref{spectralaverageKrSt}, $\left(\frac{\sup{(Z,Y)}}{\inf{(Z,Y)}}\right)^{\frac{k_{g}-1}{2}+1}$ is large in the unbalanced case (especially when the weight of $g$ is large).
\end{itemize}
Our next applications will require the use of both theorems depending on the range of the parameters $Y$ and $Z$.
\end{remark}

\subsubsection{Flavour of the proof of theorem \ref{spectralaverageKrSt}}
When solving the shifted convolution problem on average via the spectral method on average, the main issue is to deal with smooth sums of Fourier coefficients of automorphic forms of the following shape
\begin{equation}
\label{sumh}
\sum_{h<0}\frac{\sqrt{q\vert h\vert}\,\,\overline{\rho_j(-qh)}}{(q\vert h\vert)^{s-\frac{1}{2}}}\widehat{\Psi}(qh,s)
\end{equation}
with $\Re{(s)}=\frac{1}{2}+\theta+\varepsilon$. Up to harmless factors, such a sum equals
\begin{equation*}
\frac{1}{2i\pi}\int_{(\frac{1}{2}-\theta)}L\left(\overline{\widetilde{u_j}},z+s-\frac{1}{2}\right)q^{-(s+z-1)}\widetilde{\Psi}(z,s)\mathrm{d}z
\end{equation*}
where $\widetilde{u_j}$ is the underlying primitive form of $u_j$ of level at most $Da_1a_2$ and $\widetilde{\Psi}$ is an integral transform of $\widehat{\Psi}$. Thus, bounding sums of Fourier coefficients like \eqref{sumh} turns out to bounding $L$-functions like $L\left(\overline{\widetilde{u_j}},.\right)$ on the critical line in the level aspect. Of course, the maximal saving would come from Lindelöf hypothesis but as we average over a family of Maass forms of level $Da_1a_2$ large sieve inequalities will achieve Lindelöf hypothesis on average.

\subsubsection{Proof of theorem \ref{spectralaverageKrSt}}
We set
\begin{equation*}
\Sigma_{r}(\Psi,g;a_1,a_2):=\Sigma_{r}^{\text{disc}}(\Psi,g;a_1,a_2)+\Sigma_{r}^{\text{cont}}(\Psi,g;a_1,a_2)
\end{equation*}
where
\begin{multline}
\label{discrete}
\Sigma_{r}^{\text{disc}}(\Psi,g;a_1,a_2):=\frac{1}{2i\pi\sqrt{a_1a_2}}\int_{\left(\frac{1}{2}+\theta+\varepsilon\right)}2^{s+k_g-2}\pi^{k_g-\frac{1}{2}} \\
\times\sum_{j\geq 1}\frac{\Gamma\left(\frac{\vrule width 0pt height 10pt depth 6pt s-\frac{1}{2}+ir_j}{\vrule width 0pt height 10pt depth 2pt 2}\right)\Gamma\left(\frac{\vrule width 0pt height 10pt depth 6pt s-\frac{1}{2}-ir_j}{\vrule width 0pt height 10pt depth 2pt 2}\right)}{\Gamma(s+k_g-1)}\left( u_j,\overline{V}\right)\sum_{h\neq 0}\frac{\sqrt{r\vert h\vert} \, \overline{\varepsilon_j}\,\overline{\rho_j(rh)}}{\vert rh\vert^{s-\frac{1}{2}}}\widehat{\Psi}(rh,s)\mathrm{d}s
\end{multline}
is the contribution of the discrete spectrum and
\begin{multline}
\label{continuous}
\Sigma_{r}^{\text{cont}}(\Psi,g;a_1,a_2):=\frac{1}{8i\pi^2\sqrt{a_1a_2}}\int_{\left(\frac{1}{2}+\theta+\varepsilon\right)}2^{s+k_g-2}\pi^{k_g-\frac{1}{2}} \\
\times\int_{t=-\infty}^{+\infty}\frac{\Gamma\left(\frac{\vrule width 0pt height 10pt depth 6pt s-\frac{1}{2}+it}{\vrule width 0pt height 10pt depth 2pt 2}\right)\Gamma\left(\frac{\vrule width 0pt height 10pt depth 6pt s-\frac{1}{2}-it}{\vrule width 0pt height 10pt depth 2pt 2}\right)}{\Gamma(s+k_g-1)}\sum_{\kappa\in \text{Cusp}\left(\varGamma_{0}(Da_1a_2)\right)}\left( E_{\kappa}\left(.,\frac{1}{2}+it\right),\overline{V}\right) \\
\sum_{h\neq 0}\frac{\sqrt{r\vert h\vert}\, \overline{\rho_\kappa(rh,t)}}{(r\vert h\vert)^{s-\frac{1}{2}}}\widehat{\Psi}(rh,s)\mathrm{d}t\mathrm{d}s
\end{multline}
is the contribution of the continuous spectrum. We will only give some details for the estimate of $\Sigma_{r}^{\text{disc}}(\Psi,g;a_1,a_2)$ but the same method is available for $\Sigma_{r}^{\text{cont}}(\Psi,g;a_1,a_2)$ (confer \textbf{\cite{Ri}}). As $q$ is coprime with $Da_1a_2$, \eqref{Maassmulti} implies that
\begin{equation}
\label{MaassMaass}
\forall h\in\mathbb{Z}^*, \quad\sqrt{r\vert h\vert}\rho_j(rh)=\sqrt{\tilde{r}\vert h\vert}\rho_j(\tilde{r}h)\lambda_{j}(q^\alpha)-\delta_{q\mid r}\delta_{q\mid h}\sqrt{\tilde{r}\frac{h}{q}}\rho_j\left(\tilde{r}\frac{h}{q}\right).
\end{equation}
We study only the contribution coming from the first term of \eqref{MaassMaass}. By Strirling's formula and Cauchy-Schwarz inequality, the contribution of the discrete spectrum is bounded by
\begin{multline*}
\ll\frac{1}{\sqrt{a_1a_2}}\int_{\left(\frac{1}{2}+\theta+\varepsilon\right)}\left(1+\vert\Im{(s)}\vert\right)^{-k_g} \\
\times\left(\sum_{\vert r_j\vert\leq 1+\vert\Im{(s)}\vert}\left\vert\left( u_j,\overline{V}\right)\right\vert^2\cosh{(\pi r_j)}\right)^{\frac{1}{2}} \\
\times\left(\sum_{\vert r_j\vert\leq 1+\vert\Im{(s)}\vert}\frac{1}{\cosh{(\pi r_j)}}\left\vert\sum_{\vert\tilde{r}h\vert\ll\frac{\sup{(Z,Y)}}{q^\alpha}}\sqrt{\tilde{r}\vert h\vert} \, \overline{\rho_j(\tilde{r}h)}\frac{\overline{\varepsilon_j}\,\lambda_j(q^\alpha)\widehat{H}(q^\alpha\tilde{r}h,s)}{(q^\alpha\tilde{r}\vert h\vert)^{s-\frac{1}{2}}}\right\vert^2\right)^{\frac{1}{2}}\mathrm{d}s.
\end{multline*}
According to \eqref{averagesarnak} but with the refinement of B. Krötz and R.J. Stanton ($\mathbf{x=0}$), the first square-root contributes as
\begin{equation*}
\ll_{g, \varepsilon}(a_1a_2)^{1+\varepsilon}\left(1+\vert\Im{(s)}\vert\right)^{\frac{\mathbf{x}}{2}+k_g+\varepsilon}
\end{equation*} 
for any $\varepsilon>0$. The second square-root equals
\begin{equation*}
\left(\sum_{\vert r_j\vert\leq 1+\vert\Im{(s)}\vert}\frac{1}{\cosh{(\pi r_j)}}\left\vert\sum_{\vert n\vert\ll\frac{\sup{(Z,Y)}}{q^\alpha}}a_nn^{\frac{1}{2}}\overline{\rho_j(n)}\right\vert^2\right)^{\frac{1}{2}}
\end{equation*}
where one has set
\begin{equation*}
a_n:=\begin{cases}
0 & \text{if} \quad \tilde{r}\nmid n, \\
\frac{1}{(q^\alpha\tilde{r}\vert h\vert)^{s-\frac{1}{2}}}\overline{\varepsilon_j}\,\lambda_j(q^\alpha)\widehat{H}(q^\alpha\tilde{r}h,s)& \text{if} \quad n=\tilde{r}h.
\end{cases}
\end{equation*}
The large sieve inequality for the Fourier coefficients of Maass forms of weight $0$ (confer \eqref{sievedisc}) entails that this second square-root is bounded by
\begin{equation*}
\ll_\varepsilon\left((1+\vert\Im{(s)}\vert)^2+\frac{\sup{(Z,Y)}^{1+\varepsilon}}{q^{\alpha+\varepsilon}a_1a_2}\right)^{\frac{1}{2}}\vert\vert a\vert\vert_2.
\end{equation*}
According to lemma \ref{FFhat}, this is bounded by
\begin{multline*}
\ll_{g,\varepsilon,\eta}q^\varepsilon\left(\frac{\sup{(Z,Y)}}{\inf{(Z,Y)}}\right)^{\frac{k_{g}-1}{2}+\eta -1}\frac{\sup{(Z,Y)}^{1+\theta+\varepsilon}}{q^{\frac{\alpha}{2}}\tilde{r}^{\frac{1}{2}+\theta}}\sup{\left(1,\frac{\sup{(Z,Y)}^{\frac{1}{2}+\varepsilon}}{q^{\frac{\alpha}{2}+\varepsilon}\sqrt{a_1a_2}}\right)} \\
\times P^\eta\frac{\left(1+\vert\Im{(s)}\vert\right)^{\frac{\mathbf{x}}{2}+1}}{\vert s\vert^{\eta}}
\end{multline*}
and we choose $\eta=\frac{\mathbf{x}}{2}+1+\left(1+\varepsilon\right)=2+\varepsilon$ to make convergent the $s$-integral in $\Sigma_{r}^{\text{disc}}(\Psi,g;a_1,a_2)$.
\begin{flushright}
$\blacksquare$
\end{flushright}
\begin{remark}
We want to take $\mathbf{x}$ as small as possible in \eqref{averagesarnak} because following large sieve inequality, there appears a power of $P$ and $P$ may be large in our next applications. This power is precisely the number of integration by parts we have to do and grows linearly with $\mathbf{x}$. This feature puts the stress on the fact that the spectral method is not really smooth in the length of the spectrum aspect (namely in the $R$-aspect in \eqref{averagesarnak}).
\end{remark}


\section{Proofs of Proposition D and Theorem B}

Following the results of section \ref{scp}, we prove Proposition D in subsection \ref{propDcorF} and Theorem B in subsection \ref{theB}. These proofs are based on some better bound for  $\text{Errtwist}(q,\ell;\mu)$ than the one given in \eqref{errtwist}. Remember that the bound \eqref{errtwist} was obtained in \textbf{\cite{KoMiVa}} by implementing the $\delta$-symbol method. If we use the spectral method on average described in the previous section for certain ranges which depend on the weight functions instead of the $\delta$-symbol method, we can get better bounds. Once again, a key ingredient is a uniform estimate of P. Sarnak and some technical issues involve verifying that weight functions can be handled appropriately (subsection \ref{weightfunc}).

\subsection{Description of $\text{Errtwist}(q,\ell;\mu)$}

\label{weightfunc}

In \textbf{\cite{KoMiVa}}, the authors are looking for asymptotic formulas for the harmonic twisted second moments $\mathcal{M}_g^h(\mu;\ell)$. By a standard approximate functional equation for Rankin-Selberg $L$-functions (Theorem 5.3. page 98 of \textbf{\cite{IwKo}}), they are reduced to estimate sums of the form (equation (4.16) page 138 of \textbf{\cite{KoMiVa}})
\begin{multline}
\widetilde{\mathsf{M}_{g}}((\alpha,\beta);\ell):=\sum_{m,n\geq 1}\frac{\lambda_{g}(m)\lambda_{g}(n)}{\sqrt{mn}}V_{g,\alpha}\left(\frac{m}{qD}\right)V_{g,\beta}\left(\frac{n}{qD}\right) \\
\sum_{f\in S_{k}^p(q)}^{\quad h}\psi_{f}(m)\overline{\psi_{f}(n)}\lambda_{f}(\ell).
\end{multline}
where $(\alpha,\beta)=(\pm\mu,\pm\overline{\mu})$ and
\begin{equation*}
\forall z\in\left\{\pm\mu,\pm\overline{\mu}\right\}, \forall y\in\mathbb{R}_+,\quad V_{g,z}(y):=\frac{1}{2i\pi}\int_{(3)}H_{g,z}(s)\zeta^{(qD)}(1+2s)y^{-s}\frac{\mathrm{d}s}{s-z}
\end{equation*}
satisfies
\begin{equation*}
\label{proprioV}
\forall z\in\left\{\pm\mu,\pm\overline{\mu}\right\}, \forall y\in\mathbb{R}_+, \forall A>0, \quad V_{g,z}(y)\ll_A\left(1+\vert\Im{(\mu)}\vert\right)^By^{-A}
\end{equation*}
for some $B>0$. Applying Petersson's formula (remember that there are no old forms in their case) and some dyadic partitions of unity to the $m$ and $n$ sums which appears in the non-diagonal term leads to the following term (formula (7.5) of \textbf{\cite{KoMiVa}})
\begin{equation}
\label{18}
\text{Errtwist}(q,\ell;\mu)=\frac{2\pi}{i^{k}}\sum_{M,N\geq 1}\sum_{\tilde{e}e=\ell}\frac{\varepsilon_q(\tilde{e})}{\sqrt{\tilde{e}}}\sum_{ab=\tilde{e}}\frac{\mu(a)\varepsilon_{D}(a)}{\sqrt{a}}\lambda_{g}(b)T_{M,N}
\end{equation}
with
\begin{eqnarray*}
T_{M,N} & = & \sum_{\substack{c\in\mathbb{N}^* \\
c\equiv 0\mod q}}\frac{1}{c^{2}}T_{M,N}(c), \\
T_{M,N}(c) & = & c\sum_{m, n\geq 1}\lambda_g(m)\lambda_g(n)S(m,aen;c)F_{M,N}(m,n)J_{k-1}\left(\frac{4\pi\sqrt{aemn}}{c}\right).
\end{eqnarray*}
Here, $S(m,aen;c)$ is a Kloosterman sum, $J_{k-1}$ is a Bessel function of the first kind and (page 151 of \textbf{\cite{KoMiVa}})
\begin{equation*}
F_{M,N}(x,y):=\frac{1}{\sqrt{xy}}V_{g,\alpha}\left(\frac{x}{qD}\right)V_{g,\beta}\left(\frac{y}{qD}\right)\eta_M(x)\eta_N(y)
\end{equation*}
for some smooth function $\eta_M$ compactly supported in $[M/2,2M]$ satisfying for any $i\geq 0$, $x^i\eta_M^{(i)}(x)\ll 1$ and such that
\begin{equation*}
\sum_{M\geq 1}\eta_M(x)=\begin{cases}
0 & \text{ if } x\leq\frac{1}{2}, \\
1 & \text{ if } x\geq 1
\end{cases}
\end{equation*}
with $\sum_{M\leq X}1\ll\log{X}$. Thus, $F_{M,N}$ is a compactly supported function in $\left[\frac{M}{2},2M\right]\times\left[\frac{N}{2},2N\right]$ which depends on $\mu$ and satisfies (formula (7.6) of \textbf{\cite{KoMiVa}}):
\begin{equation}
\label{FMN}
\forall(i,j)\in(\mathbb{N}^*)^2, \quad x^{i}y^{j}\frac{\partial^{i+j}F_{M,N}}{\partial x^{i}\partial y^{j}}(x,y)\ll(1+\vert t\vert)^{B}(MN)^{-\frac{1}{2}}(\log{q})^{i+j}.
\end{equation}
Truncating at an admissible cost the $M$ and $N$ sums to $M, N\ll_\varepsilon (qD)^{1+\varepsilon}$ for any $\varepsilon>0$ and applying Voronoi's formula to the $m$-sum, \eqref{18} becomes
\begin{multline}
\text{Errtwist}(q,\ell;\mu)=\frac{2\pi}{i^{k}}\sum_{M,N\ll_\varepsilon(qD)^{1+\varepsilon}}\sum_{\tilde{e}e=\ell}\frac{\varepsilon_q(\tilde{e})}{\sqrt{\tilde{e}}}\sum_{ab=\tilde{e}}\frac{\mu(a)\varepsilon_{D}(a)}{\sqrt{a}}\lambda_{g}(b)T_{M,N}^{-} \\
+O_{g, \varepsilon, A}\left((1+\vert t\vert)^B\left(\frac{1}{q^A}+q^\varepsilon\frac{\sigma_g(\ell)}{\sqrt{\ell}}\right)\right)
\end{multline}
with $\sigma_g(\ell):=\sum_{d\mid\ell}\vert\lambda_g(l)\vert$ and
\begin{equation*}
T_{M,N}^{-}=\sum_{\substack{c\in\mathbb{N}^* \\
c\equiv 0\mod q}}\frac{T_{M,N}^{-}(c)}{c^2}
\end{equation*}
with
\begin{equation*}
T_{M,N}^{-}(c)=\frac{\eta_{g}(D_{2})}{\sqrt{D_{2}}}\sum_{h\neq 0}r(-h\overline{D_{2}};c)T_{h}^{-}(c)
\end{equation*}
where
\begin{eqnarray*}
T_{h}^{-}(c) & = & \sum_{m-(aeD_{2})n=h}\lambda_{g}(m)\lambda_{g}(n)G^{-}\left(\frac{m}{D_{2}},n\right), \\
G^{-}(z,y) & = & 2\pi i^{k_{g}}\int_{0}^{+\infty}J_{k_{g}-1}\left(\frac{4\pi\sqrt{zu}}{c}\right)J_{k-1}\left(\frac{4\pi\sqrt{aeyu}}{c}\right)F_{M,N}(u,y)\mathrm{d}u.
\end{eqnarray*}
Here, $D_2:=\frac{D}{D\wedge c}$, $\overline{D_{2}}$ stands for the inverse of $D_{2}$ modulo $c$, $\eta_g(D_2)$ for an Atkin-Lehner eigenvalue (of modulus one), $r$ for the Ramanujan sum and $J_{k_g-1}$ for a Bessel function of the first kind. In fact, $T_{M,N}^{-}=T_{M,N}^-(1)+T_{M,N}^-(2)$ where
\begin{equation*}
T_{M,N}^{-}(1)=\sum_{\substack{c\in\mathbb{N}^* \\
q\mid\mid c}}c^{-2}T_{M,N}^{-}(c).
\end{equation*}
We will only deal with the first term as the same method works for the second one with better results. So, $c=qc^{\prime}$ with $c^{\prime}\wedge q=1$. Expanding the Ramanujan sum leads to
\begin{equation}
\label{Th-}
T_{M,N}^{-}(c)=\frac{\eta_g(D_2)}{\sqrt{D_2}}\sum_{\widehat{q}\in\left\{1,q\right\}}\varepsilon(\widehat{q})\widehat{q}\sum_{d\mid c^{\prime}}d\mu\left(\frac{c^{\prime}}{d}\right)\sum_{h\neq 0}T_{\widehat{q}dh}^{-}(c)
\end{equation}
with
\begin{equation*}
\varepsilon(\widehat{q})=\begin{cases}
-1 & \text{if} \quad \widehat{q}=1, \\
+1 & \text{if} \quad \widehat{q}=q.
\end{cases}
\end{equation*}
Let $F$ be the following function:
\begin{equation*}
F(z,y)=2\pi i^{k_{g}}D_{2}\int_{0}^{+\infty}J_{k_{g}-1}\left(\frac{4\pi\sqrt{zx}}{c}\right)J_{k-1}\left(\frac{4\pi\sqrt{yx}}{c}\right)F_{M,N}\left(D_{2}x,\frac{y}{a_2}\right)\mathrm{d}x
\end{equation*}
which is compactly supported with respect to $y$ in $\left[\frac{NaeD_2}{2},2NaeD_2\right]$. The shifted convolution problem on average which has to be solved is 
\begin{equation*}
\Sigma_{\widehat{q}d}(F,g;1,aeD_2).
\end{equation*}
In order to get some estimates of the function $F$, one sets $Y:=Nae$, $Z_1:=\frac{c^2}{M}$, $P:=1+\sqrt{\frac{Y}{Z_1}}$ and $\mathbf{Z}:= Z_1P^2\geq\sup{\left(Z_1,Y\right)}$. We need some results about the Bessel functions which can be found in \cite{Wa}. We know that:
\begin{equation}
\label{bessel1}
\forall j\in\mathbb{N},\quad\left(\frac{x}{1+x}\right)^{j}J_{k}^{(j)}(x)\ll_{j,k}\frac{1}{\left(1+x\right)^{\frac{1}{2}}}\left(\frac{x}{1+x}\right)^{k}.
\end{equation}
More precisely, $J_{k}(x)=\exp{(ix)}V_{k}(x)+\exp{(-ix)}\overline{V_{k}}(x)$ where
\begin{equation}
\label{bessel2}
\forall j\in\mathbb{N},\quad x^{j}V_{k}^{(j)}(x)\ll_{j,k}\frac{1}{\left(1+x\right)^{\frac{1}{2}}}\left(\frac{x}{1+x}\right)^{k}.
\end{equation}
We prove:
\begin{lemma}
\label{FF}
For any natural integers $\alpha$, $\beta$, any real numbers $A_1, A_2, A_3>0$ and any non-negative real numbers $z$ and $y$,
\begin{multline*}
z^{\alpha}y^{\beta}\frac{\partial^{\alpha+\beta}F}{\partial z^{\alpha}\partial y^{\beta}}(z,y)\ll_{k, k_{g}, \alpha, \beta, A_1, A_2, A_3}\left(1+\vert t\vert\right)^{B}(\log{q})^{\alpha+\beta+A_1+A_2+A_3} \\
\times P^{\alpha+\beta}\sqrt{\frac{M}{N}}\left(\frac{\sqrt{\frac{Y}{Z_1}}}{1+\sqrt{\frac{Y}{Z_1}}}\right)^{k-1}\frac{1}{\left(1+\sqrt{\frac{Y}{Z_1}}\right)^{\frac{1}{2}}}\left(\frac{\sqrt{\frac{z}{Z_1}}}{1+\sqrt{\frac{z}{Z_1}}}\right)^{k_g-1}\frac{1}{\left(1+\sqrt{\frac{z}{Z_1}}\right)^{\frac{1}{2}}} \\
\times\frac{1}{\left(1+\frac{z}{\mathbf{Z}}\right)^{A_1}\left(1+\frac{y}{Y}\right)^{A_2}\left(1+\frac{Y}{\left(\sqrt{Z_1}+\sqrt{z}\right)^2}\right)^{A_3}}.
\end{multline*}
\end{lemma}
\noindent{\textbf{Proof of lemma \ref{FF}.}}
We give only the proof for the case $\alpha=\beta=0$. If $z<\mathbf{Z}$ then we trivially have
\begin{multline}
\label{Ftrivial}
F(z,y)\ll_{k, k_{g}}\left(1+\vert t\vert\right)^{B}\sqrt{\frac{M}{N}}\left(\frac{\sqrt{\frac{Y}{Z_1}}}{1+\sqrt{\frac{Y}{Z_1}}}\right)^{k-1}\frac{1}{\left(1+\sqrt{\frac{Y}{Z_1}}\right)^{\frac{1}{2}}} \\
\times\left(\frac{\sqrt{\frac{z}{Z_1}}}{1+\sqrt{\frac{z}{Z_1}}}\right)^{k_g-1}\frac{1}{\left(1+\sqrt{\frac{z}{Z_1}}\right)^{\frac{1}{2}}}.
\end{multline}
If $z\geq\mathbf{Z}$ then $l\geq 1$ integrations by parts lead to
\begin{equation*}
F(z,y)=4\pi i^{k_{g}}\int_{\sqrt{\frac{M}{2D_{2}}}}^{\sqrt{\frac{2M}{D_{2}}}}\left(\frac{\exp{\left(i\frac{4\pi}{c}\sqrt{z}x\right)}}{\left(i\frac{4\pi}{c}\sqrt{z}\right)^l}f^{(l)}(x)+\frac{\exp{\left(-i\frac{4\pi}{c}\sqrt{z}x\right)}}{\left(-i\frac{4\pi}{c}\sqrt{z}\right)^l}\overline{f}^{(l)}(x)\right)\mathrm{d}x
\end{equation*}
where $f(x)=xV_{k_{g}-1}\left(\frac{4\pi}{c}\sqrt{z}x\right)J_{k-1}\left(\frac{4\pi}{c}\sqrt{y}x\right)F_{M,N}\left(D_{2}x^2,\frac{y}{a_2}\right)$ satisfies
\begin{multline*}
f^{(l)}(x)\ll_{k,k_{g},l}\left(1+\vert t\vert\right)^{B}(\log{q})^{l}\frac{P^l}{\sqrt{MN}}\left(\frac{\sqrt{\frac{Y}{Z_1}}}{1+\sqrt{\frac{Y}{Z_1}}}\right)^{k-1}\frac{1}{\left(1+\sqrt{\frac{Y}{Z_1}}\right)^{\frac{1}{2}}} \\
\times\left(\frac{\sqrt{\frac{z}{Z_1}}}{1+\sqrt{\frac{z}{Z_1}}}\right)^{k_g-1}\frac{1}{\left(1+\sqrt{\frac{z}{Z_1}}\right)^{\frac{1}{2}}}\frac{1}{x^{l-1}}.
\end{multline*}
As a consequence,
\begin{multline}
\label{Fnontrivial}
F(z,y)\ll_{k, k_{g}}\left(1+\vert t\vert\right)^{B}\sqrt{\frac{M}{N}}\left(\frac{\sqrt{\frac{Y}{Z_1}}}{1+\sqrt{\frac{Y}{Z_1}}}\right)^{k-1}\frac{1}{\left(1+\sqrt{\frac{Y}{Z_1}}\right)^{\frac{1}{2}}} \\
\times\left(\frac{\sqrt{\frac{z}{Z_1}}}{1+\sqrt{\frac{z}{Z_1}}}\right)^{k_g-1}\frac{1}{\left(1+\sqrt{\frac{z}{Z_1}}\right)^{\frac{1}{2}}}P^l\left(\frac{z}{\mathbf{Z}}\right)^{-2l}.
\end{multline}
We conclude by collecting \eqref{Ftrivial} and \eqref{Fnontrivial} and by remarking that one can repeat the same procedure with $J_{k-1}$ instead of $J_{k_{g}-1}$ if $Y$ is large.
\begin{flushright}
$\blacksquare$
\end{flushright}
Let $\rho:\mathbb{R}\rightarrow\mathbb{R}$ compactly supported in $[1,2]$ satisfying  $\sum_{a\in\mathbb{N}}\rho\left(2^{-a}x\right)=1$. We set
\begin{equation*}
F_{Z}(z,y):=\rho\left(2^{-a}z\right)F(z,y)
\end{equation*}
where $Z:=2^a$. $F_{Z}$ is compactly supported in $\left[Z,2Z\right]\times \left[\frac{Y}{2},2Y\right]$ and we remark that $F(z,y)=\sum_{Z=2^a}F_{Z}(z,y)$. Lemma \ref{FF} gives
\begin{multline}
\label{Fa}
z^\alpha y^\beta\frac{\partial^{\alpha+\beta}F_{Z}}{\partial z^{\alpha}\partial y^{\beta}}(z,y)\ll_{k, k_{g}, \alpha, \beta, A_1, A_2, A_3}\left(1+\vert t\vert\right)^{B}(\log{q})^{\alpha+\beta+A_1+A_2+A_3} \\
\times P^{\alpha+\beta}\sqrt{\frac{M}{N}}\left(\frac{\sqrt{\frac{Y}{Z_1}}}{1+\sqrt{\frac{Y}{Z_1}}}\right)^{k-1}\frac{1}{\left(1+\sqrt{\frac{Y}{Z_1}}\right)^{\frac{1}{2}}}\left(\frac{\sqrt{\frac{Z}{Z_1}}}{1+\sqrt{\frac{Z}{Z_1}}}\right)^{k_g-1}\frac{1}{\left(1+\sqrt{\frac{Z}{Z_1}}\right)^{\frac{1}{2}}} \\
\times\frac{1}{\left(1+\frac{Z}{\mathbf{Z}}\right)^{A_1}\left(1+\frac{y}{Y}\right)^{A_2}\left(1+\frac{Y}{\left(\sqrt{Z_1}+\sqrt{Z}\right)^2}\right)^{A_3}}
\end{multline}
for any natural integers $\alpha$, $\beta$ and for any real numbers $A_1, A_2, A_3>0$.

\subsection{Improvement of the bound of $\text{Errtwist}(q,\ell;\mu)$ given in \eqref{errtwist}}

\label{propDcorF}

Let us recall that we want to estimate
\begin{equation*}
\text{Errtwist}^\prime(q,\ell;\mu):=\frac{2\pi}{i^{k}}\sum_{M,N\ll_\varepsilon(qD)^{1+\varepsilon}}\sum_{\tilde{e}e=\ell}\frac{\varepsilon_q(\tilde{e})}{\sqrt{\tilde{e}}}\sum_{ab=\tilde{e}}\frac{\mu(a)\varepsilon_{D}(a)}{\sqrt{a}}\lambda_{g}(b)T_{M,N}^-(1)
\end{equation*}
where
\begin{equation}
\label{estimfin}
T_{M,N}^-(1)=\frac{\eta_g(D_2)}{\sqrt{D_2}}\sum_{\widehat{q}\in\left\{1,q\right\}}\sum_{\substack{(c,d)\in\mathbb{N}^{*2} \\
q\mid\mid c \\
d\mid\frac{c}{q}}}\frac{d\widehat{q}}{c^2}\mu\left(\frac{c}{dq}\right)\sum_{\substack{Z\geq 1 \\
Z=2^a \\
a\in\mathbb{N^*}}}\Sigma_{\widehat{q}d}\left(F_Z,g;1,aeD_2\right)
\end{equation}
and that Theorem \ref{spectralaverageKrSt} implies:
\begin{theorem}
\label{spectraluni}
Let $c\geq 1$ and $d\geq 1$ some natural integers satisfying $q\mid\mid c$, $d\mid\frac{c}{q}$, $\widehat{q}\in\left\{1,q\right\}$ and $Z\geq 1$. If $\,\theta$ is admissible then:
\begin{multline*}
\Sigma_{\widehat{q}d}(F_Z,g;1,aeD_2)\ll_{\varepsilon ,g, A_1, A_3, \eta}q^\varepsilon\left(1+\vert t\vert\right)^{B}\sqrt{\frac{M}{N}}(ae)^{\frac{1}{2}}\left(\frac{\sqrt{\frac{Y}{Z_1}}}{1+\sqrt{\frac{Y}{Z_1}}}\right)^{k-1}\frac{1}{\left(1+\sqrt{\frac{Y}{Z_1}}\right)^{\frac{1}{2}}} \\
\left(\frac{\sqrt{\frac{Z}{Z_1}}}{1+\sqrt{\frac{Z}{Z_1}}}\right)^{k_g-1}\frac{1}{\left(1+\sqrt{\frac{Z}{Z_1}}\right)^{\frac{1}{2}}}\frac{1}{\left(1+\frac{Z}{\mathbf{Z}}\right)^{A_1}\left(1+\frac{Y}{\left(\sqrt{Z_1}+\sqrt{Z}\right)^2}\right)^{A_3}} \\
\left(\frac{\sup{(Z,Y)}}{\inf{(Z,Y)}}\right)^{\frac{k_{g}-1}{2}+1+\varepsilon}P^{2+\varepsilon}\frac{\sup{(Z,Y)}^{1+\theta+\varepsilon}}{\widehat{q}^{\frac{1}{2}}d^{\frac{1}{2}+\theta}}\sup{\left(1,\frac{\sup{(Z,Y)}^{\frac{1}{2}+\varepsilon}}{\widehat{q}^{\frac{1}{2}+\varepsilon}\sqrt{ae}}\right)}
\end{multline*}
for any real numbers $\eta\geq 0$ and $A_1$, $A_3$, $\varepsilon>0$.
\end{theorem}
Now, we finish the proof of Proposition D.
\noindent{\textbf{Proof of Proposition D.}} Let $0<\alpha<1$ be some real number. Setting $\mathcal{C}_{\widehat{q}d}(x,y):=\sum_{x \leq Z\leq y}\Sigma_{\widehat{q}d}\left(F_Z,g;1,aeD_2\right)$, we split the $Z$-sum occurring in \eqref{estimfin} as follows
\begin{equation}
\label{splitfin}
\sum_{\substack{Z\geq 1 \\
Z=2^a \\
a\in\mathbb{N^*}}}\Sigma_{\widehat{q}d}\left(F_Z,g;1,aeD_2\right)=\mathcal{C}_{\widehat{q}d}(1,Z_1^\alpha)+\mathcal{C}_{\widehat{q}d}(Z_1^\alpha,\textbf{Z})+\mathcal{C}_{\widehat{q}d}(\textbf{Z},+\infty)
\end{equation}
and we refer to the first (respectively second, third) term in the right-hand side of \eqref{splitfin} as the short (respectively median, long) range terms. The first point is that $F_Z$ is small when $k_g$ is large for $1\leq Z\leq Z_1^\alpha$ because proposition \ref{Fhat} implies that
\begin{multline}
\label{FZpetitstermes}
F_{Z}(z,y)\ll_{k, k_{g}}\left(1+\vert t\vert\right)^{B}q^{\varepsilon} \\
\times\sqrt{\frac{M}{N}}\left(\frac{\sqrt{\frac{Y}{Z_1}}}{1+\sqrt{\frac{Y}{Z_1}}}\right)^{k-1}\frac{1}{\left(1+\sqrt{\frac{Y}{Z_1}}\right)^{\frac{1}{2}}}\left(\frac{Z}{Z_1}\right)^{\frac{k_g-1}{2}}
\end{multline}
for $1\leq Z\leq Z_1^\alpha$ and for any $\varepsilon>0$. As a consequence, the short range terms do not restrict the length of the mollifier at least when $k_g$ is large. More precisely, if $k_g>1+\frac{5}{2(1-\alpha)}$ then one gets thanks to theorem \ref{delta} (that is to say the $\delta$-method symbol):
\begin{equation}
\label{petitstermes}
\sum_{\widehat{q}\in\left\{1,q\right\}}\sum_{\substack{(c,d)\in\mathbb{N}^{*2} \\
q\mid\mid c \\
d\mid\frac{c}{q}}}\frac{d\widehat{q}}{c^2}\mu\left(\frac{c}{dq}\right)\mathcal{C}_{\widehat{q}d}(1,Z_1^\alpha)\ll_{k, g}(1+\vert t\vert)^B\frac{1}{q^\delta}
\end{equation}
for some $\delta>0$. The long range terms "weakly" restrict the length of the mollifier. This is mainly caused by the factor $\left(\frac{\mathbf{Z}}{Z}\right)^{A_1}$ in theorem \ref{spectraluni} with $A_1$ as large as needed. Applying this theorem (that is to say the spectral method on average), one gets for $\theta$ admissible and $k>k_g+\frac{21}{4}+\frac{\theta}{2}$:
\begin{equation}
\label{grandstermes}
\sum_{\widehat{q}\in\left\{1,q\right\}}\sum_{\substack{(c,d)\in\mathbb{N}^{*2} \\
q\mid\mid c \\
d\mid\frac{c}{q}}}\frac{d\widehat{q}}{c^2}\mu\left(\frac{c}{dq}\right)\mathcal{C}_{\widehat{q}d}(\mathbf{Z},+\infty)\ll_{k, g, \varepsilon}(1+\vert t\vert)^Bq^\varepsilon\left(\frac{\ell^{\frac{5}{4}+\frac{\theta}{2}}}{q^{\frac{1}{2}-\theta}}+\frac{\ell^{2+\theta}}{q^{\frac{1}{2}-\theta}}\right)
\end{equation}
for any $\varepsilon>0$. The main restriction comes from the median range terms. Applying theorem \ref{spectraluni} (that is to say the spectral method on average), one gets for $\theta$ admissible and $k>k_g+\frac{21}{4}+\frac{\theta}{2}$
\begin{multline}
\label{termesmedians}
\sum_{\widehat{q}\in\left\{1,q\right\}}\sum_{\substack{(c,d)\in\mathbb{N}^{*2} \\
q\mid\mid c \\
d\mid\frac{c}{q}}}\frac{d\widehat{q}}{c^2}\mu\left(\frac{c}{dq}\right)\mathcal{C}_{\widehat{q}d}(Z_1^\alpha,\mathbf{Z})\ll_{k, g, \varepsilon}(1+\vert t\vert)^Bq^\varepsilon \\
\times\left(\frac{\ell^{\frac{5}{4}+\frac{\theta}{2}}}{q^{\frac{1}{2}-\theta}}+\frac{\ell^{2+\theta}}{q^{\frac{1}{2}-\theta}}+\frac{\ell^{\frac{9}{4}+\frac{\theta}{2}-\alpha}}{q^{\alpha-\frac{1}{2}-\theta}}\right)
\end{multline}
for any $\varepsilon>0$. Collecting these estimates, one gets for $\theta$ admissible, $k_g>1+\frac{5}{2(1-\alpha)}$ and $k>k_g+\frac{21}{4}+\frac{\theta}{2}$
\begin{equation}
\text{Errtwist}^\prime(q,\ell;\mu)=O_{\varepsilon,k,g}\left((q\ell)^{\varepsilon}(1+\vert\Im{(\mu)}\vert)^{B}\left(\frac{\ell^{\frac{5}{4}+\frac{\theta}{2}}}{q^{\frac{1}{2}-\theta}}+\frac{\ell^{2+\theta}}{q^{\frac{1}{2}-\theta}}+\frac{\ell^{\frac{9}{4}+\frac{\theta}{2}-\alpha}}{q^{\alpha-\frac{1}{2}-\theta}}\right)\right)
\end{equation}
and so
\begin{equation}
\text{Errsec}(q,L;\mu)=O_{\varepsilon,k,g}\left((qL)^{\varepsilon}(1+\vert\Im{(\mu)}\vert)^{B}\left(\frac{L^{5+2\theta}}{q^{\frac{1}{2}-\theta}}+\frac{L^{8+2\theta}}{q^{\frac{1}{2}-\theta}}+\frac{L^{\frac{11}{2}+\theta-2\alpha}}{q^{\alpha-\frac{1}{2}-\theta}}\right)\right).
\end{equation}
Thus, every $\Delta$ strictly less than
\begin{equation*}
\inf{\left(\frac{1-2\theta}{4(5+2\theta)},\frac{2\alpha-1-2\theta}{2(11+2\theta-4\alpha)}\right)}
\end{equation*}
is effective provided $k$ and $k_g$ are large enough. We choose $\alpha:=\frac{7}{8}+\frac{\theta}{6}+\frac{\theta^2}{6}$ to maximize the last quantity.
\begin{flushright}
$\blacksquare$
\end{flushright}


\subsection{A new subconvexity bound}

\label{theB}
As a consequence of the improvement of the bound of $\text{Errtwist}(q,\ell;\mu)$, we prove the new subconvexity bound of Rankin-Selberg $L$-functions given in Theorem B by applying the amplification method. Setting (as in \cite{KoMiVa}) for $L\geq 1$ an integer and $\overrightarrow{x}=\left(x_\ell\right)_{1\leq\ell\leq L}$ a sequence of complex numbers satisfying $x_\ell=0$ if $q\mid\ell$
\begin{equation*}
L_g(\mu,\overline{\mu};L;\overrightarrow{x}):=\sum_{1\leq\ell\leq L}x_\ell\widetilde{\mathsf{M}_g}((\mu,\overline{\mu});\ell),
\end{equation*}
one has according to \cite{KoMiVa} (page 151)
\begin{equation*}
L_g(\mu,\overline{\mu};L;\overrightarrow{x})=\sum_{1\leq\ell\leq L}x_\ell\text{Errtwist}(q,\ell;\mu)+O_{\varepsilon,g}\left((1+\vert t\vert)^Bq^\varepsilon\sum_{1\leq\ell\leq L}\vert x_\ell\vert\frac{\sigma_{g}(\ell)}{\sqrt{\ell}}\right)
\end{equation*}
which leads to:
\begin{proposition}
\label{averagepropoD}
Let $\alpha\in\left]0,1\right[$. Let $g$ be a primitive cusp form of square-free level $D$, weight $k_g>1+\frac{5}{2(1-\alpha)}$ and trivial nebentypus and $\mu\in\mathbb{C}$. Assume that $q$ is a prime coprime with $D$ and that $k\geq k_g+6$. If $\theta$ is admissible and $\left\vert\Re{(\mu)}\right\vert\ll\frac{1}{\log{q}}$ then for any $1\leq L<q$,
\begin{multline*}
L_g(\mu,\overline{\mu};L;\overrightarrow{\vrule width 0pt height 6pt depth 7pt x})\ll_{\varepsilon, k, g}(qL)^{\varepsilon}\left(1+\vert t\vert\right)^B\Bigg(\sum_{1\leq\ell\leq L}\vert x_\ell\vert\frac{\sigma_{g}(\ell)}{\sqrt{\ell}} \\
+\left(\frac{L^{2+\theta}}{q^{\frac{1}{2}-\theta}}+\frac{L^{\frac{9}{4}+\frac{\theta}{2}-\alpha}}{q^{\alpha-\frac{1}{2}-\theta}}\right)\vert\vert\overrightarrow{\vrule width 0pt height 6pt depth 7pt x}\vert\vert_1\Bigg)
\end{multline*}
for any $\varepsilon>0$.
\end{proposition}
\noindent{\textbf{Proof of Theorem B.}}
As in \cite{KoMiVa}, let $\mathcal{Q}(.)$ be the following quadratic form:
\begin{equation*}
\mathcal{Q}(\overrightarrow{x}):=\sum_{f\in S_{k}^p(q)}^{\quad h}\left\vert\sum_{1\leq\ell\leq L}x_\ell\lambda_{f}(\ell)\right\vert^2\left\vert L\left(f\times g,\frac{1}{2}+\mu\right)\right\vert^2
\end{equation*}
for $L<\sqrt{q}$. We define $\overrightarrow{X}:=\left(X_\ell\right)_{1\leq\ell\leq L^2}$ with:
\begin{equation*}
X_\ell:=\sum_{d\geq 1}\sum_{\substack{\ell_1\ell_2=\ell \\ 1\leq \ell_1,\ell_2\leq \frac{L}{d}}}\overline{x_{d\ell_1}}x_{d\ell_2}.
\end{equation*}
It is proved in \textbf{\cite{KoMiVa}} that $\mathcal{Q}(\overrightarrow{x})\ll_{g}L_{g}(\mu,\overline{\mu};L^2;\overrightarrow{X})$. This leads to:
\begin{multline*}
\left\vert\sum_{1\leq \ell\leq L}x_\ell\lambda_{f}(\ell)\right\vert^2\left\vert L\left(f\times g,\frac{1}{2}+\mu\right)\right\vert^2\ll_{\varepsilon,k,g}(qL)^{\varepsilon}(1+\vert t\vert)^Bq\Bigg(\vert\vert \overrightarrow{\vrule width 0pt height 6pt depth 7pt x}\vert\vert_2^2 \\
+\left(\frac{L^{4+2\theta}}{q^{\frac{1}{2}-\theta}}+\frac{L^{\frac{9}{2}+\theta-2\alpha}}{q^{\alpha-\frac{1}{2}-\theta}}\right)\vert\vert\overrightarrow{\vrule width 0pt height 6pt depth 7pt x}\vert\vert_1\Bigg)
\end{multline*}
for any $\varepsilon>0$. We choose the following classical lacunary $GL(2)$-amplifier:
\begin{equation*}
x_\ell:=\begin{cases}
-1 & \text{if} \quad \ell=p^2 \quad \text{with} \quad p\in\mathcal{P},  \quad p\leq\sqrt{L}, \\
\lambda_{f}(p) & \text{if} \quad \ell=p \quad \text{with} \quad p\in\mathcal{P},  \quad p\leq\sqrt{L}, \\
0 & \text{else.}
\end{cases}
\end{equation*}
With such a choice,
\begin{equation*}
\left\vert L\left(f\times g,\frac{1}{2}+\mu\right)\right\vert^2\ll_{\varepsilon,k,g}(qL)^{\varepsilon}(1+\vert t\vert)^B\Bigg(\frac{q}{\sqrt{L}}+q^{\frac{1}{2}+\theta}L^{4+2\theta}+q^{\frac{3}{2}+\theta-\alpha}L^{\frac{9}{2}+\theta-2\alpha}\Bigg).
\end{equation*}
Setting $L=q^{2x}$ with $0<x<\frac{1}{4}$, we have:
\begin{multline*}
\left\vert L\left(f\times g,\frac{1}{2}+\mu\right)\right\vert^2\ll_{\varepsilon,k,g}q^{\varepsilon}(1+\vert t\vert)^B\inf_{0<x<\frac{1}{4}}\Bigg(q^{1-x}+q^{\frac{1}{2}+\theta+(8+4\theta)x} \\
+q^{\frac{3}{2}+\theta-\alpha+(9+2\theta-4\alpha)x}\Bigg).
\end{multline*}
Finally, we choose $x=\frac{1-2\theta}{2(9+4\theta)}$ and
$\alpha:=\frac{19}{22}+\frac{2}{11}\theta+\frac{2}{11}\theta^2$ to minimize the right-hand side which achieves the proof of Theorem B for $j=0$ in a neighbourhood of the critical line. The other cases ($j\neq 0$) follow from Cauchy's inequalities.
\begin{flushright}
$\blacksquare$
\end{flushright}


\appendix

\section{The harmonic mollified second moment away from the critical point}

\label{away}

The aim of this part is to prove the following bound of $\mathcal{W}^{h}(g;\mu)$ when $\mu$ is on the right of the origin:
\begin{theorem}
\label{propoE}
Let $g$ be a primitive cusp form of square-free level $D$ and trivial nebentypus and $f$ be a non-negative function satisfying:
\begin{eqnarray*}
\lim_{q\rightarrow +\infty}f(q) & = & +\infty, \\
f(q) & = & o(\log{q}).
\end{eqnarray*}
Assume that $q$ is prime, coprime with $D$. If $\Re{(\mu)}\geq\frac{f(q)}{\log{q}}$ and $\Delta$ is effective then for any $0<a<4\Delta(1-\Upsilon)$, we have
\begin{equation}
\mathcal{W}^{h}(g;\mu)=A_q^h[1]+O_{k,g}\left((1+\vert\Im{(\mu)}\vert)^{B}q^{-a\Re{(\mu)}}\right)
\end{equation}
for some absolute constant $B>0$.
\end{theorem}
We only give a sketch of the proof of this theorem based on two lemmas and a classical convexity argument. As usual, $\mu$ is a complex number and $\tau:=\Re{(\mu)}$ and $t:=\Im{(\mu)}$. On one hand, just on the right of the critical point, we have:
\begin{lemma}
\label{A}
If $\tau=\frac{f(q)}{\log{q}}>0$ where $f$ is a non-negative function satisfying:
\begin{eqnarray*}
\lim_{q\rightarrow +\infty}f(q) & = & +\infty, \\
f(q) & = & o(\log{q})
\end{eqnarray*}
and $\Delta$ is effective then
\begin{equation}
\mathcal{W}^{h}(g;\mu)\ll_g\left(1+\vert t\vert\right)^{B}.
\end{equation}
\end{lemma}
\noindent{\textbf{Proof of lemma \ref{A}.}}
According to remark \ref{extension}, $\mathcal{W}_g^h(\mu)\ll_{k,g} 1$ if $\vert t\vert\ll 1$. So, we may assume that $\vert t\vert\gg 1$. According to proposition \ref{secondstep} and its proof, we know that up to an admissible error term
\begin{multline*}
\!\!\!\sum_{(\alpha,\beta)=(\pm\mu,\pm\overline{\mu})}\!\!\!\varepsilon_{f\times g}(\alpha,\beta)\mathsf{W}_{g}(\alpha,\beta)\!\!=\!\!\!\sum_{(\alpha,\beta)=(\pm\mu,\pm\overline{\mu})}\Psi(\alpha,\beta)\frac{1}{(2i\pi)^{2}}\int_{(3)}\int_{(3)}h_g(\alpha,\beta,s_1,s_2)\\
\times n_g(s_1,\mu_1,\alpha,\beta)n_g(s_2,\mu_2,\alpha,\beta)L(g\times g,1+s_1+s_2+2\tau)\frac{\mathrm{d}s_{1}}{s_{1}^{2}}\frac{\mathrm{d}s_{2}}{s_{2}^{2}}
\end{multline*}
with for $z\in\left\{\mu_1,\mu_2\right\}$,
\begin{multline*}
n_g(s,z,\alpha,\beta):=\frac{1}{\log{L}}L^{(1-\Upsilon)s}\left(\widehat{P^{\prime}_{L}}(s)L^{\Upsilon s}-\frac{1}{1-\Upsilon}\widehat{R^{\prime}_{L^{1-\Upsilon}}}(s)\right) \\
\frac{L(g\times g,1+s+2z)}{\zeta^{(D)}(1+s+2z)L(g\times g,1+\alpha+s+z)L(g\times g,1+\beta+s+z)}.
\end{multline*}
We are going to evaluate each term occurring in the previous equality separately. One should remark that
\begin{equation*}
\Psi(\alpha,\beta)\ll(1+\vert t\vert)^B\log^A{(q)}(qD)^{-2\tau+\alpha+\beta}
\end{equation*}
for some absolute constants $A$ and $B$ and also that:
\begin{equation*}
(qD)^{-2\tau+\alpha+\beta}\ll\begin{cases}
1 & \text{if } (\alpha,\beta)=(\mu,\overline{\mu}), \\
\exp{(-2f(q))} & \text{else.}
\end{cases}
\end{equation*}
So, we are going to give details only for the worst case which is $(\alpha,\beta)=(\mu,\overline{\mu})$. We shift the $s_1$-contour and the $s_2$-contour to $\left(+\frac{c_1}{\log{q}}\right)$ without crossing any poles ($c_1>0$). Then, we shift the $s_1$ contour to $\left(-\frac{c_2}{\log{q}}\right)$ with $c_1<(1-\Upsilon)c_2<2(1-\Upsilon)$ hitting some poles at $s_1=0$ and $s_1=s_2-2\tau$. The residual integral is bounded by $\log^A{(q)}\exp{(-4\Delta((1-\Upsilon)c_2-c_1) f(q))}$ for some $A$ which is admissible. The contribution of the pole at $s_1=s_2-2\tau$ is bounded by $\exp{(-2\Delta(2(1-\Upsilon)-c_1)f(q))}$ $\times\log^A{(q)}$ for some $A$ which is admissible. The contribution of the pole at $s_1=0$ is given by
\begin{multline*}
\frac{\Psi(\mu,\overline{\mu})}{\zeta^{(D)}(1+2\mu)L(g\times g,1+2\tau)}\frac{1}{2i\pi\log{L}}\int_{\left(+c_1\frac{f(q)}{\log{q}}\right)}h_g(\mu,\overline{\mu},0,s_2) \\
\times\frac{L(g\times g,1+s_2+2\overline{\mu})}{\zeta^{(D)}(1+s_2+2\overline{\mu})L(g\times g,1+2\overline{\mu}+s_2)} \\
\times L^{(1-\Upsilon)s_2}\left(\widehat{P^{\prime}_{L}}(s_2)L^{\Upsilon s_2}-\frac{1}{1-\Upsilon}\widehat{R^{\prime}_{L^{1-\Upsilon}}}(s_2)\right)\frac{\mathrm{d}s_2}{s_2^2}.
\end{multline*}
We shift the $s_2$-contour to $\left(-\frac{f(q)}{\log{q}}\right)$ hitting only a pole at $s_2=0$. The residual integral is bounded by $\log^A{(q)}\exp{(-2\Delta f(q))}$ for some $A$ which is still admissible and the contribution of the pole is given by:
\begin{equation*}
\frac{\Psi(\mu,\overline{\mu})}{\zeta^{(D)}(1+2\mu)\zeta^{(D)}(1+2\overline{\mu})L(g\times g,1+2\tau)}=\frac{\varphi(q)}{q}\frac{\zeta_{(q)}(1+2\overline{\mu})}{\zeta^{(D)}(2(1+2\tau))}h_g(\mu,\overline{\mu},0,0)
\end{equation*}
which is bounded.
\begin{flushright}
$\blacksquare$
\end{flushright}
On the other hand, very far away $\frac{1}{2}$ in the domain of absolute convergence, we have:
\begin{lemma}
\label{B}
If $\tau>\frac{1}{2}+\varepsilon$ then
\begin{equation}
A_q^{h}\left[\left\vert\mathcal{L}\left(.\times g,\frac{1}{2}+\mu\right)-1\right\vert^{2}\right]\ll_{\varepsilon} q^{-4\Delta(1-\Upsilon)\left(\tau-\left(\frac{1}{2}+\varepsilon\right)\right)}
\end{equation}
for any $\varepsilon>0$.
\end{lemma}
This lemma is an easy consequence of \ref{loin} as we are in the domain of absolute convergence of Rankin-Selberg $L$-functions.\newline
\noindent{\textbf{Proof of theorem \ref{propoE}.}}
Lemma \ref{A} and \ref{B} together with a Phragmen-Lindelöf type principle for subharmonic functions which can be found in \cite{Ko} give
\begin{equation*}
A_q^{h}\left[\left\vert\mathcal{L}\left(.\times g,\frac{1}{2}+\mu\right)-1\right\vert^{2}\right]\ll_{\varepsilon} (1+\vert t\vert)^Bq^{\alpha(\tau)}
\end{equation*}
where $\alpha$ is the affine function satisfying $\alpha\left(\tau_0\right)=-4\Delta(1-\Upsilon)\left(\tau_0-\left(\frac{1}{2}+\varepsilon\right)\right)$ and $\alpha\left(\frac{f(q)}{\log{q}}\right)=0$. This leads to
\begin{equation*}
A_q^{h}\left[\left\vert\mathcal{L}\left(.\times g,\frac{1}{2}+\mu\right)-1\right\vert^{2}\right]\ll_{\varepsilon} (1+\vert t\vert)^Bq^{-\frac{4\Delta(1-\Upsilon)\left(\tau_0-\left(\frac{1}{2}+\varepsilon\right)\right)}{\tau_0}\tau}
\end{equation*}
which concludes the proof by choosing $\varepsilon$ small enough and $\tau_{0}$ large enough.
\begin{flushright}
$\blacksquare$
\end{flushright}

\section{Bounding the contribution of old forms}
\label{kgrd}

The main purpose of this appendix is to deal with the eventual existence of old forms in $S_k(q)$ (when $k\geq 12$). In other words, we prove that \eqref{twist} still holds even if there are some old forms in $S_{k}(q)$. Let $N\geq 1$. We define for every integers $m, n\geq 1$ the operator $\Delta_N$ by:
\begin{equation}
\label{DeltaN}
\Delta_{N}(m,n):=\delta_{m=n}+\frac{2\pi}{i^{k}}\sum_{\substack{c\in\mathbb{N}^* \\ c\equiv 0 \mod N}}\frac{S(m,n;c)}{c}J_{k-1}\left(\frac{4\pi\sqrt{mn}}{c}\right).
\end{equation}
where $S(m,n;c)$ is the Kloosterman sum for which we recall Weil's bound (confer \cite{We}):
\begin{equation}
\label{weil}
\vert S(m,n;c)\vert\leq\tau(c) (m,n,c)^{\frac{1}{2}}\sqrt{c}.
\end{equation}
Then, Petersson trace formula expresses this operator as an average over an orthogonal basis $\mathcal{B}_{k}(N)$ of $S_{k}(N)$:
\begin{equation}
\label{petersson}
\Delta_{N}(m,n)=\sum_{h\in\mathcal{B}_{k}(N)}\omega_N(h)\psi_{h}(m)\overline{\psi_{h}(n)}
\end{equation}
where $\omega_h(N)\ll_k\frac{\log{N}}{N}$ uniformly with respect to $h$ according to \cite{GoHoLi}. H. Iwaniec, W. Luo and P. Sarnak have restricted themselves in \cite{IwLuSa} to average over primitive forms:
\begin{theorem}[H. Iwaniec-W. Luo-P. Sarnak (2001)]
Let $N\geq 1$  be a square-free number. 
\begin{equation}
\label{peterssonnew}
\sum_{h\in S_{k}^p(N)}\omega_{N}(h)\lambda_{h}(n)\lambda_h(n)=\frac{1}{N}\sum_{LM=N}\frac{\mu(L)M}{\nu(n\wedge L)}\sum_{l\mid L^{\infty}}\frac{1}{l}\Delta_{M}(ml^2,n).
\end{equation}
\end{theorem}
The authors showed in \textbf{\cite{KoMiVa}} using Petersson trace formula (confer section \ref{weightfunc} or page 138 of \textbf{\cite{KoMiVa}}) that if there are no old forms in $S_{k}(q)$ and if $1\leq\ell<q$ then
\begin{equation*}
(qD)^{2\Re{(\mu)}}\mathcal{M}_{g}^{h}(\mu;\ell)=\sum_{(\alpha,\beta)=(\pm\mu,\pm\overline{\mu})}\varepsilon_{f\times g}(\alpha,\beta)\widetilde{\mathsf{M}_{g}}((\alpha,\beta);\ell)
\end{equation*}
and that for any $(\alpha,\beta)=(\pm\mu,\pm\overline{\mu})$:
\begin{multline}
\label{start}
\widetilde{\mathsf{M}_{g}}((\alpha,\beta);\ell)=\sum_{\tilde{e}e=\ell}\frac{\varepsilon_q(\tilde{e})}{\sqrt{\tilde{e}}}\sum_{ab=\tilde{e}}\frac{\mu(a)\varepsilon_D(a)}{\sqrt{a}}\lambda_g(b) \\
\sum_{m,n\geq 1}\frac{\lambda_{g}(m)\lambda_{g}(n)}{\sqrt{mn}}V_{g,\alpha}\left(\frac{m}{qD}\right)V_{g,\beta}\left(\frac{a\tilde{e}n}{qD}\right)\Delta_q(m,aen).
\end{multline}
In our case, there are some old forms as the weight $k$ may be large but their contribution is small.
\begin{proposition}
\label{newresult}
Let $g$ be a primitive cusp form of square-free level $D$ and trivial nebentypus. Assume that $q$ is prime, coprime with $D$. If $\mu\in\mathbb{C}$ and $1\leq\ell<q$ then
\begin{multline}
\widetilde{\mathsf{M}_{g}}((\alpha,\beta);\ell)=\sum_{\tilde{e}e=l}\frac{\varepsilon_q(\tilde{e})}{\sqrt{\tilde{e}}}\sum_{ab=\tilde{e}}\frac{\mu(a)\varepsilon_D(a)}{\sqrt{a}}\lambda_g(b) \\
\sum_{m,n\geq 1}\frac{\lambda_{g}(m)\lambda_{g}(n)}{\sqrt{mn}}V_{g,\alpha}\left(\frac{m}{qD}\right)V_{g,\beta}\left(\frac{a\tilde{e}n}{qD}\right)\Delta_q(m,aen) \\
+O_{\varepsilon, k, g}\left((q\ell)^\varepsilon(1+\vert\Im{(\mu)}\vert)^B\frac{\sqrt{\ell}}{q}\right)
\end{multline}
for some $B>0$ and for any $\varepsilon>0$.
\end{proposition}
As a consequence, if $1\leq\ell<q$ then \eqref{twist} is still valid even if there are some old forms. We will need the following easy lemma:
\begin{lemma}
\label{deltaq}
Let $N\geq 1$. For every integers $m, n\geq 1$, we have:
\begin{equation}
\Delta_{N}(m,n)\ll_{\varepsilon}(Nmn)^\varepsilon\frac{\sqrt{mn}}{N}.
\end{equation}
\end{lemma}
\noindent{\textbf{Proof of proposition \ref{newresult}.}}
The multiplicative properties of Hecke eigenvalues of $f$ and $g$ lead to:
\begin{multline}
\label{start0}
\widetilde{\mathsf{M}_{g}}((\alpha,\beta);\ell)=\sum_{\tilde{e}e=\ell}\frac{\varepsilon_q(\tilde{e})}{\sqrt{\tilde{e}}}\sum_{ab=\tilde{e}}\frac{\mu(a)\varepsilon_D(a)}{\sqrt{a}}\lambda_g(b) \\
\sum_{m,n\geq 1}\frac{\lambda_{g}(m)\lambda_{g}(n)}{\sqrt{mn}}V_{g,\alpha}\left(\frac{m}{qD}\right)V_{g,\beta}\left(\frac{a\tilde{e}n}{qD}\right)\sum_{f\in S_k^p(q)}^{\quad h}\lambda_f(m)\lambda_f(aen).
\end{multline}
We split the summation as follows:
\begin{equation}
\label{splitt}
\sum_{\substack{q\nmid m \\ q^2\nmid n}}\cdots+\sum_{\substack{q\nmid n \\ q^2\nmid m}}\cdots+\sum_{q^2\mid mn}\cdots:=I+II+III.
\end{equation}
The reader may check using mainly \eqref{proprioV} and \eqref{petersson} that $III\ll_{\varepsilon}q^{\varepsilon}\frac{\sqrt{ae}}{q}$. For the first term in \eqref{splitt} (the same analysis works for the second one), one can apply (\ref{peterssonnew}) which gives
\begin{multline}
\label{first1}
I=\sum_{\substack{q\nmid m \\ q^2\nmid n}}\frac{\lambda_g(m)\lambda_g(n)}{\nu(q\wedge n)\sqrt{mn}}V_{g,\alpha}\left(\frac{m}{qD}\right)V_{g,\beta}\left(\frac{a\tilde{e}n}{qD}\right)\Delta_q(m,aen) \\
-\frac{1}{q}\sum_{\tilde{q}\mid q^\infty}\frac{1}{\tilde{q}}\sum_{\substack{q\nmid m \\ q^2\nmid n}}\frac{\lambda_g(m)\lambda_g(n)}{\nu(q\wedge n)\sqrt{mn}}V_{g,\alpha}\left(\frac{m}{qD}\right)V_{g,\beta}\left(\frac{a\tilde{e}n}{qD}\right)\Delta_1(m\tilde{q}^2,aen):=Ia-Ib.
\end{multline}
Petersson trace formula \eqref{petersson} leads to
\begin{multline}
\label{first3}
Ib=\frac{1}{q}\sum_{\tilde{q}\mid q^\infty}\frac{1}{\tilde{q}}\sum_{h\in S_k^p(1)}\omega_1(h)\left(\sum_{\substack{m\geq 1 \\ q\nmid m}}\frac{\lambda_g(m)\lambda_h(\tilde{q}^2m)}{\sqrt{m}}V_{g,\alpha}\left(\frac{m}{qD}\right)\right) \\
\times\left(\sum_{\substack{n\geq 1 \\ q^2\nmid n}}\frac{\lambda_g(n)\lambda_h(aen)}{\sqrt{n}\nu(q\wedge n)}V_{g,\beta}\left(\frac{a\tilde{e}n}{qD}\right)\right) \\
:=\frac{1}{q}\sum_{\tilde{q}\mid q^\infty}\frac{1}{\tilde{q}}\sum_{h\in S_k^p(1)}\omega_1(h)Ib1\times Ib2.
\end{multline}
Let us study $Ib2$ (the same works for $Ib1$):
\begin{multline*}
Ib2=\sum_{\substack{n\geq 1 \\ q\nmid n}}\frac{\lambda_g(n)\lambda_h(aen)}{\sqrt{n}}V_{g,\beta}\left(\frac{a\tilde{e}n}{qD}\right) \\
+\frac{1}{\nu(q)}\sum_{\substack{n\geq 1 \\ q\mid\mid n}}\frac{\lambda_g(n)\lambda_h(aen)}{\sqrt{n}}V_{g,\beta}\left(\frac{a\tilde{e}n}{qD}\right):=Ib21+Ib22.
\end{multline*}
We limit ourselves to give an estimate of $Ib21$. Mellin's inversion formula entails that
\begin{equation}
\label{first4}
Ib21=\frac{1}{2i\pi}\int_{(2)}L\left(ae,q^2;\frac{1}{2}+z\right)\left(\frac{qD}{a\tilde{e}}\right)^z\widetilde{V_{g,\beta}}(z)\mathrm{d}z
\end{equation}
with $L(ae,q^2;z):=\sum_{\substack{n\geq 1 \\ q^2\nmid n}}\lambda_h(aen)\lambda_g(n)n^{-z}$ and
\begin{equation*}
\forall z\in\mathbb{C}, \quad \widetilde{V_{g,\beta}}(z):=\int_0^{+\infty}x^{z-1}V_{g,\beta}(z)\mathrm{d}z\ll(1+\vert t\vert)^B\vert z\vert^{-2}.
\end{equation*}
As usual, $L(ae,q^2;z)=R(ae;z)\frac{L^{(aeq^2)}(h\times g,z)}{\zeta^{(aeq^2)}(2z)}$ where $R(ae;z):=\sum_{n\mid(ae)^\infty}\frac{\lambda_h(aen)\lambda_g(n)}{n^z}$ converges on $\Re{(z)}>0$ and satisfies over there $R(ae;z)\ll\tau(ae)\ll_{\varepsilon}(ae)^\varepsilon$ for any $\varepsilon>0$. Shifting the $z$-contour to $\Re{(z)}=\varepsilon$ in \eqref{first4}, the convexity bound for $L(h\times g,.)$ implies that $Ib21\ll_{\varepsilon} (qae)^\varepsilon(1+\vert t\vert)^B$ for any $\varepsilon>0$. The same lines give $Ib22\ll_{\varepsilon} \frac{(qae)^\varepsilon}{\sqrt{q}\nu(q)}(1+\vert t\vert)^B$ for any $\varepsilon>0$. Finally, $Ib\ll_{\varepsilon, g}\frac{(aeq)^{\varepsilon}}{q}$ for any $\varepsilon>0$ and we have prove that:
\begin{multline}
\label{first5}
\widetilde{\mathsf{M}_{g}}(\alpha,\beta;l)=\sum_{\tilde{e}e=l}\frac{\varepsilon_q(\tilde{e})}{\sqrt{\tilde{e}}}\sum_{ab=\tilde{e}}\frac{\mu(a)\varepsilon_D(a)}{\sqrt{a}}\lambda_g(b) \\
\sum_{m,n\geq 1}\frac{\lambda_{g}(m)\lambda_{g}(n)}{\sqrt{mn}}V_{g,\alpha}\left(\frac{m}{qD}\right)V_{g,\beta}\left(\frac{a\tilde{e}n}{qD}\right)\Delta_q(m,aen) \\
-\sum_{\tilde{e}e=l}\frac{\varepsilon_q(\tilde{e})}{\sqrt{\tilde{e}}}\sum_{ab=\tilde{e}}\frac{\mu(a)\varepsilon_D(a)}{\sqrt{a}}\lambda_g(b) \\
\sum_{\substack{m,n\geq 1 \\ q^2\mid mn}}\frac{\lambda_{g}(m)\lambda_{g}(n)}{\sqrt{mn}}V_{g,\alpha}\left(\frac{m}{qD}\right)V_{g,\beta}\left(\frac{a\tilde{e}n}{qD}\right)\Delta_q(m,aen)+O_{\varepsilon, g}\left((q\ell)^\varepsilon(1+\vert t\vert)^B\frac{\sqrt{\ell}}{q}\right).
\end{multline}
The second term in \eqref{first5} is bounded by $\ll_{\varepsilon, g}(q\ell)^\varepsilon\frac{\sqrt{\ell}}{q}$ for any $\varepsilon>0$ thanks to lemma \ref{deltaq}.
\begin{flushright}
$\blacksquare$
\end{flushright}

\section{A review of Maass forms}

\label{Mas}
In this appendix, we only give the minimal knowledge about Maass forms in order to follow the notations which are used in this paper. The reader may see \cite{DuFrIw2} for all the details. Let $N\geq 1$ be a natural integer. A function $f:\mathbb{H}\rightarrow\mathbb{C}$ is said to be $\varGamma_0(N)$\emph{-automorphic} of weight $0$ and trivial nebentypus if it satisfies $f(\gamma.z)=f(z)$ for any $\gamma\in\varGamma_0(N)$. We denote by $\mathcal{L}_0(N)$ the space of square-integrable $\varGamma_0(N)$-automorphic functions with respect to the scalar product:
\begin{equation*}
(f,g):=\int_{\varGamma_0(N)\backslash\mathbb{H}}f(z)\overline{g(z)}y^{-2}\mathrm{d}x\mathrm{d}y.
\end{equation*}
The Laplacian $\Delta_0:=y^2\left(\frac{\partial^2}{\partial^2x}+\frac{\partial^2}{\partial^2y}\right)$ acts on $\mathcal{L}_0(N)$ and splits it in eigenspaces. There are two components: a discrete one spanned by the so-called \emph{Maass cusp forms} and a continuous one spanned by the \emph{Eisenstein series} which are given for any cusp $\kappa$ of $\varGamma_0(N)$ by
\begin{equation*}
E_\kappa(z,s):=\sum_{\gamma\in\left(\varGamma_0(N)\right)_\kappa\backslash\varGamma_0(N)}\left(\Im{(\sigma_\kappa^{-1}\gamma.z)}\right)^s
\end{equation*}
where $\sigma_\kappa$ is a scaling matrix for the cusp $\kappa$. The Eisenstein series is holomorphic on $Re{(s)}>1$, admits meromorphic continuation to $\mathbb{C}$ with only one  pole on $\Re{(s)}\geq\frac{1}{2}$ at $s=1$ and are eigenfunctions of the Laplacian: $\left(\Delta_0+\lambda(s)E_\kappa(.,s)\right)=0$ with $\lambda(s):=s(1-s)$ and $s=\frac{1}{2}+ir$ ($r\in\mathbb{C}$). They admit the following Fourier expansion
\begin{equation*}
E_\kappa\left(z,\frac{1}{2}+ir\right)=\delta_{\kappa=\infty}+\phi_\kappa\left(\frac{1}{2}+ir\right)y^{\frac{1}{2}-ir}+2\sqrt{y}\sum_{n\in\mathbb{Z}^*}\rho_\kappa(n,r)\vert n\vert^{\frac{1}{2}}K_{ir}(2\pi\vert n\vert y)e(nx)
\end{equation*}
with for any $n\in\mathbb{Z}^*$
\begin{equation}
\label{rhokappan}
\rho_\kappa(n,r)=\frac{\pi^s\vert n\vert^{ir-\frac{1}{2}}}{\varGamma\left(\frac{1}{2}+ir\right)}\left(\frac{\text{gcd}\left(w,\frac{N}{w}\right)}{wN}\right)^{\frac{1}{2}+ir}\sum_{\text{gcd}\left(\gamma,\frac{N}{w}\right)=1}\!\!\gamma^{-1-2ir}\sum_{\substack{\delta\mod (\gamma w) \\ \text{gcd}(\delta,\gamma w)=1 \\ \delta\gamma\equiv u\mod\left(w\wedge\frac{N}{w}\right)}}\!\! e\left(-n\frac{\delta}{\gamma w}\right)
\end{equation}
in the space $\Im{(r)}<0$ for $\kappa=\frac{u}{w}$ with $w\mid N$, $\text{gcd}(u,w)=1$, $1\leq u\leq\text{gcd}\left(w,\frac{N}{w}\right)$. Here, $\phi_\kappa\left(\frac{1}{2}+ir\right)$ is some explicit complex number. Let $\mathcal{E}_0(N)$ be the closure for $(.,.)$ in $\mathcal{L}_0(N)$ of the $\mathbb{C}$-vector space spanned by:
\begin{equation*}
\left\{\sum_{\gamma\in\left(\varGamma_0(N)\right)_\kappa\backslash\varGamma_0(N)}\psi\left(\Im{(\sigma_\kappa^{-1}\gamma.z)}\right), \psi\text{ compactly supported in }\mathbb{R}_+\right\}.
\end{equation*}
$\Delta_0$ has a continuous spectrum on $\mathcal{E}_0(N)$ which is $\left[\frac{1}{4},+\infty\right[$ and its multiplicity is the number of cusps of $\varGamma_0(N)$. Moreover, if $f$ belongs to $\mathcal{E}_0(N)$ then
\begin{equation*}
f(z)=(f,u_0)u_0(z)+\sum_{\kappa\in\text{Cusp}(\varGamma_0(N))}\frac{1}{4\pi}\int_{-\infty}^{+\infty}\left(f(.),E_\kappa\left(.,\frac{1}{2}+ir\right)\right)E_\kappa\left(z,\frac{1}{2}+ir\right)\mathrm{d}r
\end{equation*}
where $u_0$ is the constant function of value $\left(\text{Vol}(X_0(N))\right)^{-\frac{1}{2}}$. Let $\mathcal{C}_0(N)$ be the $(.,.)$-orthogonal of $\mathcal{E}_0(N)$: it is the space generated by the Maass cusp forms. The Fourier expansion of a Maass cusp form $f$ at infinity is
\begin{equation*}
f(z)=2\sqrt{y}\sum_{n\in\mathbb{Z}^*}\rho_f(n)\vert n\vert^{\frac{1}{2}}K_{ir_f}(2\pi\vert n\vert y)e(nx)
\end{equation*}
where $\left(\Delta_0+\lambda_f\right)f=0$ and $\lambda_f:=\lambda(s_f):=\lambda\left(\frac{1}{2}+ir_f\right)$. Let $\left(u_j\right)_{j\geq 1}$ be an orthonormal basis of $\mathcal{C}_0(N)$ made of Maass cusp forms. If $f$ belongs to $\mathcal{C}_0(N)$ then
\begin{equation*}
f(z)=\sum_{j\geq 1}(f,u_j)\, u_j(z).
\end{equation*}
J.-M. Deshouillers and H. Iwaniec established in \cite{DeIw} the following large sieve inequalities for all the previous Fourier coefficients
\begin{equation}
\label{sievedisc}
\sum_{\vert r_j\vert\leq R}\frac{1}{\cosh{(\pi r_j)}}\left\vert\sum_{1\leq m\leq M}a_m m^{\frac{1}{2}}\rho_{j}(m)\right\vert^2\ll_{\varepsilon}\left(R^2+\frac{M^{1+\varepsilon}}{N}\right)\vert\vert a\vert\vert_2^2,
\end{equation}
\begin{multline}
\label{sievecont}
\sum_{\kappa\in \text{Cusp}(\varGamma_0(N))}\int_{-R}^{+R}\left\vert\varGamma\left(\frac{1}{2}+ir\right)\right\vert^2\left\vert\sum_{1\leq m\leq M}a_m m^{\frac{1}{2}}\rho_{\kappa}(m,r)\right\vert^2\mathrm{d}r\ll_{\varepsilon} \\
\left(R^2+\frac{M^{1+\varepsilon}}{N}\right)\vert\vert a\vert\vert_2^2
\end{multline}
for $R\geq 1$, any $\varepsilon>0$ and any sequence of complex numbers $\left(a_m\right)_{1\leq m\leq M}$. The Hecke operators $\left(T_n\right)_{n\geq 1}$ also act on $\mathcal{L}_0(N)$, commute with $\Delta_0$ and are hermitian if $\text{gcd}(n,N)=1$. A \emph{Hecke-Maass} cusp form is a Maass cusp form which is also an eigenfunction of the $T_n$ for $\text{gcd}(n,N)=1$. A \emph{Hecke eigenbasis} is an orthonormal basis of $\mathcal{C}_0(N)$ made of Hecke-Maass cusp forms. For $f$ a Hecke-Maass cups form of Hecke eigenvalues $\left(\lambda_f(n)\right)_{\text{gcd}(n,N)=1}$, one has for any $\text{gcd}(mn,N)=1$:
\begin{equation}
\lambda_f(m)\lambda_f(n)=\sum_{d\mid m\wedge n}\varepsilon_N(d)\lambda(mnd^{-2}),
\end{equation}
\begin{equation}
\lambda_f(mn)=\sum_{d\mid m\wedge n}\mu(d)\varepsilon_N(d)\lambda_f(m/d)\lambda_f(n/d).
\end{equation}
The action of Hecke operators on the Fourier expansion of a Hecke-Maass cusp form $f$ is known:
\begin{equation}
\sqrt{m}\rho_f(m)\lambda_f(n)=\sum_{d\mid m\wedge n}\varepsilon_N(d)\rho_f(mnd^{-2})\sqrt{\frac{mn}{d^2}},
\end{equation}
\begin{equation}
\label{Maassmulti}
\sqrt{mn}\rho_f(mn)=\sum_{d\mid m\wedge n}\mu(d)\varepsilon_N(d)\rho_f(m/d)\sqrt{\frac{m}{d}}\lambda_f(n/d)
\end{equation}
for any $m, n\geq 1$ with $\text{gcd}(n,N)=1$.

\section{The computation of an Euler product}

The purpose of this appendix is to prove that the arithmetical constants which appear in the asymptotic formulas of the harmonic mollified second moment equal one. More precisely, we prove that equations \eqref{calcul1} and \eqref{calcul2} hold.\newline\newline
Remember that, according to lemma \ref{lemmaII}, $h_2(\mu+\overline{\mu},\mu,\overline{\mu})$ is an absolutely convergent Euler product when the real part of $\mu$ is greater than a small negative real number (say $10^{-6}$) namely $h_{2}(\mu+\overline{\mu},\mu,\overline{\mu})=\prod_{p\in\mathcal{P}}h_{2,p}(\mu+\overline{\mu},\mu,\overline{\mu})$ with (confer \eqref{h2bis}):
\begin{multline}
\label{Eulerapp}
\forall p\in\mathcal{P}, h_{2,p}(\mu+\overline{\mu},\mu,\overline{\mu})=\prod_{z\in\left\{\mu,\overline{\mu}\right\}}\left(\frac{L_p(g\times g,1+z+\mu)L_p(g\times g,1+z+\overline{\mu})}{L_p^{(q)}(\text{Sym}^2(g),1+2z)}\right) \\
L_p(g\times g,1+\mu+\overline{\mu}){L}_{p}(\mu+\overline{\mu},\mu,\overline{\mu})
\end{multline}
where $L_{p}(\mu+\overline{\mu},\mu,\overline{\mu})$ is defined in \eqref{tildeL}. Firstly, we need to have an idea of the shape of $\nu_g(p^k;u,v)$ (see \eqref{nnnnug} for its definition) for any prime number $p$, any natural integer $1\leq k\leq 4$ and any complex numbers $u$ and $v$.
\begin{lemma}
\label{vug}
Let $p$ be a prime number, $1\leq k\leq 4$ be some natural integer and $u$, $v$ be some complex numbers. If $Q:=p^{-1}$, $U:=p^{-u}$ and $V:=p^{-v}$ then it turns out that $\nu_g(p^k;u,v)$ is of the following shape:
\begin{equation*}
\nu_g(p^k;u,v)=(1+\varepsilon_D(p)QUV)^{-1}P_{g,k}(Q,U,V)
\end{equation*}
where $P_{g,k}(Q,U,V)$ is some explicit polynomial in three variables whose coefficients depend on $\varepsilon_D(q)$ and on $\lambda_g(p)^i$ ($1\leq i\leq k$).
\end{lemma}
\noindent{\textbf{Proof of lemma \ref{vug}.}}
We set for any natural integers $a,b\geq 0$:
\begin{equation*}
S_g(a,b;u,v):=\sum_{k\geq 0}\lambda_g(p^{k+a})\lambda_g(p^{k+b})Q^kU^kV^k.
\end{equation*}
The relationships beetween Hecke eigenvalues of $g$ enable us to express $S_g(a,b;u,v)$ in function of $S_g(0,0;u,v)$ by induction. More precisely, it shows that:
\begin{equation*}
S_g(a,b;u,v)=(1+\varepsilon_D(p)QUV)^{-1}R_{g,a,b}(Q,U,V)
\end{equation*}
for some explicit polynomial in three variables $R_{g,a,b}(Q,U,V)$ whose coefficients depend on $\varepsilon_D(q)$ and on $\lambda_g(p)^i$ ($1\leq i\leq a+b$). In addition, one remarks that:
\begin{eqnarray*}
S_g(0,0;u,v)\nu_g(p^1;u,v) & = & (U+V)S_g(1,0;u,v), \\
S_g(0,0;u,v)\nu_g(p^2;u,v) & = & (U^2+V^2)S_g(2,0;u,v)+UVS_g(1,1;u,v), \\
S_g(0,0;u,v)\nu_g(p^3;u,v) & = & (U^3+V^3)S_g(3,0;u,v)+(U^2V+UV^2)S_g(2,1;u,v), \\
S_g(0,0;u,v)\nu_g(p^4;u,v) & = & (U^4+V^4)S_g(4,0;u,v)+(U^3V+UV^3)S_g(3,1;u,v) \\
& & +U^2V^2S_g(2,2;u,v).
\end{eqnarray*}
Both previous remarks lead to the result.
\begin{flushright}
$\blacksquare$
\end{flushright}
\begin{remark}
The proof of lemma \ref{vug} also gives the explicit procedure we used for computing $\nu_g(p^k;u,v)$ for any prime number $p$, any natural integer $1\leq k\leq 4$ and any complex numbers $u$ and $v$.
\end{remark}
Having this in mind, we can compute the local factor $h_{2,p}(\mu+\overline{\mu},\mu,\overline{\mu})$ at each prime $p$ which does not divide $q$:
\begin{lemma}
\label{partialEuler}
Let $\mu$ be a complex number. We have:
\begin{equation*}
\forall p\nmid q,\quad \frac{h_{2,p}(\mu+\overline{\mu},\mu,\overline{\mu})}{\zeta_p^{(D)}(2(1+\mu+\overline{\mu}))}=1.
\end{equation*}
\end{lemma}
\noindent{\textbf{Proof of lemma \ref{partialEuler}.}}
Once again, we set $Q:=p^{-1}$, $U:=p^{-\mu}$ and $V:=p^{-\overline{\mu}}$. With these notations and knowing the local factors of Rankin-Selberg $L$-functions and symmetric-square $L$-functions, one computes that
\begin{multline}
\label{111}
\prod_{z\in\left\{\mu,\overline{\mu}\right\}}\left(\frac{L_p(g\times g,1+z+\mu)L_p(g\times g,1+z+\overline{\mu})}{L_p^{(q)}(\text{Sym}^2(g),1+2z)}\right)\frac{L_p(g\times g,1+\mu+\overline{\mu})}{\zeta_p^{(D)}(2(1+\mu+\overline{\mu}))}= \\
(1-\lambda_g(p)^2QUV)^{-1}
\end{multline}
if $p\mid D$ and $p\nmid q$ and that
\begin{multline}
\label{222}
\prod_{z\in\left\{\mu,\overline{\mu}\right\}}\left(\frac{L_p(g\times g,1+z+\mu)L_p(g\times g,1+z+\overline{\mu})}{L_p^{(q)}(\text{Sym}^2(g),1+2z)}\right)\frac{L_p(g\times g,1+\mu+\overline{\mu})}{\zeta_p^{(D)}(2(1+\mu+\overline{\mu}))}= \\
\frac{1+QUV}{(1-QU^2)(1-QV^2)(1-QUV)(1+2QUV-\lambda_g(p)^2QUV+Q^2U^2V^2)}
\end{multline}
if $p\nmid (Dq)$. According to lemma \ref{vug} and its definition, $L_p(\mu+\overline{\mu},\mu,\overline{\mu})$ for any prime number $p$ is, a priori, a rational fraction in three variables $Q$, $U$ and $V$; namely it looks like
\begin{equation*}
L_p(\mu+\overline{\mu},\mu,\overline{\mu})=\frac{P_1(Q,U,V)}{P_2(Q,U,V)}
\end{equation*}
for some polynomials $P_1$ and $P_2$ of total degrees less than $20$ whose coefficients depend on $\varepsilon_q(p)$, $\varepsilon_D(p)$ and on the Hecke eigenvalues of $g$ at the powers of $p$ namely $\lambda_g(p)^k$ for $1\leq k\leq 4$. To factor this fraction we have used a computational algebra system; for instance the scripts of this computation (vg.mws, ctemumubar.mws and cte.mws) are available at\newline
http://www.dms.umontreal.ca/$\sim$ricotta. We 
obtain
\begin{equation}
\label{aaa}
L_p(\mu+\overline{\mu},\mu,\overline{\mu})=1-\lambda_g(p)^2QUV
\end{equation}
if $p\mid D$ and $p\nmid q$ and
\begin{equation}
\label{bbb}
L_p(\mu+\overline{\mu},\mu,\overline{\mu})=\frac{(1-QU^2)(1-QV^2)(1-QUV)(1+2QUV-\lambda_g(p)^2QUV+Q^2U^2V^2)}{1+QUV}
\end{equation}
if $p\nmid Dq$. Note that the computations above are purely formal (no numerical approximation is made); in fact, once the above factorizations have been obtained, it is possible (but lenghtly) to check them directly by hand. Then we finish the proof of lemma D.2 by simplifying \eqref{111} with \eqref{aaa} and \eqref{222} with \eqref{bbb}.
\begin{flushright}
$\blacksquare$
\end{flushright}
We can now state the main result:
\begin{proposition}
\label{arithconst}
Let $\mu$ be a complex number with $\tau:=\Re{(\mu)}\geq -\gamma$ for some $\gamma>0$ small enough. We have:
\begin{equation*}
\frac{h_{2}(\mu+\overline{\mu},\mu,\overline{\mu})}{\zeta^{(D)}(2(1+\mu+\overline{\mu}))}=1+O_g\left(\frac{1}{q^\delta}\right)
\end{equation*}
for some $\delta>0$.
\end{proposition}
\noindent{\textbf{Proof of proposition \ref{arithconst}.}}
The proof is an immediate consequence of the previous lemma as the various Euler products are absolutely convergent under the assumption made on $\mu$. The admissible error term comes from the local factor of the Euler product at the primes which divide $q$.
\begin{flushright}
$\blacksquare$
\end{flushright}

\strut\newline
\noindent{\textit{G. Ricotta}\newline}
Université Montpellier II, Département de Mathématiques, Case Courier 051, 34095 Montpellier CEDEX 05, France; ricotta@math.univ-montp2.fr; current: Université de Montréal, Département de Mathématiques et de Statistique, CP 6128 succ Centre-Ville, Montréal QC   H3C 3J7, Canada; ricotta@dms.umontreal.ca.
 
\end{document}